\documentclass[11pt]{article}
\usepackage{amsmath,amssymb,amsthm,amsfonts,authblk,biblatex,graphicx,enumerate,verbatim,appendix,subcaption,overpic}
\usepackage[labelfont=bf]{caption}
\usepackage[margin=1in]{geometry}
\usepackage{xcolor}
\usepackage{hyperref}
\renewbibmacro{in:}{}

\newcommand{\abs}[1]{\left\lvert #1 \right\rvert}
\newcommand{\norm}[1]{\left\lVert #1 \right\rVert}

\newcommand{\hfrac}[2]{#1/#2}

\def\R{\mathbb{R}}
\def\O{\mathcal{O}}

\def\Tr{\text{Tr}}

\def\B{\mathcal{B}}

\def\E{\mathbb{E}}
\def\G{\mathbb{G}}
\def\eps{\varepsilon}
\newcommand{\vecspan}[1]{\text{span}\left\{ #1 \right\}}

\numberwithin{equation}{section}

\usepackage{multicol}
\usepackage{fancyhdr}

\usepackage{comment}

\newtheorem{lemma}{Lemma}[section]
\newtheorem{remark}{Remark}[section]
\newtheorem{proposition}{Proposition}[section]

\title{Stable and unstable spatially-periodic canards created in singular subcritical Turing bifurcations in the Brusselator system}
\author[1]{Robert Jencks}
\author[2]{Arjen Doelman}
\author[1]{Tasso J. Kaper}
\author[3]{Theodore Vo}

\affil[1]{\footnotesize Department of Mathematics and Statistics, Boston University, Boston, MA 02215, USA}
\affil[2]{\footnotesize Mathematisch Instituut, Universiteit Leiden, 2300 RA Leiden, the Netherlands}
\affil[3]{\footnotesize School of Mathematics, Monash University, Clayton, Victoria 3800, Australia}
\date{\today}

\begin{document}
\maketitle

\begin{abstract}
In this article, we study the canonical Brusselator partial differential equation (PDE) from the field of pattern formation in the limit in which the diffusivity of the activator is much smaller than that of the inhibitor. 
The PDE robustly exhibits a subcritical Turing bifurcation that, in this limit, is labeled as a {\em singular Turing bifurcation}. 
We show that families of spatially-periodic canard solutions emerge from this subcritical singular Turing bifurcation. 
Then, right after they emerge, the solutions lose their purely sinusoidal structure ($e^{ik_Tx}$, where $k_T$ is the critical wavenumber at the Turing bifurcation) and gain a distinct multi-scale spatial structure.
They consist of segments along which the components vary gradually in space, interspersed with short intervals on which the activator component exhibits pulses and steep gradients.
The branches of these spatially-periodic canards undergo a saddle-node bifurcation.
Some of the large-amplitude patterns on the upper branches are attractors of the PDE, and unstable patterns with small pulses that exist below the folds appear to guide the evolution of data to the attractors.
We also analyze the spatial ordinary differential equations (ODEs) that govern time-independent solutions.
We show that the spatial ODE system has a folded singularity known as a reversible folded saddle-node of type II (RFSN-II) point that coincides with the Turing bifurcation in the singular limit.
We demonstrate that, for parameter values close to the Turing bifurcation, the true and faux canards of the RFSN-II point are responsible for generating the spatially-periodic canards, and for parameter values away from the Turing value there is a reversible folded saddle point whose true and faux canards generate the spatially-periodic canard solutions. 
In short, these solutions are new examples of ``les canards de Turing", first discovered in the van der Pol PDE in \cite{VDK2025}.
Overall, we identify the RFSN-II and RFS folded singularities and their canards as new selection mechanisms for subcritical bifurcations.
\end{abstract}

\vspace{1em}


\section{Introduction}
In the prototypical van der Pol partial differential equation (PDE) with widely separated diffusivities for the activator and inhibitor, new families of Turing patterns were recently discovered that consist of spatially-periodic canard solutions \cite{VDK2025}.
Just as for classical Turing bifurcations \cite{CH1993,E1965,EK2005,T1952,W1997}, these patterns emerge as stationary, small-amplitude, spatially-periodic solutions from a  homogeneous state.
However, they have canard segments in space ({\it i.e.,} long segments on which they are near unstable states), which had not previously been observed in Turing patterns.
Moreover, these patterns continue (in the bifurcation parameter) to families of large-amplitude spatially-periodic canards that have distinct fast-slow structure with spike layers and that are stationary attractors of the PDE.

These families of small- and large-amplitude canard patterns were discovered based on the observation --also reported in \cite{VDK2025}-- that the Turing bifurcation in the van der Pol PDE corresponds to the limit of a folded singularity with reversibility symmetry known as a reversible folded saddle-node of type II (RFSN-II) 
in the associated system of spatial ordinary differential equations (ODEs),
see \cite{SW2001} for folded singularities and general FS points and \cite{Krupa2010} for general FSN-II points.
We showed that the RFSN-II points and their true and faux canards are the mechanisms responsible for generating the canard segments of the spatially-periodic patterns in the PDE.

In \cite{VDK2025}, it was conjectured that the RFSN-II points, this selection mechanism, and families of spatially-periodic canard patterns exist not just in the van der Pol PDE but more generally in coupled systems of reaction-diffusion equations of activator-inhibitor type in which the diffusivities are widely separated.
The goal of this article is to present a second prototypical PDE system that possesses these folded singularities as well as classes of small-amplitude and large-amplitude spatially-periodic canard patterns.
Furthermore, we have chosen the second example to have monostable kinetics, to show that the bistability of the van der Pol kinetics is not required, {\it i.e.,} that the spatial canard solutions also exist in systems with monostable kinetics.

In this article, we analyze the Brusselator PDE,
\begin{equation}  \label{Brusselator}
\begin{split}
    \partial_{t}u &= d\Delta_{x} u-(B+1)u+A+u^{2}v, \\
    \partial_{t}v &= \Delta_{x} v+Bu-u^{2}v.
\end{split}
\end{equation} 
Introduced in \cite{PL1968}, the Brusselator PDE is a prototypical system in pattern formation, see \cite{EP1998,GS1997,Pena2001,SZJV2018,VWDB1992,W1997} for example.
Here, $x \in$ \(\R\), $t \ge 0$, $u=u(x,t)$ and $v=v(x,t)$ denote the concentrations of the activator and inhibitor species, 
and \(A\) and \(B\) are positive rate constants. The parameter \(d\) is a diffusion coefficient, which we consider as a small parameter: $0< d \ll 1$, so that the diffusivities of the activator and inhibitor are widely separated. 

The Brusselator PDE has a unique spatially homogeneous equilibrium state, 
\begin{equation}
\label{equilib}
(u,v)=\left(A,\frac{B}{A}\right).
\end{equation}
For all $A>0$ and $d>0$, this state undergoes a Turing bifurcation at $B_T=(1+\sqrt{d}A)^2$ (see \cite{Pena2001,SZJV2018,VWDB1992}, and the derivation is presented below in Section \ref{ss:linstab}).
In this article, we focus on the case $0<d\ll 1$, which implies that the Turing bifurcation is singular and subcritical (see Section \ref{ss:weaklynonlinear}).
We show that the spatially-periodic solutions emerging from the singular Turing bifurcation are spatially-periodic canard patterns and that a folded singularity asymptotically close to the Turing bifurcation point is responsible for generating the canards.

\begin{figure}[ht!]
    \centering
    \includegraphics[width=5in]{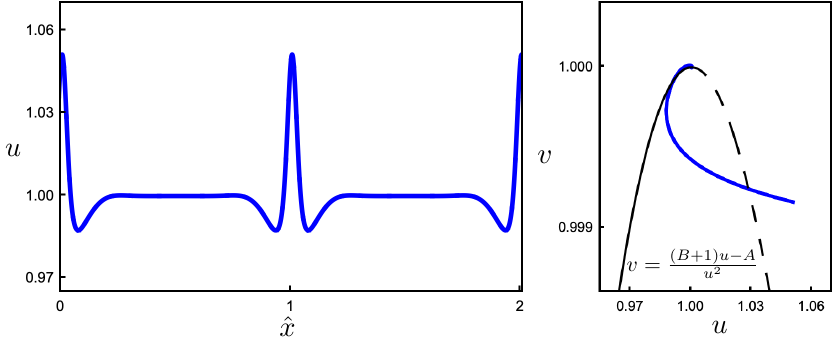}
  \put(-364,140){(a)}
  \put(-130,140){(b)}
    \caption{
    A representative stationary, spatially-periodic solution (blue) of the Brusselator PDE \eqref{Brusselator}. 
    It has wavenumber $k=2$ with period $T=2\pi/k\sqrt{d} \approx 314.829$, and it exhibits one small-amplitude pulse each period. 
    (a) Spatial profile of the $u$-component as a function of the scaled spatial variable $\hat{x} = x/T$. 
    (b) Projection of the solution into the $(u,v)$ phase plane. 
    Here, $A=1$, $B=1$, and $\sqrt{d}=0.01$.
    }
    \label{fig:representative-SAO}
\end{figure}

\begin{figure}[h!tbp]
    \centering
    \includegraphics[width=5in]{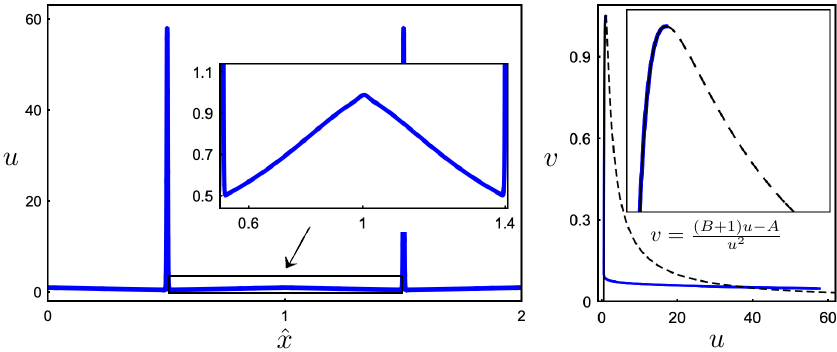}
    \put(-366,142){(a)}
    \put(-130,142){(b)}
    \caption{A second representative stationary, spatially-periodic solution (blue) of \eqref{Brusselator}.
    It has a large period ($T=665$) and a pronounced fast-slow structure.
    (a) $u$ as a function of the scaled spatial variable $\hat{x} = x/T$. 
    (b) Projection of the solution into the $(u,v)$ plane. 
    Here, $A=1$, $B=1.05$, and $\sqrt{d}=0.01$.}
    \label{fig:representative}
\end{figure}

Examples of spatially-periodic canard patterns of \eqref{Brusselator} are presented in Figs.~\ref{fig:representative-SAO} and \ref{fig:representative}.
The pattern in Fig.~\ref{fig:representative-SAO} has wavenumber $k=2$ and small, $\mathcal{O}(\sqrt{d})$ amplitude.
It does not have a purely sinusoidal profile as would be expected from spatial patterns that emerge from a Turing bifurcation. 
Instead, it has fast-slow structure.
The fast segment is a short spatial interval where $u$ exhibits a pulse (Fig.~\ref{fig:representative-SAO}(a)).
In phase space, this pulse corresponds to a rapid excursion that $u$ makes around the right branch of the $u$-nullcline (dashed black curve in Fig.~\ref{fig:representative-SAO}(b)).
The slow segment consists of the long, complementary portion of the pattern.
On it, $u$ is close to equilibrium $u=A$, and near the left branch (solid black curve) of the $u$-nullcline in the phase space.

The second sample spatially-periodic canard orbit, shown in Fig.~\ref{fig:representative}(a), also has fast-slow structure.
However, the pulses in $u$ have much a larger amplitude.
In phase space (Fig.~\ref{fig:representative}(b)), $u$ makes a large-amplitude excursion, with a steep spatial gradient, around the right branch of the $u$-nullcline (dashed black curve), where $v$ is small and constant to leading order. 
The slow segments are much longer than those in Fig.~\ref{fig:representative-SAO}.
Indeed, the orbit has spatial period $T=665$, and the wavenumber $k \approx 0.945$ is well below $k_T$ (where $k_T$, with $k_T^2 = A/\sqrt{d}$, is the critical wavenumber associated to the Turing instability, so that $k_T = 10$ here).
On the slow segments, $u$ gradually increases and then decreases (inset in Fig.~\ref{fig:representative}(a)), which in the phase space corresponds to the long, slow, upward and downward drift close to the left branch of the $u$-nullcline, from the small value of $v$ corresponding to the jump value all the way up to the local maximum (knee) of the nullcline, and then back.

In this article, the overall goal is to study the geometry and quantitative properties of several types of spatially-periodic canard solutions of the Brusselator PDE \eqref{Brusselator}, of the kind shown in Figs.~\ref{fig:representative-SAO} and \ref{fig:representative}, as well as other types.
First, we show that, for $B=1$, which is asymptotically close to the Turing bifurcation $B_T$, the spatial ODEs governing the time-independent solutions of \eqref{Brusselator} have a folded singularity known as a reversible folded saddle-node of type II (RFSN-II). 
Then, for each $B>1 + \mathcal{O}(\sqrt{d})$, the system has a reversible folded saddle (RFS). 
We study the RFSN-II and RFS points using geometric desingularization \cite{DR1996}.
Second, we construct families of spatially-periodic canard solutions using geometric singular perturbation theory  \cite{F1979,J1995}.
We show that, for $B=1$ and for each $B$ close to one, the canard patterns are created by the true and faux canards of the RFSN-II point.
Then, for each $B>1+\mathcal{O}(\sqrt{d})$, they are created instead by the RFS true and faux canards. 
Third, we analyze the near self-similarity of the canards.
Fourth, the stability of these spatially-periodic canards in the PDE is studied numerically, both on periodic domains and on large domains with homogeneous Neumann boundary conditions to study side-band instabilities.
The data evolve to stable canard patterns with large-amplitude pulses in the Busse balloon \cite{B1978}, and the transients consist of unstable canard patterns with small-amplitude spikes.

The new canard patterns in \eqref{Brusselator} are created in subcritical Turing bifurcations, as in the van der Pol PDE \cite{VDK2025}.
Also, there is a fold of spatially-periodic canard patterns, along which branches of unstable and stable patterns meet.
The analysis here (and in \cite{VDK2025}) reveals that folded singularities and their associated canards serve as a new pattern selection mechanism in the subcritical case.

Our results complement the literature on sub-critical Turing bifurcations generating spatially localized patterns (see \cite{CASBMGMVW2021, GZK2018, K2015}) that typically also exhibit saddle-node bifurcations of periodic patterns (see \cite{CH1993, S2003}). 
These periodic patterns have the sinusoidal nature of `classical' Turing patterns. 
Then, to analytically establish the saddle-node bifurcation, one assumes that the parameters are (relatively) close to the co-dimension 2 point at which the bifurcation changes to supercritical --which is not the case for the canard patterns here.  
On the other hand, unlike in the studies \cite{CASBMGMVW2021, GZK2018, K2015}, it is crucial here that the reaction-diffusion equation is singularly perturbed -- and thus that the Turing bifurcation is singular. 
Therefore, together with \cite{VDK2025}, our analysis provides a complementary approach to understanding the dynamics of patterns near subcritical Turing bifurcations.

We choose the Brusselator as the second example of ``les Canards de Turing", because the Turing and $k=0$ temporal Hopf bifurcations are well separated, by $\mathcal{O}(1)$ in the main parameter (cf. van der Pol PDE \cite{VDK2025}, where they are asymptotically close).
Hence, the Turing bifurcation can be studied more fully by itself. 
Furthermore, the results extend the existence of the families of spatially-periodic canards to the general case in which the spatial ODE system is fully 4-dimensional and does not have a conserved quantity, where in contrast the van der Pol has one.

The article is organized as follows.
In Sec.~\ref{s:TH}, the Turing and Hopf bifurcations and the Ginzburg-Landau amplitude equation of \eqref{Brusselator} are reviewed. 
Also, we apply the normal form theory for reversible, 1:1 resonant (spatial) Hopf bifurcation points in the system of spatial ODEs that govern time-independent solutions of \eqref{Brusselator}. 
In Sec.~\ref{s:3}, we present some of the main spatially-periodic canard patterns of the spatial ODE system of \eqref{Brusselator}.
Next, in Secs.~\ref{s:spatialODE}--\ref{s:K1}, we analyze the spatial ODE system and identify the key features responsible for generating the spatial canards.
In particular, in Sec.~\ref{s:spatialODE}, we analyze the dynamics of the two-dimensional fast and slow systems of the spatial ODE, identifying the critical manifold, the key RFSN-II point, and its true and faux canards.
Then, Secs.~\ref{s:desingularization} and \ref{s:K1} contain the geometric desingularization analysis of the central RFSN-II point and the analysis of the dynamics and invariant manifolds in the rescaling chart and the entry/exit chart, respectively.
In Sec.~\ref{s:geoconstruction}, we use
the desingularization  results to present the geometric construction of families of small-amplitude and large-amplitude spatially-periodic canard patterns in the spatial ODE.
In Sec.~\ref{s:selfsimilarity}, we focus on the nearly self-similar dynamics exhibited by some classes of spatially-periodic canards.
Secs.~\ref{s:PDE-numerics} and \ref{s:PDE-largedomains} contain the results of the simulations of \eqref{Brusselator} on periodic domains and large domains, respectively.
Conclusions, discussion comparing and contrasting the spatial canards with classical temporal canards in fast-slow systems, and some open problems about spatial canards are presented in Sec.~\ref{sec:conclusions}.
The appendices contain several useful results.

\section{Brief Review: Turing and Hopf Bifurcations, the Ginzburg-Landau Equation, and the Normal Form of the Spatial ODE} \label{s:TH}

\subsection{Linear stability analysis}
\label{ss:linstab}

Linearization of the system~\eqref{Brusselator} around the homogeneous equilibrium~\eqref{equilib} yields the operator
\begin{equation*}
    L=\begin{bmatrix}
        d \Delta_{x}+B-1 & A^{2}\\
        -B & \Delta_{x}-A^{2}
    \end{bmatrix}.
\end{equation*}
After a Fourier transform, we find
\begin{equation*}
    \hat{L}(k)=\begin{bmatrix}
        B-1-dk^{2} & A^{2}\\
        -B & -A^{2}-k^{2}
    \end{bmatrix}.
\end{equation*}
Hence, the trace and determinant of $\hat{L}(k)$ are given by
\begin{align*}
    \Tr \,\hat{L}(k)=-(1+d)k^{2}+B-A^{2}-1, \qquad 
    \det \hat{L}(k)=dk^{4}+(-B+dA^{2}+1)k^{2}+A^{2}.
\end{align*}
Therefore, the Turing bifurcation occurs at 
\begin{equation} \label{kTBT}
    B_{T}=\left( 1+\sqrt{d}A \right)^2, \qquad 
    k_{T}^{2}=\frac{A}{\sqrt{d}}, 
\end{equation} 
where the determinant and its partial derivative with respect to $k^2$ both vanish.
At the Turing point, the homogeneous state is marginally unstable to spatially-periodic perturbations of the form $e^{i k_T x}$, \cite{E1965,EK2005,S2003,T1952,W1997}.
(We note that the determinant also vanishes at \(B=(1-\sqrt{d}A)^2\), however there \(k^{2}<0\). 
This corresponds to a Belyakov-Devaney point
\cite{B1974,D1977}.) 

From the matrix $\hat{L}(k)$, we also see that the homogeneous steady state \eqref{equilib}  undergoes a ($k=0$) temporal Hopf bifurcation at 
\begin{align}
\label{BH}
    B_{H}&=A^{2}+1.
\end{align}
For the PDE problem on the real line, which is the main problem we study, a long wave instability occurs for $B>B_H$.
We add that, on finite domains, additional modes
are unstable for $B>B_H(k)=A^2 + 1 + k^2 (1 + d)$.

\subsection{Weakly nonlinear stability}  \label{ss:weaklynonlinear}

The Ginzburg-Landau amplitude equation for the Brusselator~\eqref{Brusselator} is
\begin{equation} \label{GL}
       \left(\tfrac{1+A\sqrt{d}}{1-d}\right)\partial_{\tau}M_{1}=\tfrac{4}{(1+A\sqrt{d})^{2}}\partial_{\chi}^{2}M_{1}+\left(\tfrac{B-B_{T}}{B_{T}}\right) M_{1} +\left(\tfrac{8-38A\sqrt{d}-5A^{2}d+8A^{3}d\sqrt{d}}{9A^{3}\sqrt{d}(1+A\sqrt{d})}\right) \abs{M_{1}}^{2}M_{1},
\end{equation}
where we have used
\begin{align} 
        \begin{bmatrix}
        u(x,t)\\
        v(x,t)
    \end{bmatrix}=&\begin{bmatrix}
        A\\
        \tfrac{B}{A}
    \end{bmatrix}+\delta M_{1}(\delta^{2}t,\delta x)e^{ik_{T}x}\begin{bmatrix}
        1\\
        -\tfrac{\sqrt{d}}{A}(1+A\sqrt{d})
    \end{bmatrix}+c.c.+\O(\delta^{2}). \nonumber
\end{align}
Here, $c.c.$ denotes complex conjugate, $0<\delta \ll 1$, 
\(\tau=\delta^{2}t\), and \(\chi=\delta x\). 
This follows from known results, see for example \cite{Pena2001,SZJV2018,VWDB1992}, and we refer to \cite{E1965,SU2017} for Ginzburg-Landau theory.

Since $0<d \ll 1$, the coefficient on the time derivative in the left member of \eqref{GL} is positive. 
The coefficient on the cubic term in \eqref{GL}, which is often referred to as the Landau coefficient, is a cubic function of $A$.
It vanishes at $A=\frac{21 \pm \sqrt{313}}{16\sqrt{d}}$ and at $A= -2/\sqrt{d}$.
Hence, the Landau coefficient is positive for $A \in \left( 0, \frac{21-\sqrt{313}}{16\sqrt{d}}\right)$ and $A \in \left(\frac{21+\sqrt{313}}{16\sqrt{d}}, \infty \right)$,
and it is negative for $A \in \left( \frac{21-\sqrt{313}}{16\sqrt{d}}, \frac{21 +\sqrt{313}}{16\sqrt{d}}\right)$.
We focus on the first interval of $A$ values,
since $0<d\ll 1$, and the Turing bifurcation is subcritical.

\begin{remark} With \(\sqrt{d}\) small the Ginzburg-Landau equation is given to leading order by
\begin{equation}
    \partial_{\tau}M_{1}=4\partial_{\chi}^{2}M_{1}+\left(\tfrac{B-B_{T}}{B_{T}}\right) M_{1} + \tfrac{8}{9A^{3}\sqrt{d}} \abs{M_{1}}^{2}M_{1},
\end{equation}
where 
\begin{align}
    \begin{bmatrix}
        u(x,t)\\
        v(x,t)
    \end{bmatrix}=&\begin{bmatrix}
        A\\
        \tfrac{B}{A}
    \end{bmatrix}+\delta M_{1}(\delta^{2}t,\delta x)e^{ik_{T}x}\begin{bmatrix}
        1\\
        0
    \end{bmatrix}+c.c.+\O(\delta^{2}). \nonumber
\end{align}
\end{remark}

The criticality of the Turing bifurcation may also be determined by examining the system of ODEs governing the spatial dynamics of steady solutions,
\begin{equation}  \label{spatialODE-x}
\begin{split}
    \eps u_x &= p \\
    \eps p_x &= -A + (B+1) u - u^2 v \\
    v_x &= q \\
    q_x &= -Bu + u^2 v.
\end{split}
\end{equation}
Here, $\eps=\sqrt{d}$.
For all $\eps>0$, an equivalent formulation of this system is given by
\begin{equation}  \label{spatialODE}
\begin{split}
    u_y &= p \\
    p_y &= -A + (B+1) u - u^2 v \\
    v_y &= \eps q \\
    q_y &= \eps (-Bu + u^2 v),
\end{split}
\end{equation}
where $y=x/\eps$.
We use both formulations at different points in the analysis. 
The former will be useful for studying the `slow' spatial dynamics (regions where solutions vary gradually in space) and the latter for the `fast' spatial dynamics (intervals with large spatial gradients).

The spatial ODE system \eqref{spatialODE} (equivalently \eqref{spatialODE-x}) has a reversibility symmetry, induced by the $x \to -x$ symmetry of \eqref{Brusselator}. 
Let
\begin{align}
    \mathcal{R} = 
    \begin{bmatrix}
        1 & 0 & 0 & 0 \\
        0 & -1 & 0 & 0 \\ 
        0 & 0 & 1 & 0 \\
        0 & 0 & 0 & -1
    \end{bmatrix}, \ \ \ \
    {\bf u} =
    \begin{bmatrix}
        u \\ p \\ v \\ q
    \end{bmatrix}, \ \ \ 
    {\bf F} = 
    \begin{bmatrix}
          p \\
          -A + (B+1) u - u^2 v\\
          \eps q\\
          \eps (-Bu + u^2 v) 
    \end{bmatrix}.
\end{align}
Here, $\mathcal{R}$ is a reversibility symmetry because it anti-commutes with the vector field:
\begin{equation}\label{reversibility}
    \mathcal{R} {\bf F} ({\bf u}) = - {\bf F}( \mathcal{R}({\bf u})).
\end{equation}
Also, let $\mathcal{L}$ denote the operator obtained by linearizing the vector field ${\bf F}$ about the homogeneous state \eqref{equilib}.
As a consequence of the reversibility symmetry, the spectrum of $\mathcal{L}$ is symmetric with respect to the real-axis and to the imaginary-axis in the $\lambda$ plane, since ${\bf u}$ is real-valued and since we have $\mathcal{R}\mathcal{L} = - \mathcal{L} \mathcal{R}$ and $(\lambda I + \mathcal{L} )^{-1} \mathcal{R} = \mathcal{R} (\lambda I - \mathcal{L} )^{-1}$. 

\begin{figure}[h!tbp]
\centering
\includegraphics[width=5in]{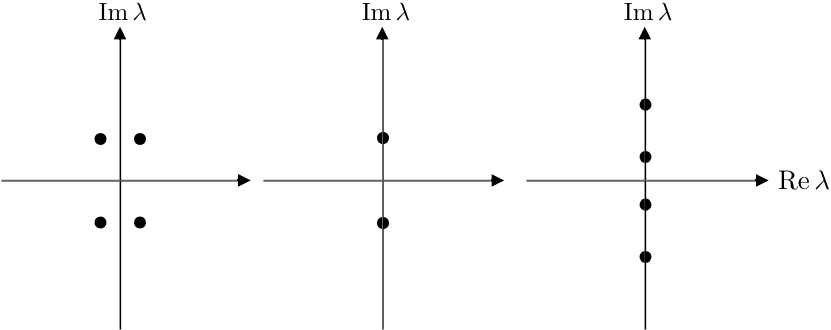}
\put(-360,136){(a)}
\put(-244,134){(b)}
\put(-134,134){(c)}
\caption{
Quartet of spatial eigenvalues $\lambda$ \eqref{quartet}. (a) $(1-\eps A)^2 < B<B_T$, (b) $B=B_T$, and (c) $B>B_T$.
}
\label{f-quartet}
\end{figure}

The (spatial) eigenvalues of the ODE system \eqref{spatialODE} are
\begin{equation}\label{quartet}
    \lambda= \pm \frac{1}{\sqrt{2}}
    \sqrt{1 - B + \eps^2 A^2 \pm \sqrt{ (1-B+\eps^2 A^2)^2 - 4 \eps^2 A^2} } \, .
\end{equation}
Therefore, as shown in Fig.~\ref{f-quartet}, 
for $B<B_T$, there is a quartet of eigenvalues with non-zero real and imaginary parts, symmetric about the axes in the $\lambda$ plane. 
These merge at $B=B_T$ to form two identical pairs of pure imaginary eigenvalues in a reversible 1:1 resonant (spatial) Hopf point, which is the Turing bifurcation point in the PDE \eqref{Brusselator}, consistent with general theory \cite{HI2011,IMD1989}.
Then, for $B>B_T$, the eigenvalues consist of two pure imaginary pairs with different imaginary parts.

Now, the normal form theory for reversible, 1:1 resonant (spatial) Hopf bifurcations guarantees the existence of families of periodic solutions, homoclinic solutions, and other solutions for \eqref{spatialODE}. 
Specifically, application of Theorem 3.21 in Chapter 4.3.3 of \cite{HI2011} yields the following proposition:

\bigskip
\begin{proposition} \label{prop:haragusiooss}
For the spatial ODE system \eqref{spatialODE} of the Brusselator, with $A>0$, $B>0$, $\eps>0$, and $B-B_T$ sufficiently small, the following statements hold:
\begin{enumerate}[(i)]
\item For all $B-B_T$ small and for $A \in \left( 0, \frac{21-\sqrt{313}}{16 \eps} \right) \cup \left( \frac{21+\sqrt{313}}{16\eps}, \infty\right)$, there is a symmetric equilibrium, a one-parameter family of periodic orbits, and a two-parameter family of quasi-periodic orbits located on KAM tori.
\item For $B > B_T$ and $A \in \left( \frac{21-\sqrt{313}}{16 \eps}, \frac{21+\sqrt{313}}{16\eps} \right)$, there is a symmetric equilibrium, a one-parameter family of periodic orbits, and a two-parameter family of quasi-periodic orbits located on KAM tori. Also, there is a pair of homoclinic orbits to the periodic orbits.
\item For $B < B_T$ and $A \in \left( 0, \frac{21-\sqrt{313}}{16 \eps} \right) \cup \left( \frac{21+\sqrt{313}}{16\eps}, \infty\right)$, 
there is a pair of homoclinic orbits to the symmetric equilibrium.
\item For $ B < B_T$ and $A \in \left( \frac{21-\sqrt{313}}{16 \eps}, \frac{21+\sqrt{313}}{16\eps} \right)$, there is a symmetric equilibrium, and there are no other bounded solutions.
\end{enumerate}
\end{proposition}

\noindent
The proof of this proposition is presented in App. A. 

\bigskip
We focus on the regime $A \in \left(0, \frac{21-\sqrt{313}}{16\eps}\right)$, {\it i.e.},  cases (i) and (iii) of Prop.~\ref{prop:haragusiooss}. 
System \eqref{spatialODE} has one-parameter families of small-amplitude, spatially-periodic solutions for all values of $B$ near $B_T$ (on both sides of $B_T$), along with the limiting homoclinics.
The Landau coefficient is positive, and the Turing bifurcation is subcritical, in agreement with the Ginzburg-Landau results \eqref{GL}.

\section{Spatially-Periodic Canard Solutions born in Turing Bifurcations}
\label{s:3}

In this section, we present many of the spatially-periodic solutions observed in the system of spatial ODEs \eqref{spatialODE} using numerical continuation \cite{HKO2018}  implemented via the software package \textsc{Auto} \cite{auto07p}.

\begin{figure}[h!tbp]
  \centering
  \includegraphics[width=4in]{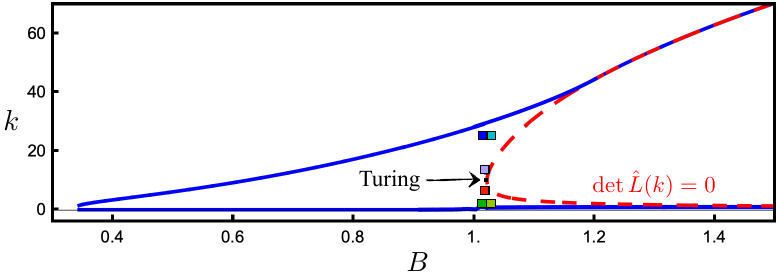}
  \caption{Regime in the $(B,k)$ plane in which spatially-periodic canard solutions were observed using a numerical continuation method, for $A=1$ and $\eps = 0.01$.
  The linear stability boundary (dashed red curve) of the equilibrium is the curve $B_T(k^2)= dk^2 + 1 + d A^2 + \frac{A^2}{k^2}$ where $\det \hat{L}(k)=0$. The local minimum is the Turing bifurcation point at $(B_T,k_T) = \left((1+\eps A)^2, \sqrt{A/\eps} \right) = (1.0201,10)$.
  The square markers correspond to the parameter values of six different spatially-periodic canards shown below in Figs.~\ref{fig:kabovekT} and \ref{fig:kbelowkT}.
  }
  \label{fig:existenceballoon}
\end{figure}

In Fig.~\ref{fig:existenceballoon}, we show the regime in the $(B,k)$ plane where spatially-periodic solutions of \eqref{spatialODE} have been obtained for $A=1$ and $\eps=0.01$ using numerical continuation. 
The solutions have spatial wavenumber $k = \tfrac{2\pi}{ \eps T}$, where $T$ is the period in \eqref{spatialODE}. 
For each $B$ in the interval shown in Fig.~\ref{fig:existenceballoon}, spatially-periodic orbits have been observed for $k$ values between the two (blue) curves.

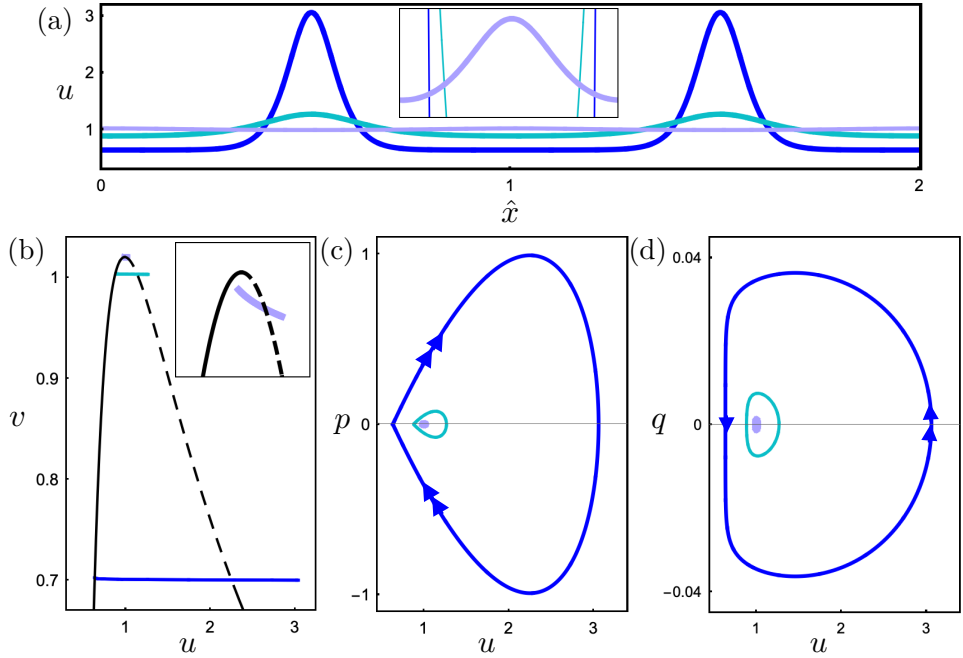
\begin{figure}[h!tbp]
  \centering
  \includegraphics[width=5in]{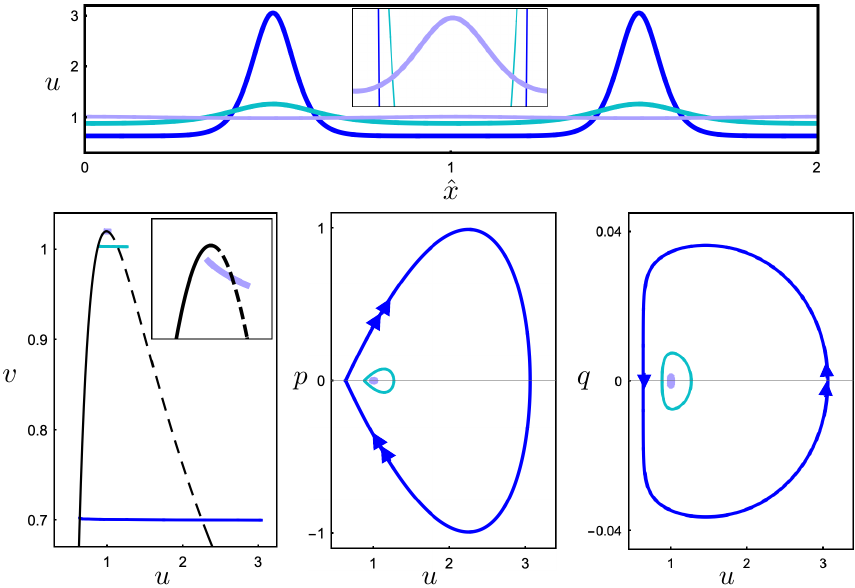}
  \put(-350,237){(a)}
  \put(-361,150){(b)}
  \put(-243,150){(c)}
  \put(-126,150){(d)}
  \caption{Spatially-periodic solutions with $B=B_T$ and $k > k_T = \sqrt{A/\eps}$ corresponding to the blue, cyan, and lavender squares in Fig.~\ref{fig:existenceballoon}. 
  (a) $u(\hat{x})$, where $\hat{x} = x/T$ and $T=25$ (blue and cyan) and $T=50$ (lavender).
  Near the Turing bifurcation, solutions can be approximately sinusoidal (lavender curve, $k\approx 12.566$) or feature pulses (blue, $k \approx 25.133$). 
  (b) $(u,v)$ projections and the $u$-nullcline (black curves). 
  Inset: magnified view of the almost sinusoidal solution. 
  (c) $(u,p)$ projections. 
  (d) $(u,q)$ projections. 
  Here, $A=1$ and $\eps=0.01$.
  }
  \label{fig:kabovekT}
\end{figure}

For $B=B_T$, some solutions with $k$ close to and {\it greater} than $k_T$ are shown in Fig.~\ref{fig:kabovekT}. 
These solutions can be close to sinusoidal (lavender solution in Fig.~\ref{fig:kabovekT}), as expected classically at Turing bifurcations. 
Also, for these same $k$, the solutions can exhibit pulsatile profiles.
Indeed, with $k$ slightly larger than $k_T$, spatially-periodic solutions are observed with fast-slow structure, exhibiting localized pulses of small-amplitude or $\mathcal{O}(1)$ amplitude, interspersed with longer intervals on which the solution slowly varies near the left branch of the $u$ nullcline (see the blue and cyan solutions in Fig.~\ref{fig:kabovekT}). 
Moreover, solutions of these types have been observed with $\mathcal{O}(1)$ periods ({\it i.e.,} with $k=\mathcal{O}(1)$) and with asymptotically small periods ({\it i.e.,} with $k\gg 1$).

\begin{figure}[h!tbp]
  \centering
  \includegraphics[width=5in]{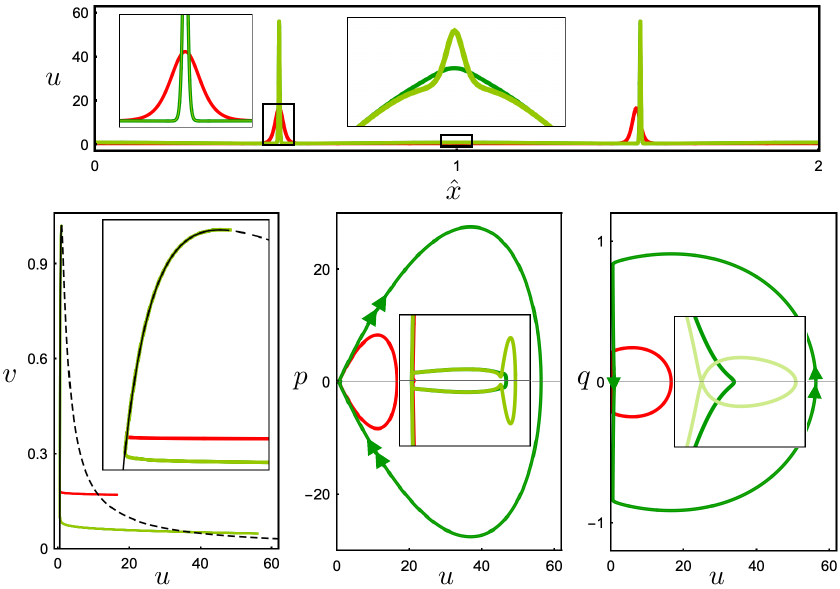}
  \put(-348,243){(a)}
  \put(-361,155){(b)}
  \put(-236,155){(c)}
  \put(-116,155){(d)}
  \caption{Spatially-periodic solutions with fast-slow structure for $B=B_T$ and $k < k_T = \sqrt{A/\eps}$ corresponding to the red, green, and olive square markers in Fig.~\ref{fig:existenceballoon}. (a) $u(\hat{x})$, with $T=100$ (red) and $T=700$ (green and olive). 
    The fast segments are narrow intervals on which the solutions exhibit sharp pulses of $\mathcal{O}\left( \eps^{-1} \right)$ amplitude.  
    The pulses consist of oscillations about the right branch (dashed black curve) of the $u$-nullcline.
  Along the long slow segments, the amplitudes are small and gradually varying.  
  The profiles of the green and olive orbits differ visually essentially only in the cusp region (middle inset).
  (b) Projection into the $(u,v)$ plane. 
  The slow segments closely follow the left branch (solid black curve) of the $u$-nullcline (inset).  
  (c) Projection into the $(u,p)$ plane. Inset: zoom on the origin.  
  (d) Projection into the $(u,q)$ plane. Inset: the (green and olive) small wavenumber solutions ($k \approx 0.898$) exhibit cusp-like behavior near the fold.
  }
   \label{fig:kbelowkT}
\end{figure}
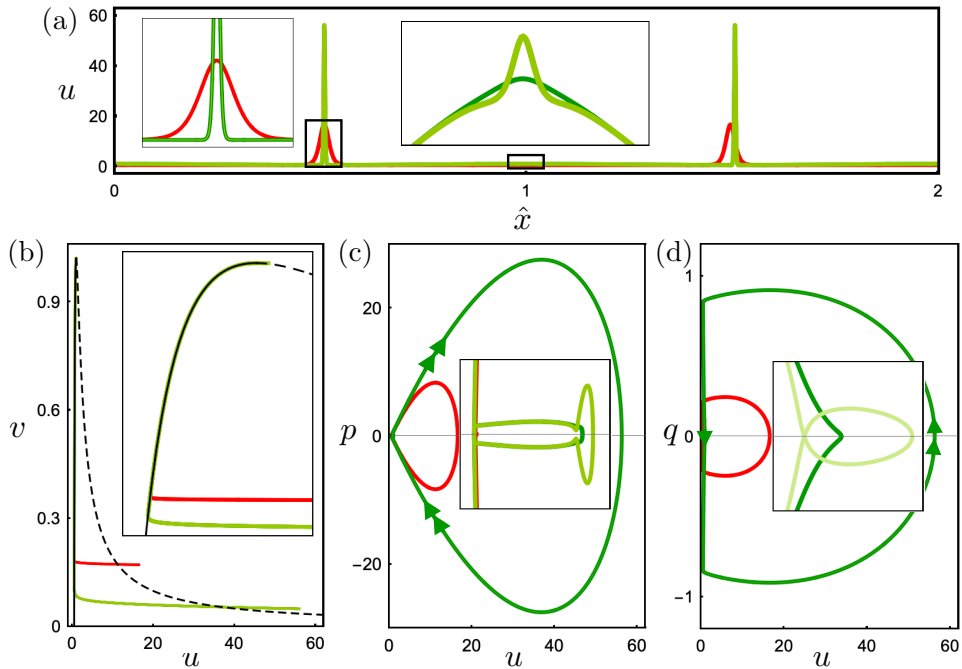

For $k$ close to but {\it less} than $k_T$, the spatially-periodic canards can also have fast-slow structure.
See the red, green, and olive orbits in Fig.~\ref{fig:kbelowkT}.
However, these solutions have asymptotically large periods, with long slow segments that closely follow the left branch of the $u$-nullcline up to the turning point, interspersed with fast ({\it i.e.}, steep) pulses in $u$ with amplitudes that range from being asymptotically small to being asymptotically large.
$v$ is constant to leading order on the narrow pulse intervals as the orbits make excursions about the right branch. 
These solutions can also spend long times near the cusp point (on the fold curve of the $u$-nullcline) and exhibit nearly self-similar dynamics. 
Also, the tiny spike centered at $\hat{x}=1$ in the green orbit will be seen to play a central role in the nucleation of new pulses during the evolution to attractors in the PDE.

Similar behavior is seen along other fixed slices in $B$. 
That is, for each fixed $B$ on either side of $B_T$, there are intervals of $k$ values in which numerical continuation reveals the existence of spatially-periodic pulse solutions with fast-slow structure. 
Concerning the fast segments of these solutions, we observe that the amplitudes decrease as $k$ is increased.
Conversely, the amplitudes can grow to $\mathcal{O}(1/\eps)$ size as $k$ is decreased.
Then, concerning the slow (flat) segments, the dynamics depend on whether $B<B_T$ or $B>B_T$.
In the former case ($B<B_T$), the slow segments converge to the stable and unstable manifolds of the saddle equilibrium in the limit as $k$ approaches its minimal value.
In contrast, for $B>B_T$, the slow (flat) segments of solutions converge to the true and faux canards of a key reversible folded saddle singularity in the limit as $k$ approaches its minimal value.

 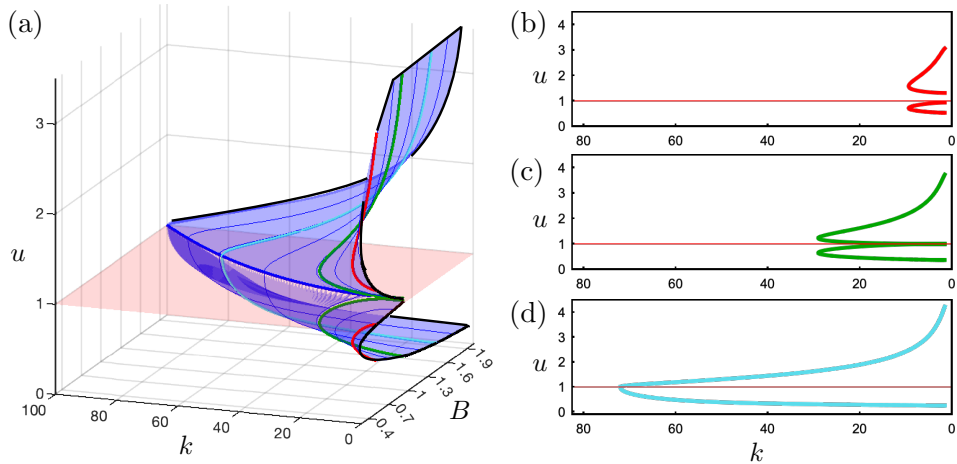
\begin{figure}[h!tbp]
    \centering
    \includegraphics[width=5in]{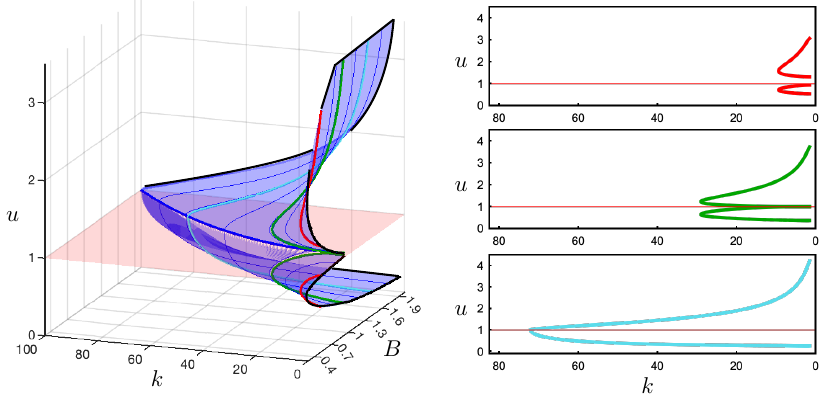}
    \put(-360,164){(a)}
    \put(-170,164){(b)}
    \put(-170,108){(c)}
    \put(-170,54){(d)}
    \caption{The existence of a saddle-node bifurcation of spatially-periodic solutions in the Brusselator spatial ODEs for $A=1$ and $\eps = 0.01$.
    (a) Projection into the $(B,k,u)$ space. 
    The red surface corresponds to the homogeneous steady state $(u,p,v,q) = (A,0,\tfrac{B}{A},0)$.
    The blue surface corresponds to a manifold of spatially periodic solutions, represented here by the extreme $u$-values for each spatially periodic solution, scaled for visual clarity. 
    That is, the upper blue surface (above the red $\{u=A\}$ plane) is obtained by the scaling $\left(\max u \right)^{1/3}$.
    The lower blue surface (below the red $\{u=A\}$ plane) is obtained by the scaling $\left( \min u \right)^{3/2}$. 
    This manifold of periodics is foliated by the slices $\{ B = {\rm constant} \}$, shown here as thin blue curves. 
    The specific cross-sections $\{ B=0.52 \}, \{ B=1.01 \}$, and $\{ B=1.52 \}$ are shown in red, green, and cyan, respectively.
    Their projections into the $(k,u)$ plane are shown in panels (b)--(d).  
    This manifold of periodics has a cusp bifurcation at the Turing point $(B_T,k_T) = ((1+\eps A)^2,\sqrt{A/\eps}) = (1.0201,10)$.
    For $B < B_T$, this manifold of periodics is folded, and both the upper and lower blue surfaces consist of two sheets separated by the saddle-node bifurcation curve (see panels (b) and (c)). 
    For $B>B_T$, the manifold consists of a single sheet on each side of the red surface, {\it i.e.}, there is no saddle-node bifurcation of periodics (see panel (d)). 
    Note that the thin black border in (a) has been added for visual clarity.
    Also, for later reference, we observe that the points $(B,k) \approx (0.34,1)$ and $(0.53,7.2)$ lie on the saddle-node curve.
    }
    \label{fig:bifurcation3d}
\end{figure}

To further elucidate the organization of the solutions in parameter space, we report (see Fig.~\ref{fig:bifurcation3d}) on the existence of saddle-node bifurcations of spatially-periodic canard solutions.
In Fig.~\ref{fig:bifurcation3d}, we show a manifold of spatially-periodic canard solutions (blue surface) is shown as a function of $(B,k)$ for $A=1$, obtained from numerical continuation. 
The manifold has a cusp point at $(B,k)=(B_T,k_T)$ for this (and other) values of $A$.
For $B < B_T$, this manifold of periodics is folded, and both the upper and lower blue surfaces consist of two sheets separated by the saddle-node bifurcation curve (see Figs.~\ref{fig:bifurcation3d}(b) and (c)). 
For $B>B_T$, the manifold consists of a single sheet on each side of the red surface, {\it i.e.}, there is no saddle-node bifurcation of periodics (see Fig.~\ref{fig:bifurcation3d}(d)). 
Similar saddle-nodes of spatially-periodic canards are observed for a wide range of values of $A$ and $0<\eps\ll 1$.

Overall, the numerics confirm Proposition~\ref{prop:haragusiooss}(i), and they extend the results in $B$ beyond where the theory holds.
Moreover, the numerics show that multiple spatially-periodic canard solutions can exist for the same $(A,B,d,k)$ parameters. 
The analysis in the next sections will elucidate these structures and geometrically deconstruct the dynamics of spatially-periodic canards.

\section{The Layer Problem, Critical and Slow Manifolds of \texorpdfstring{\eqref{Brusselator}}{Lg}, Desingularized Reduced System, and Folded Singularities}
\label{s:spatialODE}

In this section, we study the layer problem of \eqref{spatialODE}, which governs the fast spatial dynamics (pulses with steep gradients in the $u$-component).
Also, we identify the critical manifold, derive the desingularized reduced slow flow on it, and identify the key folded singularity that is responsible for generating the canard solutions of \eqref{spatialODE}.

\subsection{Layer problem}

We begin by setting $\eps=0$ in \eqref{spatialODE}.
This yields the layer system (or reduced fast system),
\begin{equation} \label{layer}
\begin{split}
    u_y &= p \\
    p_y &= - A + (B+1) u - u^2 v,
\end{split}
\end{equation}
in which $v$ and $q$ are constants.
For $v=0$, the layer problem has an equilibrium point at $\left( \frac{A}{B+1},0\right)$, and then for each $v \in \left( \left. 0,\frac{(1+B)^2}{4A}\right. \right)$, there is a pair of equilibria at 
\begin{equation} \label{usc}
(u,p)=(u_{s,c}(v), 0), \qquad {\rm where}  \quad 
    u_{s,c}(v) = \frac{1 + B \mp \sqrt{(1+B)^2 -4Av}}{2v}.
\end{equation}
The equilibria $(u_s(v),0)$ are saddle fixed points of \eqref{layer} for each $v \in \left(0,\frac{(1+B)^2}{4A}\right)$,
and $(u_c(v),0)$ are center fixed points.
For $v = \frac{(1+B)^2}{4A}$, the equilibria $u_s(v)$ and $u_c(v)$ coincide, and the eigenvalues of the linearization are both zero. 
In addition, the layer problem is Hamiltonian, with
\begin{equation}\label{H_f}
H_f (u,p;v) = \frac{1}{2} p^2 + Au - \frac{B+1}{2} u^2 + \frac{1}{3} u^3 v.
\end{equation}
For each $v \in \left(0,\frac{(1+B)^2}{4A}\right)$, the saddle point $(u_s(v),0)$ has an orbit homoclinic to it, with the center equilibrium $(u_c(v),0)$ inside. 
(An explicit formula for the homoclinics is given in Sec.~\ref{s:geoconstruction}.)

\subsection{The critical manifold and its saddle and center sheets}

\begin{figure}[h!tbp]
    \centering
    \includegraphics[width=5in]{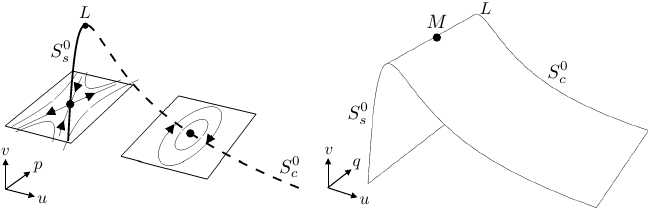}
    \put(-372,106){(a)}
    \put(-190,106){(b)}
    \caption{Plots of the saddle sheet $S^0_s$ and center sheet $S^0_c$ of the critical manifold, along with the fold curve $L$ that separates them and the RFSN-II point $M$ that lies on $L$ (a) in the $(u,p,v)$ space and (b) in the $(u,q,v)$ space. 
    Each point on $S^0_s$ is a saddle point of the two-dimensional layer problem in the $(u,p)$ plane, and each point on $S^0_c$ is a center point of the layer problem.
   }
    \label{fig:criticalmanifold}
\end{figure}

Taking the unions of the fixed points of the layer problem \eqref{usc}, we find that the critical manifold of \eqref{spatialODE} is
\begin{equation} \label{criticalmanifold}
    S^0 = \left\{
          p = 0, \,\, v = \frac{(B+1) u - A}{u^2}; \quad u > \frac{A}{B+1}, \, \, q \in \mathbb{R}  
        \right\} \, .
\end{equation}
It consists of two sheets:  $S^0  = S^0_s \bigcup S^0_c$ (see Fig.~\ref{fig:criticalmanifold}), where
\begin{equation} \label{saddle+center-sheets}
\begin{split}
    S^0_s &= \left\{
          p = 0, \, u = u_s(v); \quad 0 < v <  \frac{(1 + B)^2}{4A}  
        \right\} \\
    S^0_c &= \left\{
          p = 0, \, u = u_c(v); \quad 0 < v <  \frac{(1 + B)^2}{4A}  
        \right\}.
\end{split}
\end{equation}

Next, we examine the geometry of the saddle sheet $S^0_s$ of the critical manifold
in the four-dimensional phase space of \eqref{spatialODE} with $\eps=0$.
Over each point $(u_s(v_0),0,v_0,q_0) \in S^0_s$, there exist one-dimensional fast stable and unstable fibers. 
These are given by the local stable and unstable manifolds $W^{s,u}_{\rm loc}(u_s(v_0),0)$ of the saddle equilibrium of \eqref{layer}.
The unions of these one-dimensional fibers over all points on $S^0_s$ form the three-dimensional local stable and unstable manifolds of $S^0_s$:
\begin{eqnarray*}
    W^s(S^0_s) &=& \bigcup_{(u_s(v_0),0,v_0,q_0) \in S^0_s} W^s_{\rm loc}(u_s(v_0),0) \\ 
     W^u(S^0_s) &=& \bigcup_{(u_s(v_0),0,v_0,q_0) \in S^0_s} W^u_{\rm loc}(u_s(v_0),0).
\end{eqnarray*} 

Now, for $0<\eps \ll 1$, there is a family of slow, saddle invariant manifolds $S^\eps_s$,
\begin{equation*}
    S^\eps_s 
    = \left\{ 
 u = u_s(v) + \mathcal{O}(\eps^2), \quad 
 p = \eps \frac{ q u_s^2(v) }{ v ( u_c(v)- u_s(v) )} + \mathcal{O}(\eps^2), \quad v, q \in \mathbb{R}
    \right\},
\end{equation*}
that are $C^r$ $\mathcal{O}(\eps)$ close to $S^0_s$ for any $r>0$. 
This follows directly from Fenichel theory  \cite{F1979,J1995}, because $S^0_s$ is normally hyperbolic.
A representative $S^\eps_s$ is shown in Fig.~\ref{fig:saddleslowmanifold}.

In addition, for $0<\eps \ll 1$, the manifolds $W^s(S^0_s)$ and $W^u(S^0_s)$ persist as locally invariant stable and unstable manifolds $W^s(S^\eps_s)$ and $W^u(S^\eps_s)$ in the phase space of \eqref{spatialODE}.
We will show that $W^u(S^\eps_s)$ and $W^s(S^\eps_s)$ intersect transversely, see Sec.~\ref{s:geoconstruction}.

\subsection{Desingularized reduced vector field on the critical manifold \texorpdfstring{$S^0$}{Lg}}

\begin{figure}[h!tbp]
    \centering
    \includegraphics[width=5in]{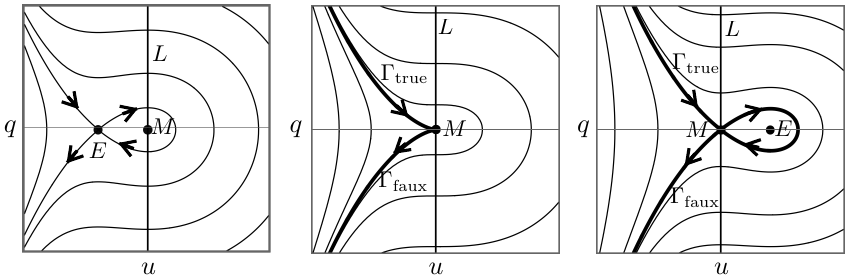}
    \put(-365,110){(a)}
    \put(-244,110){(b)}
    \put(-121,110){(c)}
    \caption{Phase planes of the desingularized reduced system \eqref{desingularizedreduced} for:(a) $0<B<1$, (b) $B=1$, and (c) $B>1$.
    The central vertical line is the fold line $L$. 
    The saddle sheet $S^0_s$ of the critical manifold lies to the left of $L$, and the center sheet $S^0_c$ to the right.
    For all $A,B>0$, \eqref{desingularizedreduced} has the equilibrium $E=\{ (u,q) = (A,0) \}$ (recall \eqref{equilib}) and a reversible folded singularity $M = \left\{ (u,q)=\left(\tfrac{2A}{B+1},0 \right) \right \}$.
    (a) For $0<B<1$, $M$ lies inside the homoclinic orbit to $E$, and there are no singular canard solutions. 
    (b) For $B=1$, $M$ and $E$ coincide in a RFSN-II point.
    Its true and faux canards, $\Gamma_{\rm true}$ and $\Gamma_{\rm faux}$, are the single-branched stable and unstable manifolds of $M$ that lie on $S^0_s$.
    They meet at $M$ in a cusp. 
    (c) For $B>1$, $M$ is an RFS point. 
    Its true and faux canards, $\Gamma_{\rm true}$ and $\Gamma_{\rm faux}$, 
    lie on the stable and unstable manifolds of $M$, and both have components on $S^0_s$ and $S^0_c$.
    On $S^0_c$, the two canards form a homoclinic to $M$ that surrounds $E$.
    }
    \label{fig:desingularized}
\end{figure}

On the critical manifold $S^0$, the equations are given by \eqref{spatialODE-x} with $\eps=0$,
\begin{equation}\label{onS0}
\begin{split}
    0 &= p \\
    0 &= -A + (B+1) u - u^2 v \\
    v_x &= q \\
    q_x &= -Bu + u^2 v.
\end{split}
\end{equation}
Let $\mathcal{F}(u,v,A,B) = - A + (B+1)u - u^2 v$.
Differentiating the constraint $\mathcal{F}=0$ with respect to $x$, we find
$\mathcal{F}_u \frac{du}{dx} + \mathcal{F}_v  \frac{dv}{dx}
=\left( (B+1) - 2 u v \right) \frac{du}{dx} - u^2 \frac{dv}{dx} = 0$ on $S^0$.
Hence, on $S^0$,
\begin{equation*}
    \begin{bmatrix}
        \frac{-\mathcal{F}_u}{\mathcal{F}_v} & 0 \\ 0 & 1 
    \end{bmatrix}
    \begin{bmatrix}
        u_x \\
        q_x
    \end{bmatrix}
    = 
    \begin{bmatrix}
        0 & 1 \\
        1 & 0 
    \end{bmatrix}
    \begin{bmatrix}
        u \\ q 
    \end{bmatrix}
    +
    \begin{bmatrix}
        0 \\
        -A
    \end{bmatrix}.
\end{equation*}
Then, left-multiplying by the adjoint matrix and introducing the dynamic independent variable $\frac{d}{dx_d} = \frac{-\mathcal{F}_u}{\mathcal{F}_v}|_{S^0} \, \, \frac{d}{dx}$, we find that the vector field on $S^0$ is
\begin{eqnarray*}
u_{x_d} &=& q \nonumber \\
q_{x_d} &=& g(u;A,B) = \frac{-\mathcal{F}_u}{\mathcal{F}_v}|_{S^0} \left( u - A \right).    
\end{eqnarray*}
Here, $\frac{-\mathcal{F}_u}{\mathcal{F}_v}|_{S^0} = \frac{2A - (B+1)u}{u^3}$.
Hence, the desingularized reduced system on $S^0$ is
\begin{equation}\label{desingularizedreduced}
\begin{split}
u_{x_d} &= q  \\
q_{x_d} &= \frac{-(B+1)}{u} + \frac{A(B+3)}{u^2} - \frac{2A^2}{u^3}.    
\end{split}
\end{equation}
The desingularized reduced vector field \eqref{desingularizedreduced} is illustrated in Fig.~\ref{fig:desingularized}.
The trajectories lie on the level sets of the Hamiltonian,
 \begin{equation} \label{eq:desingularizedham}
     H_d(u,q;A,B) = \frac{1}{2} q^2 + (B+1) \ln (u) + \frac{A(B+3)}{u} - \frac{A^2}{u^2}.
 \end{equation}

\subsection{The RFSN-II point for  \texorpdfstring{$B=1$}{Lg} and the RFS points for \texorpdfstring{$B>1$}{Lg}}

For all $A,B>0$, the desingularized reduced vector field \eqref{desingularizedreduced} vanishes at the point
\begin{equation}\label{RFSpointM}
M = \left\{ (u,q) = \left( \frac{2A}{B+1}, 0 \right) \right\}.
\end{equation}
It is a folded singularity of \eqref{desingularizedreduced}, because it lies on the fold line $L$ where $q=0$ and $\frac{-\mathcal{F}_u}{\mathcal{F}_v}|_{S^0} = 0$.
Moreover, $M$ is a reversible folded singularity
due to the reversibility symmetry of \eqref{spatialODE}.
We add that $M$ is not an equilibrium of the full slow flow.

For $0<B<1$, $M$ is a reversible folded center, while $E$ is a saddle fixed point.
See Fig.~\ref{fig:desingularized}(a).
Then, for $B>1$, the stability is reversed. 
Here, $M$ is a reversible folded saddle (RFS), 
\begin{equation}
\label{M-RFS}
M_{\rm RFS} = \left\{  (u,q;A,B) : u = \tfrac{2A}{B+1}, \, q=0 ; \, \, A>0, \, B>1 \right\},
\end{equation} 
while $E$ is a center. 
Further, the true canard, $\Gamma_{\rm true}$, of $M_{\rm RFS}$
corresponds to its (two-branched) stable manifold, and the faux canard, $\Gamma_{\rm faux}$, corresponds to its (two-branched) unstable manifold. 
See Fig.~\ref{fig:desingularized}(c).
(Here, we recall that the linear stability of the singularities $M_{\rm RFS}$ and $E$ is given by the Jacobian matrix of \eqref{desingularizedreduced},
\begin{equation*}
    \begin{bmatrix}
        0 & 1 \\
        g_u & 0
    \end{bmatrix},
\end{equation*}
where $g_u=\left( 6A^2 - 2A(3+B)u + (B+1)u^2\right) u^{-4}$. 
Note that $g_u|_{E} = \frac{1-B}{A^2}$, 
and $g_u|_{M} = \frac{(B+1)^3 (B-1)}{8A^2}$. 
Then, their nonlinear stability is determined from the Hamiltonian $H_d$, \eqref{eq:desingularizedham}.)

Exactly for $B=1$, $M$ merges with $E$ at $(A,0)$, and it is an RFSN-II point, 
\begin{equation}
\label{M-RFSNII}
M_{\rm RFSN-II} = \{ (u,q;A,B): u=A, \, q=0; \, \, A>0, \, B=1 \}.
\end{equation}
See Fig.~\ref{fig:desingularized}(b).
This RFSN-II point occurs at the singular limit ($d \to 0$) of the Turing bifurcation $B_T$ of the homogeneous state $(u,v)=\left( A, \frac{B}{A}\right)$, recall \eqref{kTBT}.
Also, the stable and unstable manifolds of this RFSN-II point have one branch each, and they are cusp shaped (see Fig.~\ref{fig:desingularized}(b)). 
They correspond to the true and faux canards, $\Gamma_{\rm true}$ and $\Gamma_{\rm faux}$, of the folded singularity.

\begin{figure}[h!tbp]
\centering
\includegraphics[width=3.5in]{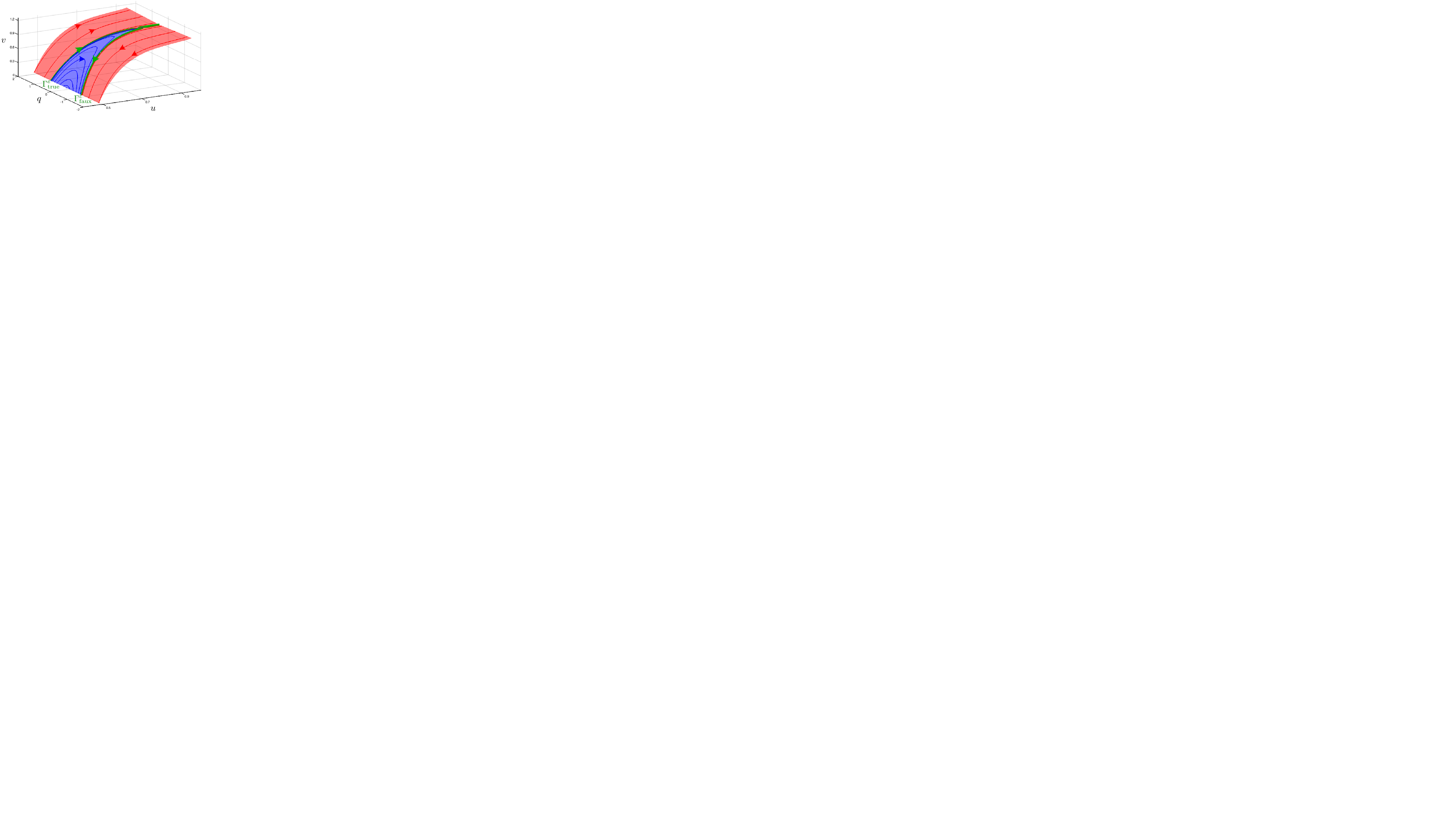}
\caption{A plot of a saddle slow manifold $S^\eps_s$ in the $(u,q,v)$ space, along with the true and faux canards (green curves). 
Solutions in the blue region enclosed by the true and faux canards head toward -and then turn away from- the fold curve. 
Solutions in the complementary red region reach the neighborhood of the fold where they exhibit a fast jump. 
Here, $A=1$, $B=1.1$, and $\eps=0.01$.}
\label{fig:saddleslowmanifold}
\end{figure}

We will show in Secs.~\ref{s:desingularization} and \ref{s:K1} that the true canard $\Gamma_{\rm true}$ persists as a maximal canard $\Gamma_{\rm true}^{\eps}$ 
for $0 < \eps \ll 1$. 
Similarly, the faux canard $\Gamma_{\rm faux}^{\eps}$ persists as a maximal canard $\Gamma_{\rm faux}^{\eps}$ for $0< \eps \ll 1$. 
See Fig.~\ref{fig:saddleslowmanifold}.
Moreover, we will show that for all $B>1+\mathcal{O}(\eps)$, the connection $\Gamma_{\rm true} \cup \Gamma_{\rm faux}$ persists as a maximal canard solution $\Gamma_{\rm true}^{\eps} \cup \Gamma_{\rm faux}^{\eps}$ for $0<\eps \ll 1$. 
Numerically, we see evidence of this in Fig.~\ref{fig:saddleslowmanifold} where the (green) maximal canard boundaries that separate the red and blue sheets of the saddle slow manifold form a single continuous curve.  

For clarity, we use the terms ``true canard'' ($\Gamma_{\rm true}$) and ``faux canard'' ($\Gamma_{\rm faux}$) to refer specifically to the stable and unstable manifolds of the RFS or RFSN-II in the singular limit $\eps \to 0$. We will show in Section~\ref{subsec:persistence} that $\Gamma_{\rm true}$ and $\Gamma_{\rm faux}$ serve as good approximations of their non-singular counterparts $\Gamma_{\rm true}^{\eps}$ and $\Gamma_{\rm faux}^{\eps}$ for $0<\eps \ll 1$.

Overall, the folded singularity $M_{\rm RFSN-II}$ \eqref{M-RFSNII} is a key organizing point in the phase space of the spatial ODE system. 
It sits on the boundary of the parameter regime in which the homogeneous steady state is normally hyperbolic (recall the transition from Fig.~\ref{f-quartet}(a) to (b)).
Along with its true and faux canards, $M_{\rm RFSN-II}$ is responsible for generating the canard segments of the periodic solutions in the neighborhood of the Turing bifurcation point $B_T$.
Then, for parameters past the Turing point, {\it i.e.,} for parameter values $B > B_T$, with $|B-B_T| = \mathcal{O}(1)$, the true and faux canards of $M_{\rm RFS}$ \eqref{M-RFS} are the mechanisms responsible for generating the canard segments of the spatially-periodic patterns. 

\begin{remark}
Since we performed the dynamic coordinate transformation $\frac{d}{dx_d} = \frac{-\mathcal{F}_u}{\mathcal{F}_v}|_{S^0} \, \, \frac{d}{dx}$, the desingularized reduced system \eqref{desingularizedreduced} is topologically equivalent to the slow flow \eqref{onS0} in the regions of phase space where $\frac{-\mathcal{F}_u}{\mathcal{F}_v}|_{S^0} >0$. 
Hence, on the saddle sheet $S^0_s$, \eqref{desingularizedreduced} is topologically equivalent to the slow flow. 
However, on the center sheet $S^0_c$, the flow of \eqref{desingularizedreduced} shown in Fig.~\ref{fig:desingularized} is in the opposite direction to the flow of the slow system \eqref{onS0}, since $\frac{-\mathcal{F}_u}{\mathcal{F}_v}|_{S^0_c} < 0$. 
\end{remark}

\section{Desingularizing \texorpdfstring{$M_{\rm RFSN-II}$ \eqref{M-RFSNII}}{Lg}: The Rescaling Chart \texorpdfstring{$K_2$}{Lg}} \label{s:desingularization}

In this section, we apply the method of geometric desingularization to the spatial ODE \eqref{spatialODE}, and we analyze the dynamics in the rescaling (or central) chart of the blown-up locus.
We identify the algebraic solution that defines the true and faux canards of $M_{\rm RFSN-II}$ in the rescaling chart.

\subsection{Translating the folded singularity to the origin}

As a preliminary step, we translate the RFSN-II at $(u,p,v,q)=\left(A, 0, \tfrac{1}{A},0 \right)$ to the origin. Let 
\begin{equation}
\label{translateFS2origina}    
d=\varepsilon^{2}, \quad  
\B=B-1, \quad 
u=A+U, \quad 
p=P, \quad
v=\tfrac{1}{A}+V, \quad
q=Q.
\end{equation}
The spatial ODE \eqref{spatialODE} becomes
\begin{equation}
\label{Fast Variable Steady State System}
\begin{split}
        U_{y}&= P \\
        P_{y}&= A\B - A^{2}V - \tfrac{1}{A}U^{2} - F(U,V;A,\B)  \\
        V_{y}&= \varepsilon Q  \\
        Q_{y}&= \varepsilon \left(-A\B + U + A^{2}V + \tfrac{1}{A}U^{2} + F(U,V;A,\B)\right), 
\end{split}
\end{equation}
where
\begin{align}
    F(U,V;A,\B)&= -\B U + 2AUV + U^{2}V.
\end{align}
It consists of the terms that will be seen to be higher order.
This is the system that we will analyze in this section, and in the next section.

\subsection{Desingularization}

To desingularize the origin, we append the equations \(\eps_{y}=0\) and \(\B_{y}=0\) to the system \eqref{Fast Variable Steady State System} and blow-up the nilpotent point at the origin in \(\R^{6}\).
Let 
\((\bar{U},\bar{P},\bar{V},\bar{Q},\bar{\varepsilon},\bar{\B})\in \mathcal{S}^{5}\) and \(r\in[0,r_{max}]\) for some \(r_{max}>0\) to be determined later. 
We define 
\begin{equation} \label{powersofr}
U=r^2\bar{U}, \quad 
P=r^3\bar{P}, \quad 
V=r^4\bar{V}, \quad 
Q=r^3\bar{Q}, \quad 
\varepsilon=r^2\bar{\varepsilon}, \quad 
\B=r^4\bar{\B}.
\end{equation}
The powers of $r$, which is the new dynamic small variable, are chosen to produce simplified dynamics in a small neighborhood of the RFSN-II point.
Specifically, the \(\B\) component is scaled by 2 powers of \(r\) more than the significant degeneration, leading to simplified \(B\) dependence in later calculations.
We will consider two charts for \(\mathcal{S}^{5}\). 
The rescaling chart, \(K_{2}\), is defined by \(\bar{\varepsilon}=1\), and the dynamics in $K_2$ are analyzed in this section.
The entry/exit chart, \(K_{1}\), is defined by \(\bar{U}=-1\), and the dynamics in chart $K_1$ are analyzed in Sec.~\ref{s:K1}.

In the chart $K_2$, where $\eps_2 = 1$ is fixed, the coordinates are given by
$U=r_2^2 U_2, P=r_2^3 P_2, V = r_2^4 V_2, Q=r_2^3 Q_2, \eps=r_2^2,$ and $\B = r_2^4 \B_2$.
The dynamics in the \(K_{2}\) chart are governed by:
\begin{equation}
\label{K2 Dynamics}
\begin{split}
        U_{2}'&= P_{2} \\
        P_{2}'&= A\B_2 - A^{2}V_{2} - \tfrac{1}{A}U_{2}^{2} - r_2^2 F_2 \\
        V_{2}'&= Q_{2} \\
        Q_{2}'&= U_2 - r_2^2 \left( A\B_{2} - A^{2}V_{2} - \tfrac{1}{A}U_{2}^{2}
                 - r_{2}^{2}F_2 \right)  \\
        r_{2}'&= 0 \\
        \B_{2}'&= 0.
\end{split}
\end{equation}
Here, the prime denotes the derivative with respect to the variable \(x_{2}\) defined by \(x_{2}=r_{2}y\), and 
\begin{align}
    F_2(U_{2},V_{2},\B_{2},r_{2};A)= -\B_2 U_2 + 2AU_{2}V_{2} + r_2^2 U_{2}^{2}V_{2}.
\end{align}

\begin{remark}
As with all systems of differential equations that are analyzed using the method of geometric desingularization, it may at first seem counter-intuitive to increase the dimension of the system by one (the new equation for $r$).
However, in each of the relevant charts, one of the variables is held constant.
Hence, the systems in the charts have the same dimension as the original system, and they are usually simpler to analyze than the original system. 
\end{remark} 

\subsection{A key algebraic solution in chart \texorpdfstring{$K_2$}{Lg}}

The following lemma (which follows by a direct calculation) establishes the existence of a special algebraic solution in chart $K_2$ in the $\{ r_{2}=0 \}$ subspace.

\begin{lemma} \label{lem-gamma0}
For every $\B_2$, there exists an algebraic solution, \(\Gamma_{0}(x_{2})\), on $\{ r_2=0 \}$ given by 
\begin{eqnarray} 
    (U_2, P_2, V_2, Q_2)=
        \left(
        \tfrac{-A^{3}}{12}x_{2}^{2}, \, 
        \tfrac{-A^{3}}{6}x_{2}, \, 
        \tfrac{-A^{3}}{144}x_{2}^{4} + \tfrac{A}{6}+\tfrac{\B_{2}}{A}, \,
        \tfrac{-A^{3}}{36}x_{2}^{3} 
        \right).
\end{eqnarray}
The orbit $\Gamma_0(x_2)$ consists of the following two halves:
\begin{equation} \label{2halvesofGamma0}
    \Gamma_{\rm true} = \Gamma_0 \vert_{x_2 \le 0} 
    \qquad{\rm and}
    \qquad
    \Gamma_{\rm faux} = \Gamma_0 \vert_{x_2 \ge 0}
\end{equation}
that correspond to the true and faux canards, respectively, of the RFSN-II singularity.
\end{lemma}

\bigskip

The true canard approaches the RFSN-II point as $x_2 \to 0^-$, and the faux canard emerges from the RFSN-II point as $x_2>0$ increases from zero.
They are cusp-shaped curves, with implicit parametrization
$Q_2 = \mp \tfrac{2}{\sqrt{3} A^{3/2}}\left(- U_2\right)^{3/2} = \pm \tfrac{2}{\sqrt{3} A^{3/2}}\vert U_2\vert  \left( -U_2 \right)^{1/2}$. 

As $x_2 \to \pm \infty$, the special solution $\Gamma_0(x_2)$ limits on equilibria on the equator of $\mathcal{S}^5$. 
These points are naturally expressed in the $K_1$ chart, and they are determined in the following lemma:

\begin{lemma}\label{lem-limitsGamma0}
The limits of $\Gamma_0(x_2)$ as $x_2 \to \pm \infty$ are the following points on the equator of $\mathcal{S}^5$:
\begin{equation}\label{limitGamma0}
    {\bf p}_{\pm} = \left(r_{1},P_{1},V_{1},Q_{1},\varepsilon_{1},\B_{1}\right)
    =\left(0,0,\tfrac{-1}{A^{3}},\tfrac{\mp 2}{\sqrt{3A^{3}}},0,0\right).
\end{equation}
\end{lemma} 

\bigskip
\noindent
{\bf Proof.}
The change of variables $\kappa_{21}: K_{2} \to K_{1}$, defined on \(U_{2}<0\), transports objects from the rescaling chart $K_2$ into the entry/exit chart $K_1$. It is given by
\begin{equation}
\label{K2 algebraic solutions}
\begin{split}
        r_{1} &= r_{2}\sqrt{-U_{2}},\,\,
        P_{1}=\frac{P_{2}}{(-U_{2})^{3/2}},\,\,
        V_{1}=\frac{V_{2}}{U_{2}^{2}}, \,\,
        Q_{1}= \frac{Q_{2}}{(-U_{2})^{3/2}},\,\,
        \varepsilon_{1}=-\frac{1}{U_{2}},\,\,
        \B_{1}=-\frac{\B_{2}}{U_{2}^{2}}.
\end{split}
\end{equation}
In the \(K_{1}\) chart, \(\Gamma_{0}\) lies in the set $\{ r_1 = 0\}$ and is given by
\begin{equation*}
\begin{split}
        P_{1}&= \tfrac{\mp4\sqrt{3}}{A^{\hfrac{3}{2}}x_{2}^2}, \quad 
        V_{1} = \left(\tfrac{-A^3}{144} x_2^4 + \tfrac{A}{6} + \tfrac{\B_{2}}{A}
        \right) \left( \tfrac{144}{A^6 x_2^4}\right), \quad 
        Q_{1}= \tfrac{\mp 2}{\sqrt{3}A^{\hfrac{3}{2}}},  \quad 
        \varepsilon_{1}= \tfrac{12}{A^3 x_2^2},  \quad 
        \B_{1}= \tfrac{144 \B_{2}}{A^6 x_2^4}.
\end{split}
\end{equation*}
Then, taking the limit as $x_2 \to \mp \infty$, we obtain \eqref{limitGamma0}. \hfill 
$\Box$

\bigskip
\noindent 
We also want to know the limits of the unit tangent vector to \(\Gamma_{0}\).

\begin{lemma}\label{lem-limitstangentvectors2Gamma0}
The limits of the unit tangent vector to $\Gamma_0$ as \(x_2 \to\pm\infty\) are
\begin{align}\label{limitGamma0Tangent}
    \lim_{x_{2}\to\pm \infty}\frac{\Gamma_{0}'(x_2)}{\norm{\Gamma_{0}'(x_2)}}=&\left(0,\pm \sqrt{1-\tfrac{3}{A^{3}+3}},0,0,-\sqrt{\tfrac{3}{A^{3}+3}},0\right)
\end{align}
\end{lemma} 

\noindent
{\bf Proof.}
We compute the derivatives of the different components in the vector $\Gamma_0(x_2)$,
\begin{align}\label{derivatives} 
    \frac{dr_1}{d{x_{2}}}&= 0, \quad  
    \frac{dP_1}{dx_{2}}= \frac{\pm8\sqrt{3}}{A^{3/2} x_{2}^{3}}, \quad 
    \frac{dV_1}{dx_{2}}= \frac{-576}{A^{6} x_{2}^{5}} \left(\frac{A}{6} + \frac{\B_{2}}{A} \right), \nonumber \\
    \frac{dQ_1}{dx_{2}}&= 0,\quad
    \frac{d\eps_1}{dx_{2}}= \frac{-24}{A^{3} x_{2}^{3}}, \quad \quad 
    \frac{d\B_1}{dx_{2}}=\frac{-576\B_{2}}{A^{6} x_{2}^{5}}.
\end{align}
Hence, we find that \(\norm{\Gamma_{0}'}\) is given by 
\begin{align*}
    \norm{\Gamma_{0}'}= \tilde{\gamma}_0 x_{2}^{-3} + \mathcal{O}\left( x_{2}^{-5} \right),
    \quad {\rm where} \quad \tilde{\gamma}_{0} =
     \tfrac{8}{A^{3}} \sqrt{3(A^{3}+3)} \, ,
\end{align*}
Now, combining this with \eqref{derivatives}, we obtain the stated results for the limits in the lemma.
\hfill $\Box$

\begin{remark}
    A similar calculation yields the formulas for the singular true and faux canards of the RFS point, $M_{\rm RFS}$ given by \eqref{M-RFS}.
    See App.~\ref{sec:App-RFS}.
\end{remark}

\section{Desingularizing the RFSN-II Point: The Entry/Exit Chart \texorpdfstring{\(K_{1}\)}{Lg}}\label{s:K1}

In the entry/exit chart $K_1$, where $U_1=-1$, the variables are defined by
$U=-r_1^2, P=r_1^3 P_1, V=r_1^4 P_1, Q=r_1^3 Q_1, \eps=r_1^2 \eps_1,$ and $\B = r_1^4 \B_1$.
We substitute these into \eqref{Fast Variable Steady State System} and find that the dynamics in the chart \(K_{1}\) are governed by the following system of ODEs:
\begin{equation}
\label{K1 dynamics}
\begin{split}
        r_{1}'&= \tfrac{-1}{2}r_{1}P_{1} \\
        P_{1}'&= \tfrac{3}{2}P_{1}^{2} + A\B_1 
        - A^{2}V_{1} - \tfrac{1}{A} - r_1^2 F_1 \\
        V_{1}'&= 2P_{1}V_{1} + \varepsilon_{1}Q_{1} \\ 
        Q_{1}'&= \tfrac{3}{2}P_{1}Q_{1} - \eps_1 
        - r_{1}^{2}\eps_{1}\big(A\B_1 
        - A^{2}V_{1} - \tfrac{1}{A} - r_1^2 F_1 \big)\\
        \varepsilon_{1}'&=\varepsilon_{1}P_{1} \\
        \B_{1}'&=2\B_{1}P_{1}.
\end{split}
\end{equation}
%
Here, the prime has been recycled to denote the derivative with respect to the variable \(x_{1}\) defined by \(x_{1}=r_{1}y\), and the function \(F_1 = F_1(-1,V_{1},\B_{1},r_{1};A)\) is given by
\begin{align}
    F_1(-1,V_{1},\B_{1},r_{1};A)&= \B_1 - 2AV_{1}
    +r_1^2 V_{1}.
\end{align}

There are a variety of invariant hyperplanes in this chart, defined by \(\left\{r_{1}=0\right\}\), \(\left\{\varepsilon_{1}=0\right\}\), and \(\left\{\B_{1}=C\eps_{1}\right\}\), where \(C\) may be any real number, including zero. 
In the analysis, we focus on the coordinate hyperplanes.
To study the full dynamics in \(K_{1}\), we find it useful to examine the dynamics in these invariant subspaces and in their intersections, systematically building out from the lowest, three-dimensional subspace to the full phase space.
In Sec.~\ref{ss:3D}, we present the analysis in the three-dimensional subspace $\{ r_1, \eps_1, \B_1 = 0 \}$, which is the intersection of the three coordinate hyperplanes.
Then, the analysis of \eqref{K1 dynamics} in the full, six-dimensional phase space is given in Sec.~\ref{ss:full6D}.
It is built upon the analyses (not shown) of the intermediate systems in the different four-dimensional subspaces 
$\{ r_1, \eps_1=0 \}$,
$\{ \eps_1, \B_1 = 0\}$,
$\{ r_1, \B_1 = 0 \}$,
and then successively also in the different five-dimensional invariant subspaces
$\{ \eps_1=0 \} $,
$\{ r_1 = 0 \}$,
and $\{ \B_1 = 0 \}$.

\subsection{The subspace \texorpdfstring{\(\left\{r_{1},\varepsilon_{1},\B_{1}=0\right\}\)}{Lg}}\label{ss:3D}
In the invariant subspace \(\left\{r_{1},\varepsilon_{1},\B_{1}=0\right\}\), the equations \eqref{K1 dynamics} are the three-dimensional system,
\begin{equation}
\begin{split}
\label{K2 reB=0 Dynamics}
        P_{1}'&=\tfrac{3}{2}P_{1}^{2}-A^{2}V_{1}-\tfrac{1}{A} \\
        V_{1}'&=2P_{1}V_{1} \\
        Q_{1}'&=\tfrac{3}{2}P_{1}Q_{1}.
\end{split}
\end{equation}
Within this subspace, the equation for \(Q_{1}\) decouples, and 
$ Q_{1}(x_{1})=Q_{1}(0)e^{\tfrac{3}{2}\int_{0}^{x_{1}}P_{1}(s)ds}.$

The system \eqref{K2 reB=0 Dynamics} has a line of saddle equilibria,
\begin{align}\label{l1 def}
    \ell_{1}=&\left\{r_{1}=0, \, P_{1}=0, \, V_{1}=\tfrac{-1}{A^{3}}, \, Q_{1}\in\R, \, \varepsilon_{1}=0, \, \B_{1}=0\right\}.
\end{align}
The eigenvalues at any point on $\ell_{1}$ are \(\lambda=\pm\sqrt{\hfrac{2}{A}},0\), and the eigenspaces are
\begin{equation*}
    \E_{\sqrt{\hfrac{2}{A}}}=\vecspan{\begin{bmatrix}
        \sqrt{\hfrac{2}{A}}\\
        \hfrac{-2}{A^{3}}\\
        \tfrac{3}{2}Q_{1}
    \end{bmatrix}},
    \quad 
    \E_{-\sqrt{\hfrac{2}{A}}}
    =\vecspan{\begin{bmatrix}
        -\sqrt{\hfrac{2}{A}}\\
        \hfrac{-2}{A^{3}}\\
        \tfrac{3}{2}Q_{1}
    \end{bmatrix}},
    \quad \E_{0}
    = Q_1 {\rm -axis}.
\end{equation*}

\begin{figure}[h!tbp]
    \centering
    \includegraphics[width=5in]{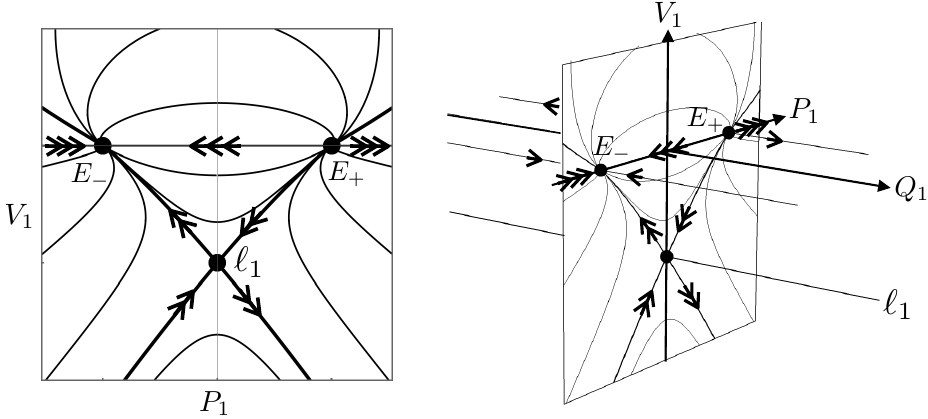}
    \put(-360,152){(a)}
    \put(-188,152){(b)}
    \caption{Dynamics in the $\{ r_1, \eps_1, \B_1=0 \}$ subspace, with the line $\ell_1$ \eqref{l1 def} of saddle fixed points, the invariant line $I$ \eqref{I def}, and the equilibria at $E_{\pm}$ \eqref{E+- def}. (a) Projection into the $(P_1,V_1)$ plane. (b) Dynamics in the three-dimensional $(P_1,V_1,Q_1)$ phase space.}
    \label{fig:chart1}
\end{figure}

The system \eqref{K2 reB=0 Dynamics} also has an invariant line,
\begin{align}\label{I def}
    I=&\left\{r_{1}=0,P_{1}\in\R,V_{1}=0,Q_{1}=0,\varepsilon_{1}=0,\B_{1}=0\right\}.
\end{align}
On $I$, the system reduces to $P_1'=\tfrac{3}{2}P_1^2 - \tfrac{1}{A}$, and there are two isolated equilibria
\begin{align}\label{E+- def}
    E_{\pm} = \left\{ r_{1}=0,P_{1}=\pm \rho, V_{1}=0,Q_{1}=0,\varepsilon_{1}=0,\B_{1}=0 \right\},
\end{align}
where we define $\rho = \sqrt{\tfrac{2}{3A}}$.
The eigenvalues at \(E_{\pm }\) are \(\lambda=\pm 3\rho,\pm 2\rho,\pm \tfrac{3}{2} \rho\), with eigenspaces
\begin{equation*}
    \E_{\pm 3\rho}= P_1{\rm -axis}, \quad 
    \E_{\pm 2\rho}=\vecspan{\begin{bmatrix}
        A^{2}\\
        \pm \rho\\
        0
    \end{bmatrix}},
    \quad 
    \E_{\pm 3\rho/2}
    =Q_1 {\rm -axis}.
\end{equation*}
Hence, $E_{\pm}$ are a source and a sink, respectively.
The dynamics of \eqref{K2 reB=0 Dynamics} are illustrated in Fig.~\ref{fig:chart1}.

Lastly, we observe that the trajectories of the $(P_1,V_1)$ system can be found explicitly with the transformation \(V_{1}=\left(\tfrac{3}{2}P_{1}^{2}-\tfrac{1}{A}\right)\tilde{V}_{1}\), which is invertible for \(P_{1}\neq \pm \rho\). 
The dynamics are given by
\begin{equation*}
\begin{split}
    P_{1}' &= \left(\tfrac{3}{2}P_{1}^{2}-\tfrac{1}{A}\right)\left(1-A^{2}\tilde{V}_{1}\right),\\
    \tilde{V}_{1}' &= P_{1}\tilde{V}_{1}\left(3A^{2}\tilde{V}_{1}-1\right).
\end{split}
\end{equation*}
The sets \(\left\{\tilde{V}_{1}=0\right\}\) and \(\left\{\tilde{V}_{1}=\frac{1}{3A^{2}}\right\}\) are invariant, corresponding to invariant parabolas in the original variables. 
Dividing the first equation by the second and separating variables, we find that the
trajectories are given by
\begin{align*}
    P_{1}=& \pm\sqrt{\tfrac{2}{3}} 
    \sqrt{\tfrac{1}{A} 
    + C\tfrac{\left(3A^{2}\tilde{V}_{1}-1\right)^{2}}{\tilde{V}_{1}^{3}}} \, ,
\end{align*}
where $C$ is an integration constant. 
Converting back to the original variables, we find 
\begin{align}\label{PofV1}
    P_{1}=\pm\sqrt{\tfrac{2}{3A}+2A^{2}V_{1}-\tfrac{2}{3\sqrt{C}}\abs{V_{1}}^{\hfrac{3}{2}}}.
\end{align}

\subsection{The dynamics of the six-dimensional system in chart \texorpdfstring{$K_1$}{Lg}} \label{ss:full6D}
The full dynamics in the entry/exit chart \(K_{1}\) are governed by the six-dimensional system \eqref{K1 dynamics}, which we recall is: 
\begin{equation}\label{K1 dynamics2}
\begin{split}
        r_{1}'&= -\tfrac{1}{2}r_{1}P_{1} \\
        P_{1}'&= \tfrac{3}{2}P_{1}^{2} + A\B_1 
        - A^2 V_{1} - \tfrac{1}{A} - r_1^2 F_1 \\
        V_{1}'&= 2P_{1}V_{1}+\varepsilon_{1}Q_{1} \\
        Q_{1}'&= \tfrac{3}{2}P_{1}Q_{1} - \eps_1 
        - r_{1}^{2}\eps_{1}\big(A\B_1 
        - A^{2}V_{1} - \tfrac{1}{A} - r_1^2 F_1 \big) \\
        \varepsilon_{1}'&=\varepsilon_{1}P_{1}  \\
        \B_{1}'&= 2\B_{1}P_{1}.
\end{split}
\end{equation}

\noindent
In this section, we establish the following lemma about the dynamics of \eqref{K1 dynamics2}:

\begin{lemma}\label{lem-K1}
System \eqref{K1 dynamics2} (equivalently \eqref{K1 dynamics}) possesses
\begin{enumerate}
    \item A pair of isolated, symmetric equilibria, \(E_{\pm}\), given by equation \eqref{E+- def}, which lie on the invariant line \(I\) defined in \eqref{I def}.
    The stable manifold of \(E_{+}\) is one-dimensional, and its unstable manifold is five-dimensional. 
    The unstable manifold of \(E_{-}\) is one-dimensional, and its stable manifold is five-dimensional.
    \item A three-dimensional manifold of equilibria, 
    \begin{equation}\label{V1 Def}
        \mathcal{V}_{1}
        =\Big\{r_{1}\in\R, \, \, 
        P_{1}=0, \, \,
        V_{1}=\tfrac{-(1-A\B_1(A-r_1^2))}{A\left(A-r_{1}^{2}\right)^{2}}, \, \, 
        Q_{1}\in\R, \, \, 
        \varepsilon_{1}=0, \, \, 
        \B_{1}\in\R\Big\}.
    \end{equation}
    This manifold corresponds to the entirety of the critical manifold $S_0$ with \(u<0\), expressed in chart $K_1$.
    Over each point in \(\mathcal{V}_{1}\), there exist one-dimensional stable and unstable fibers.
    There is also a one-dimensional neutrally stable fiber containing non-exponential dynamics.
    The line of equilibria \(\ell_{1}\) \eqref{l1 def} is contained in 
    \(\mathcal{V}_{1}\) and satisfies the relationship \(\ell_{1}=\mathcal{V}_{1}\rvert_{r_{1},\B_{1}=0}\).
\end{enumerate}
\end{lemma} 

\noindent
Before presenting the proof, we state a consequence of this lemma for the dynamics of \(\Gamma_{0}\) in chart $K_1$.
As established in Lemma~\ref{lem-limitsGamma0}, the limit points of \(\Gamma_{0}\), as $x_2 \to \mp \infty$, are the points \({\bf p}_{\pm}\) on the equator defined in \eqref{limitGamma0}.
These points are equilibria in $K_1$ and satisfy 
${\bf p}_{\pm}\subset\ell_{1}\subset\mathcal{V}_{1}$.
The limiting tangent vectors, given in \eqref{limitGamma0Tangent}, lie in the center eigenspace of \(\mathcal{V}_{1}\).
The asymptotic behavior can be decomposed into components tangent to 
\(\mathcal{V}_{1}\rvert_{r_{1}=0}\) and in the neutrally stable fiber over the limit points.
This tangency suggests that \(\Gamma_{0}\) must approach the limit points at a less than exponential rate.
This is confirmed in the proof of Lemma~\ref{lem-limitstangentvectors2Gamma0}, which shows \(\norm{\Gamma_{0}'}=\O(\abs{x_{2}}^{-3})\).

\bigskip 
\noindent
{\bf Proof.}
The formulas for the equilibria follow from direct computations.
Then, the linear stability properties follow from analyzing the Jacobian,
\begin{align}
    \tilde{L}=&\begin{bmatrix}
        \tfrac{-1}{2}P_{1} & \tfrac{-1}{2}r_{1} & 0 & 0 & 0 & 0\\
        \tilde{\mathbf{A}} & 3P_{1} & \tilde{\mathbf{C}} & 0 & 0 & \tilde{\mathbf{E}}\\
        0 & 2 V_{1} & 2P_{1} & \varepsilon_{1} & Q_{1} & 0\\
        \tilde{\mathbf{B}} & \tfrac{3}{2}Q_{1} & -r_{1}^{2}\varepsilon_{1}\tilde{\mathbf{C}} & \tfrac{3}{2}P_{1} & \tilde{\mathbf{D}} & -r_{1}^{2}\eps_{1}\tilde{\mathbf{E}}\\
        0 & \varepsilon_{1} & 0 & 0 & P_{1} & 0\\
        0 & 2\B_{1} & 0 & 0 & 0 & 2P_{1}
    \end{bmatrix},
\end{align}
where
\begin{align*}
    \tilde{\mathbf{A}}&=  -2 r_{1} (\B_{1}-2AV_{1}+2 r_{1}^{2}V_{1}),\\
    \tilde{\mathbf{B}}&=  - 2 r_{1} \eps_{1}(A \B_{1}-A^{2}V_{1}-\tfrac{1}{A} 
    - 2r_1^2(\B_1-2AV_1) - 3r_1^4 V_1) ,\\
    \tilde{\mathbf{C}}&=
    -\left(A-r_{1}^{2}\right)^{2},\\
    \tilde{\mathbf{D}}&= 
    -1 -r_{1}^{2}(A\B_{1}-A^{2}V_{1}-\tfrac{1}{A}-r_{1}^{2}{\bf F}_1),\\
    \tilde{\mathbf{E}}&=
    A-r_1^2,
\end{align*}
At the equilibria $E_\pm$, the Jacobian is upper triangular, and the eigenvalues are given by
$\lambda=\mp \tfrac{1}{2} \rho, 
\pm \rho,
\pm \tfrac{3}{2}\rho,
\pm 2\rho,
\pm 2\rho,
\pm 3\rho$, 
with eigenspaces 
\begin{align*}
    &\E_{\mp \rho/2}\left(E_{\pm}\right) = r_1 {\rm -axis},
    \quad
    \E_{\pm \rho}\left(E_{\pm}\right)=\vecspan{\begin{bmatrix}
        0\\
        0\\
        0\\
        \mp 2\\
        \rho\\
        0
    \end{bmatrix}
    }, \, \, \, 
    \E_{\pm 3\rho/2} \left(E_{\pm}\right) = Q_1 {\rm -axis}, \\
    &\E_{\pm 2\rho} \left(E_{\pm}\right) =\vecspan{\begin{bmatrix}
        0\\
        \pm A^{2}\\
        \rho\\
        0\\
        0\\
        0
    \end{bmatrix},
    \begin{bmatrix}
        0\\
        0\\
        0\\
        0\\
        0\\
        1
    \end{bmatrix}
    }, \quad 
    \E_{\pm 3 \rho}\left(E_{\pm}\right)=
    P_1 {\rm -axis}.
\end{align*}
Since the eigenvalues are bounded away from zero, the existence of the stable and unstable manifolds of $E_{\pm}$ stated in Part 1 of the lemma follows from standard invariant manifold theory \cite{HPS77}.

Next, at an arbitrary point on $\mathcal{V}_{1}$, the terms  \(\tilde{\mathbf{A}}\)-\(\tilde{\mathbf{E}}\) in the Jacobian $\tilde{L}$ reduce to
\begin{align*}
    \tilde{\mathbf{A}}&=  \tfrac{-2r_{1}}{A(A-r_1^2)} \left(2 - A\B_1(A-r_1^2) \right),\quad 
    \tilde{\mathbf{B}}=0, 
    \quad 
    \tilde{\mathbf{C}}=-\left(A-r_{1}^{2}\right)^{2}, \quad 
    \tilde{\mathbf{D}}= -1,
    \quad 
    \tilde{\mathbf{E}}= A-r_{1}^{2}.
\end{align*}
Also, the terms that are directly proportional to $\eps_1$ vanish at points on $\mathcal{V}_1$.
Hence, the eigenvalues at a point on \(\mathcal{V}_{1}\) are $\lambda=\pm\sqrt{\tfrac{2+r_{1}^{2}(r_{1}^{2}-A)\B_{1}}{r_{1}^{2}-A}}$ and $0$, with eigenspaces
\begin{align*}
    \E_{\lambda_{\pm}}\left(\mathcal{V}_{1}\right) &= \vecspan{\begin{bmatrix}
        \tfrac{-1}{2}r_{1}\\
        \lambda_{\pm}\\
        \tfrac{-2(1-A\B_{1}\big(A-r_{1}^{2})\big)}{A(A-r_{1}^{2})^{2}}\\
        \tfrac{3}{2}Q_{1}\\
        0\\
        2\B_{1}
    \end{bmatrix}}, \,\,\,
    \E_{0}\left(\mathcal{V}_{1}\right) = \vecspan{\begin{bmatrix}
        0\\
        0\\
        0\\
        1\\
        0\\
        0
    \end{bmatrix}, \begin{bmatrix}
        -\tilde{\mathbf{C}}-\tilde{\mathbf{E}}\\
        0\\
        \tilde{\mathbf{A}}\\
        0\\
        0\\
        \tilde{\mathbf{A}}
    \end{bmatrix},\begin{bmatrix}
        \tilde{\mathbf{C}}\\
        0\\
        -\tilde{\mathbf{A}}-\tilde{\mathbf{E}}\\
        0\\
        0\\
        \tilde{\mathbf{C}}
    \end{bmatrix}},\\
    \G_{0}^{2}\left(\mathcal{V}_{1}\right) &= \vecspan{\begin{bmatrix}
        0\\
        Q_{1}(A-r_{1}^{2})^{2}\\
        0\\
        0\\
        \tfrac{2-r_{1}^{2}(A-r_{1}^{2})\B_{1}}{A-r_{1}^{2}}   \\
        0
    \end{bmatrix}}.
\end{align*}
The eigenspace \(\E_{0}(\mathcal{V}_{1})\) is the tangent space of $\mathcal{V}_1$, and the generalized eigenspace \(\G_{0}^{2}(\mathcal{V}_{1})\) is the neutrally stable fiber over each point on it.
It should be noted that there are two special points in \(\mathcal{V}_{1}\) where the generalized eigenvector is actually a proper eigenvector.
These points are exactly the points \({\bf p}_{\pm}\).
The statements about the stable, unstable, and center manifolds of $\mathcal{V}_1$ made in Part 2 of the lemma follow from standard invariant manifold theory \cite{C1981,HPS77}.
\hfill $\Box$

\begin{remark} 
The generalized eigenvector was obtained by finding eigenvectors of $\tilde{L}^{2}$,  
\begin{align*}
    \tilde{L}^{2}=&\begin{bmatrix}
        \tfrac{-r_{1}}{2}\tilde{\mathbf{A}} & 0 & \tfrac{-r_{1}}{2}\tilde{\mathbf{C}} & 0 & 0 & \tfrac{-r_{1}}{2}\tilde{\mathbf{E}}\\
        0 & \tfrac{-r_{1}}{2}\tilde{\mathbf{A}}+2V_{1}\tilde{\mathbf{C}}+2\B_{1}\tilde{\mathbf{E}} & 0 & 0 & Q_{1}\tilde{\mathbf{C}} & 0\\
        2V_{1}\tilde{\mathbf{A}} & 0 & 2V_{1}\tilde{\mathbf{C}} & 0 & 0 & 2V_{1}\tilde{\mathbf{E}}\\
        \tfrac{3}{2}Q_{1}\tilde{\mathbf{A}} & 0 & \tfrac{3}{2}Q_{1}\tilde{\mathbf{C}} & 0 & 0 & \tfrac{3}{2}Q_{1}\tilde{\mathbf{E}}\\
        0 & 0 & 0 & 0 & 0 & 0\\
        2\B_{1}\tilde{\mathbf{A}} & 0 & 2\B_{1}\tilde{\mathbf{C}} & 0 & 0 & 2\B_{1}\tilde{\mathbf{E}}
    \end{bmatrix}.
\end{align*}
\end{remark}

\subsection{The true and faux canards persist as maximal canards} \label{subsec:persistence}

\begin{figure}[h!tbp]
    \centering
    \includegraphics[width=3.5in]{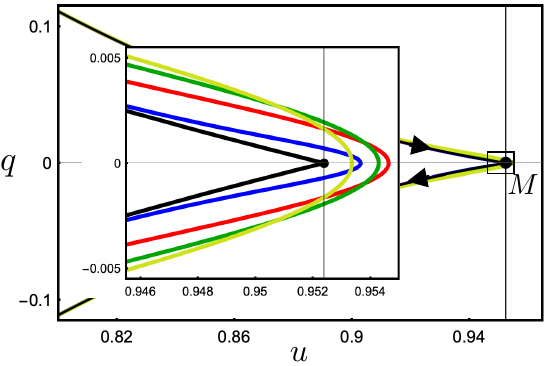}
    \caption{Persistence of the true and faux canards, and their connection, for $0< \eps \ll 1$: the singular canards (black) and the maximal canards for $\eps = 0.001$ (blue), $\eps = 0.005$ (red), $\eps=0.01$ (green), and $\eps = 0.015$ (olive). Inset: zoom on a neighborhood of the RFS point $M$. $A=1$ and $B=1.1$.}
    \label{fig:persistence}
\end{figure}

We now summarize the main results of the blow-up analysis presented in Sections~\ref{s:desingularization} and \ref{s:K1}. In the rescaling chart, we identified the key algebraic solution $\Gamma_0(x_2)$ (Lemma~\ref{lem-gamma0}), with the left-half (i.e., $x_2<0$) corresponding to the true canard $\Gamma_{\rm true}$ and the right-half (i.e., $x_2>0$) corresponding to the faux canard $\Gamma_{\rm faux}$.
This solution is forward and backward asymptotic to the equilibria ${\bf p}_{\pm} \in \ell_1 \subset \mathcal{V}_1$ on the equator of the blown-up hemisphere, i.e., in the $\left\{ \eps_1 = 0 \right\}$ subspace in the entry/exit chart $K_1$ (Lemma~\ref{lem-limitsGamma0}). 
Moreover, we have demonstrated that the special algebraic solution $\Gamma_0(x_2)$ emanates from the (attracting) 4D center manifold $W^c\left( \bf{p}_+ \right)$ and terminates in the (repelling) 4D center manifold $W^c\left( \bf{p}_- \right)$ in the $\left\{ r_1=0 \right\}$ hyperplane of chart $K_1$ (Lemma~\ref{lem-limitstangentvectors2Gamma0} and Lemma~\ref{lem-K1}).
Thus, $\Gamma_0$ lies in the 2D transverse intersection of $W^c\left( \bf{p}_+ \right)$ and $W^c\left( \bf{p}_- \right)$. We note that the 2D manifold $W^c\left( \bf{p}_+ \right) \cap W^c\left( \bf{p}_- \right)$ is foliated by the parameter $\mathcal{B}_2$ (in chart $K_2$), which is equivalent to the nontrivial foliation $r_1^4 \mathcal{B}_1 = {\rm constant}$ (in chart $K_1$). Thus, for each fixed value of $\mathcal{B}_2$ (or  $r_1^4 \mathcal{B}_1$), $\Gamma_0(x_2)$ is a 1D curve corresponding to the maximal canard solution. 
Hence, the singular canards $\Gamma_{\rm true}$ and $\Gamma_{\rm faux}$ perturb to maximal canard solutions $\Gamma_{\rm true}^{\eps}$ and $\Gamma_{\rm faux}^{\eps}$ for sufficiently small $\eps$. These persistent maximal canard solutions are shown in Fig.~\ref{fig:persistence}. 

\begin{remark}
The persistence of the connection between $\Gamma_{\rm true}$ and $\Gamma_{\rm faux}$ is a generic property of singularly perturbed spatial ODEs derived from reaction-diffusion PDEs.
Persistence is a direct consequence of transversality and the reversibility symmetry (see \eqref{reversibility}) of the system, which may be seen as follows.
First, we observe that the singular true canard $\Gamma_{\rm true}$ intersects the singular faux canard $\Gamma_{\rm faux}$  transversely in the hyperplane $\{ q=0 \}$ for all $B>1$. 
Indeed, the singular true and faux canards can be obtained as functions of $u$ from the Hamiltonian \eqref{eq:desingularizedham} of the desingularized system, i.e., by solving the following quadratic equation in $q$: 
$H_d(u,q;A,B) = H_d \left( \tfrac{2A}{1+B},0; A,B \right).$
The positive/negative root corresponds to the singular true/faux canard. Then, differentiating with respect to $u$ and evaluating at the folded singularity (where $q=0$), we find that the tangent lines to $\Gamma_{\rm true}$ and $\Gamma_{\rm faux}$ have slopes of $-\tfrac{\sqrt{B-1} (1+B)^{3/2}}{2\sqrt{2}A}$ and $\tfrac{\sqrt{B-1} (1+B)^{3/2}}{2\sqrt{2}A}$, respectively. 
Hence, for all $B>1$, the singular true canard intersects the singular faux canard transversely. 
From this, we conclude that the intersection will persist for sufficiently small $\eps$. 
Lastly, for all $\eps>0$, the reversibility symmetry (encoded in the matrix $\mathcal{R}$) enforces the constraint that $\Gamma_{\rm true}^{\eps} = \mathcal{R} \Gamma_{\rm faux}^{\eps}$, so that $\Gamma_{\rm true}^{\eps}$ and $\Gamma_{\rm faux}^{\eps}$ must intersect the $\{ q=0 \}$ hyperplane at the same point. 
\end{remark}

\section{Construction of Singular Spatially-Periodic Canards in \texorpdfstring{\eqref{spatialODE}}{Lg} }\label{s:geoconstruction}

In this section, we geometrically construct singular spatially-periodic canard solutions of \eqref{spatialODE}.
These solutions consist of alternating segments of canard orbits of the slow system on the critical manifold and pulse-shaped orbit segments of the layer problem (reduced fast system).

\subsection{Preliminary steps: Rectification of \texorpdfstring{$S^0_s$}{Lg} and locating the fast homoclinics}
\label{ss:straightenoutsaddlesheet}

It is useful to first rectify the saddle sheet of the critical manifold to the $(v,q)$ plane and to make a convenient choice of the parametrization of orbits.
In particular, to rectify $S^0_s$ we set 
\begin{equation} \label{eq:uuhat}
\hat{u} = u - u_s(v),
\end{equation} 
where we recall $u_s(v)$ and $u_c(v)$ from \eqref{usc}.
Also, since \eqref{spatialODE} is autonomous, the parametrization of orbits may be chosen freely. 
We parametrize so that the origin occurs at the symmetry point of a fast solution. 
Hence, we use $\xi$ to denote the independent variable, which is a suitable translation of the variable $y$ used in \eqref{spatialODE}, so that the leading order fast solutions have the symmetry $\xi \to -\xi$.

In the new dependent and independent variables, the spatial ODEs \eqref{spatialODE} are 
\begin{equation}\label{spatialODE-S0s-rectified}
\begin{split}
    \hat{u}_\xi &= p - \eps \frac{q u_s^2(v)}{v ( u_c(v)-u_s(v))}  \\
    p_\xi  &= v \hat{u} \left(u_c(v) - u_s(v) - \hat{u} \right)  \\
    v_\xi  &=  \eps q \\
    q_\xi  &=  \eps \left( \hat{u} + u_s(v) - v u_s(v)u_c(v) - p_y \right).
\end{split}
\end{equation}
The layer problem is now
\begin{eqnarray*}
    \hat{u}_\xi  &=& p \\
    p_\xi  &=& v \hat{u} ( u_c(v) - u_s(v) - \hat{u}),
\end{eqnarray*}
with Hamiltonian $H_f= \frac{1}{2}p^2 + \hat{u}^2 v \left( \frac{1}{3} \hat{u} - \frac{1}{2}(u_c(v)-u_s(v))\right)$. 

The saddle and center sheets of the critical manifold are now represented by
\begin{eqnarray*}
  S^0_s &= &\left\{ p=0, \, \, \hat{u}=0, \, \, 0 < v < \frac{(1+B)^2}{4A} \right\}   \\
  S^0_c &= &\left\{ p=0, \, \, \hat{u}=u_c(v)-u_s(v), \, \, 0 < v < \frac{(1+B)^2}{4A} \right\}. 
\end{eqnarray*}
All points on $S^0_s$ are on the zero level set of $H_f$.
Along this level set $p = \pm \hat{u} \sqrt{v(u_c(v)-u_s(v) - \frac{2}{3}\hat{u})}$, with $0<u_s(v) < u_c(v)$ and $v>0$.
Hence, for the homoclinics $\gamma_0(\xi)=(\hat{u}_{\rm HOM}(\xi;v_0),p_{\rm HOM}(\xi; v_0))$ that connect each point on $S^0_s$ to itself, we have the following explicit formula:
\begin{equation} \label{fasthomoclinic}
\begin{split}
    \hat{u}_{\rm HOM} (\xi; v_0) &= \frac{3}{2} (u_c(v_0)-u_s(v_0)) {\rm sech}^2 \left( \frac{1}{2} \sqrt{v_0 ( u_c(v_0)-u_s(v_0))} \, \, \xi \right) \\
    p_{\rm HOM}(\xi; v_0) &= \frac{d}{d\xi} \hat{u}_{\rm HOM}(\xi; v_0).
\end{split}
\end{equation}
The $\hat{u}$ component has the desired $\xi \to -\xi$ symmetry.

\subsection{\texorpdfstring{$W^u(S^\eps_s)$ and $W^s(S^\eps_s)$}{Lg} intersect transversely}\label{ss:transversality}

For $0<\eps \ll 1$, the splitting distance between the perturbed manifolds may be measured using the Hamiltonian. 
We calculate 
\begin{equation} \label{DeltaHf}
    \Delta H_f = \int_{-\infty}^{0} \left. \frac{dH_f}{d\xi} \right\vert_{\Gamma_\eps^u} d\xi
    +\int_{0}^\infty \left. \frac{dH_f}{d\xi} \right\vert_{\Gamma_\eps^s} d\xi,
\end{equation}
where $\Gamma_\eps^{u,s}$ denote solutions on $W^{s,u}_{\rm loc}(S^\eps_s)$.
A lengthy but straightforward calculation reveals that 
\begin{equation*}
    \frac{dH_f}{d\xi} = \eps q \hat{u} \left( \frac{1}{3} \hat{u}^2 + \hat{u} u_s(v) + u_s^2(v)\right).
\end{equation*}
The solutions on the local perturbed manifolds are represented by asymptotic expansions 
$\gamma_\eps^{u,s} (\xi) =
\gamma_0(\xi) + \eps \gamma_1^{u,s}(\xi) + \mathcal{O}(\eps^2),$
on the semi-infinite intervals, $\xi \in (-\infty,0)$ and $(0,\infty)$, respectively.

We substitute the perturbation expansions into the vector field \eqref{spatialODE-S0s-rectified}.
At $\mathcal{O}(1)$ (leading order), we recover the symmetric homoclinic orbit $\gamma_0$  \eqref{fasthomoclinic}.
Next, at $\mathcal{O}(\eps)$, the equations are
\begin{eqnarray*}
    \hat{u}_{1\xi} &=& p_1 - \frac{q_0 u_s^2(v_0)}{v_0 ( u_c (v_0-u_s(v_0))} \nonumber \\
    p_{1\xi}  &=& (\hat{u}_1 v_0 + \hat{u}_0 v_1) (u_c(v_0)-u_s(v_0) - \hat{u}_0) + \hat{u}_0 v_0 \left( v_1 \frac{d(u_c-u_s)}{dv_0} - \hat{u}_1 \right)\nonumber \\
    v_{1\xi}  &=& q_0 \nonumber \\
    q_{1\xi}  &=& \hat{u}_0 + u_s(v_0) - v_0 u_c(v_0) u_s(v_0) - \frac{dp_0}{d\xi},
\end{eqnarray*}
where $\hat{u}_0 = \hat{u}_{\rm HOM}$ and $p_0 = p_{\rm HOM}$. 
This system may be solved in closed form.
We see directly that $v_1(\xi)=v_1(0) + q_0 \xi$.
Also, $q_1(\xi) = -p_0(\xi) + (u_s(v_0)-v_0 u_c(v_0) u_s(v_0))\xi + \int_0^\xi u_0(s) ds.$
Then, the $(\hat{u}_1,p_1)$ subsystem is described by the linear inhomogeneous system
\begin{eqnarray*}
    \begin{bmatrix}
        \hat{u}_{1\xi} \\
        p_{1\xi} 
    \end{bmatrix}
    =
    \begin{bmatrix}
        0 & 1 \\
        v_0(u_c(v_0)-u_s(v_0) -2 {\hat u}_0) & 0
    \end{bmatrix}
    \begin{bmatrix}
        \hat{u}_1 \\
        p_1
    \end{bmatrix}
    + \begin{bmatrix}
        \frac{-q_0 u_s^2}{v_0(u_c - u_s)} \\
        -\hat{u}_0 v_1 \left( \frac{2u_c u_s}{u_c - u_s} + \hat{u}_0 \right)
    \end{bmatrix}
\end{eqnarray*}
The solution is
\begin{equation}\label{u1p1solution}
    \begin{bmatrix}
        \hat{u}_1 (\xi) \\
        p_1 (\xi)
    \end{bmatrix}
    =
    U(\xi) U^{-1}(0) 
    \begin{bmatrix}
        \hat{u}_1(0) \\
        p_1(0)
    \end{bmatrix}
    + U(\xi) \int_0^\xi U^{-1}(s) F(s) ds,
\end{equation}
where (in this section) $F$ denotes the vector representing the inhomogeneous term and $U$ denotes the fundamental matrix of the homogeneous problem,
\begin{equation*}
    U(\xi) = \begin{bmatrix}
        u_{11}(\xi) & u_{12}(\xi) \\
        u_{11}'(\xi) & u_{12}'(\xi)
    \end{bmatrix},
\end{equation*}
and $u_{11}(\xi)= \frac{d\hat{u}_{\rm HOM}}{d\xi}$
and $u_{12}(\xi) = u_{11}(\xi) \int_0^\xi u_{11}^{-2}(s) ds$.
The integral may be found in closed form.

Substituting the terms of the expansion into the splitting distance measurement, we find
\begin{equation}\label{DeltaHf=}
\Delta H_f = \eps D_1 + \eps^2 D_2 + \mathcal{O}(\eps^3),
\end{equation}
where 
\begin{equation}    \label{D1D2}
\begin{split}
    D_1 &= \int_{-\infty}^\infty
         q_0 \hat{u}_0 \left( \tfrac{1}{3} \hat{u}_0^2 + \hat{u}_0 u_s + u_s^2 \right) d\xi \nonumber \\
    D_2 &= \int_{-\infty}^\infty 
         \left[ 
         (q_1 \hat{u}_0 + q_0 \hat{u}_1 ) \left(  \tfrac{1}{3} \hat{u}_0^2 + \hat{u}_0 u_s + u_s^2 \right)
         + q_0 \hat{u}_0\! \left( \hat{u}_1 \left( \tfrac{2}{3} \hat{u}_0 + u_s \right)
         + v_1 (\hat{u}_0 + 2 u_s) \tfrac{du_s}{dv_0} \right)
         \right] d\xi,
\end{split}
\end{equation}
where the $v_0$-dependence in $u_s$ and $u_c$ has  been suppressed. 
We find
\begin{equation} \label{D1-evaluated}
    D_1 = \frac{6 q_0}{5 v_0} (2u_c^2 + u_c u_s + 2 u_s^2 ) 
    \sqrt{v_0(u_c-u_s )}.
\end{equation}
The factor of $q_0$ in the numerator of $D_1$ determines the sign of $D_1$, because the rest of the function is positive.
Also, $D_2$ is a positive multiple of $q_1(0)$.
Hence, we find $\Delta H_f=0$ to within $\mathcal{O}(\eps^3)$ for $q_0=0$ and $q_1(0)=0$.
Therefore, $W^u(S^\eps_s)$ and $W^s(S^\eps_s)$ intersect for these values.
Furthermore, the intersection is transverse, because the zeroes of $D_1$ and $D_2$ are simple zeroes.

\begin{figure}[h!tbp]
\centering
\includegraphics[width=4in]{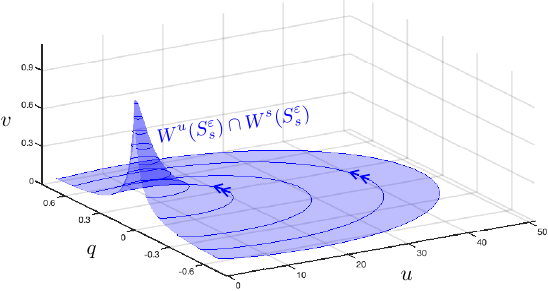}
\caption{
A plot in the $(u,q,v)$ space showing the two-dimensional (blue) surface consisting of the homoclinic orbits in the transverse intersection 
$W^u(S^\eps_s) \cap W^s(S^\eps_s)$.
The homoclinic orbits (blue curves) have been computed using numerical continuation (see App.~\ref{app:homman}).
Here, $A=1$, $B=1.03$, and $\eps=0.01$.
}
\label{fig:homoclinicmanifold}
\end{figure}

The union of the homoclinic orbits in $W^u(S^\eps_s) \cap W^s(S^\eps_s)$ forms a 2-D surface
(see Fig.~\ref{fig:homoclinicmanifold}).
It is symmetric in $\{ p=0, q=0\}$, because by the symmetry $\xi \to -\xi$ the reflection maps $W^s(S^\eps_s)$ to $W^u(S^\eps_s)$ and vice versa.
The surface intersects the hyperplane $\{ p=0, q=0 \}$ transversely.

\subsection{The takeoff and touchdown curves}
\label{ss:takeofftouchdown}

As just shown,
the transversal intersection $W^u(S^\eps_s) \cap W^s(S^\eps_s)$ is two-dimensional.
The fast unstable and stable fibers in this intersection have base points on $S^\eps_s$.
The union of the base points of the fast unstable fibers in this intersection defines the takeoff curve $T_{\rm takeoff}$, and the union of the base points of the fast stable fibers in this intersection defines the touchdown curve $T_{\rm touchdown}$ on $S^\eps_s$.

\begin{figure}[t!]
\centering
\includegraphics[width=5in]{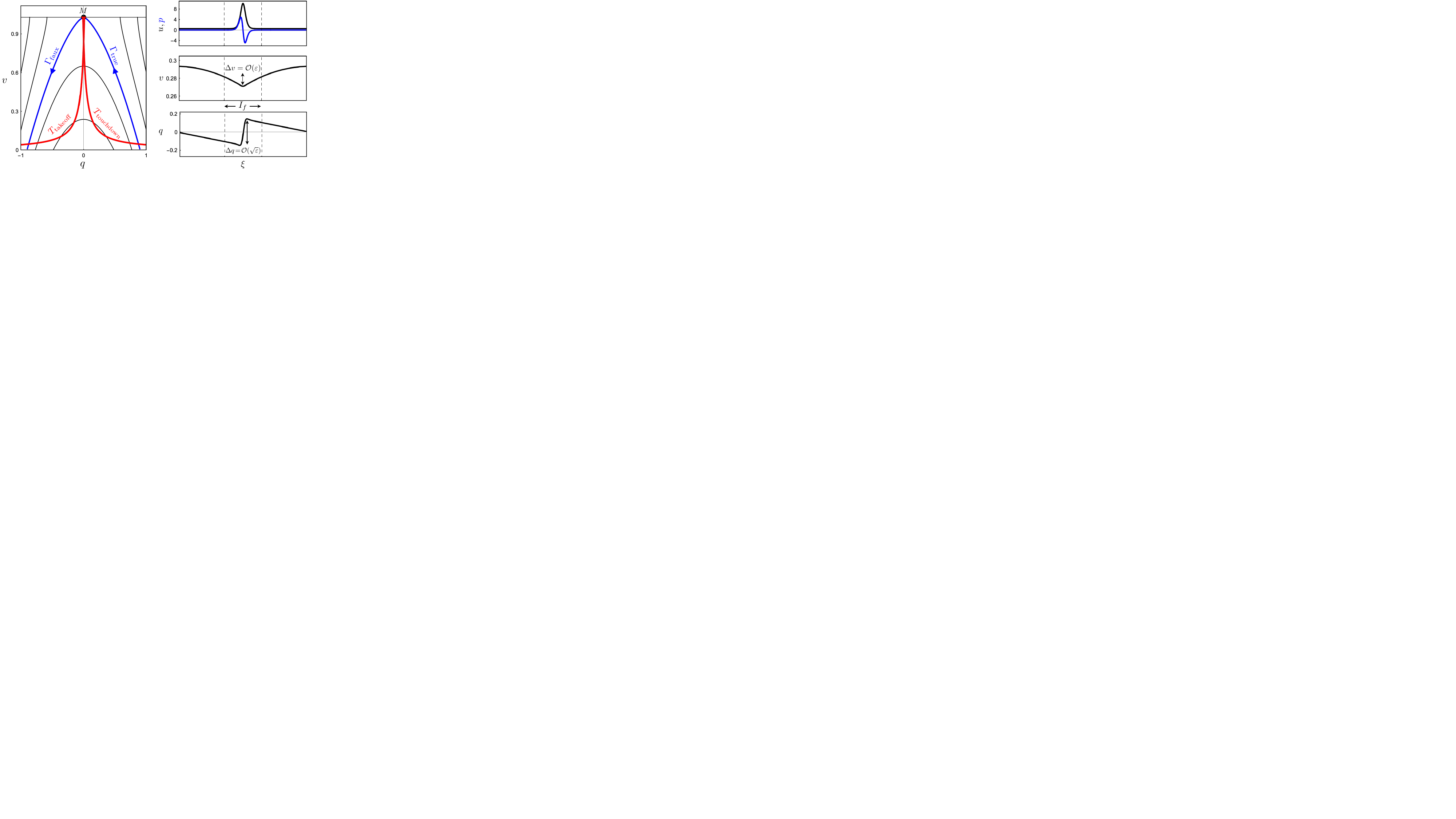}
\put(-362,192){(a)}
\put(-180,192){(b)}
\put(-180,128){(c)}
\put(-180,62){(d)}
\caption{
The takeoff and touchdown curves $T_{\rm takeoff}$ and $T_{\rm touchdown}$, \eqref{takeoff+touchdown-T}, for $A=1$, $B=1.03$, and $\eps=0.01$. 
(a) Projections into the $(q,v)$ plane of $T_{\rm takeoff}$ and $T_{\rm touchdown}$, which lie on $S_s^0$. 
We also show the desingularized reduced flow (black curves) on $S_s^0$ and $\Gamma_{\rm true}$ and $\Gamma_{\rm faux}$ of $M_{\rm RFS}$ (black curves). (b)--(d) Profiles of a spatially-periodic solution to highlight the dynamics along the spatial interval $I_f$ (between the vertical dashed lines), where the fast jump occurs. 
The fast pulse in $u$ and its spatial derivative $p$ are shown in (b) in black and blue, respectively.
(c) shows the small excursion in $v$ that occurs on $I_f$, and (d) shows the jump $\Delta q$ in $q$, recall \eqref{fastjumpinq}.
That jump bridges the gap between $T_{\rm takeoff}$ and $T_{\rm touchdown}$ at the height of $v \approx 0.29$.
Also, outside of $I_f$, on the complementary $\xi$-intervals, the solution drifts slowly along $S_s^{\eps}$ near the (black) orbit that connects $T_{\rm touchdown}$ back to $T_{\rm takeoff}$. 
(See also Fig.~\ref{fig:pulseconstruction}(b).)
}
\label{fig:takeoff+touchdown}
\end{figure}
To determine the locations of the takeoff and touchdown curves, we calculate the jump in $q$ that occurs as orbits make the excursion from --and back to-- $S^\eps_s$, {\it i.e.,} the distance between the base points of the fast unstable and stable fibers (see Fig.~\ref{fig:takeoff+touchdown}, left column).
Let 
\[ I_f = \left\{ \xi \in \mathbb{R}: -\frac{1}{\sqrt{\eps}\, \lambda(v_0)} < \xi < \frac{1}{\sqrt{\eps} \, \lambda(v_0)} \right\} \]
be the interval over which the fast jump occurs, where $\pm \lambda(v_0) = \pm \sqrt{v_0 \left( u_c(v_0)-u_s(v_0) \right)}$ corresponds to the eigenvalues of the layer problem along the saddle sheet, $S_s^0$, of the critical manifold. 

The interval, $I_f$, has been chosen so that the $u$ and $p$ components of the solution are $\mathcal{O} \left( \exp \left( -\tfrac{1}{\sqrt{\eps}} \right)\right)$ at the ends of the interval (see Fig.~\ref{fig:takeoff+touchdown}(a)). 
Moreover, $I_f$ is asymptotically large in the $\xi$ variable and asymptotically small in the $\eps \, \xi$ variable.
The jump in $q$, then, is given to leading order by
\begin{equation}\label{fastjumpinq}
    \Delta q = \int_{I_f} \frac{dq}{d\xi} d\xi 
    = \int_{I_f} \left( \frac{dq_0}{d\xi} + \eps \frac{dq_1}{d\xi} + \mathcal{O}(\eps^2) \right) d\xi,
\end{equation}
and we note that the tails for $\xi \in (-\infty,-\frac{1}{\sqrt{\eps}\, \lambda(v_0)})$ and $\xi \in (\frac{1}{\sqrt{\eps}\, \lambda(v_0)},\infty)$ contribute only at higher order.
Substituting in the formula \eqref{fasthomoclinic}, we find
\begin{equation}
 \Delta q = 2 \mathbf{T}(v_0; \eps),
 \end{equation}
 where
 \begin{equation}\label{T}
 \begin{split}
     \mathbf{T}(v_0;\eps) = \sqrt{\eps} \frac{u_s(v_0)-A}{\lambda(v_0)} + \eps \frac{3\lambda(v_0)}{v_0} \tanh \left( \frac{1}{2\sqrt{\eps}} \right) \left( 1+ \frac{1}{2} \lambda^2(v_0) \operatorname{sech}^2 \left( \frac{1}{2\sqrt{\eps}} \right) \right)+\mathcal{O}\left( \eps^{3/2} \right).
\end{split}
\end{equation}
Therefore, the takeoff and touchdown curves are located to leading order at
\begin{equation}\label{takeoff+touchdown-T}
\begin{split}
    T_{\rm takeoff} &= \{ (u,p,v,q)\in S^\eps_s: \,\, q = - \mathbf{T}(v;\eps) \}, \\
    T_{\rm touchdown} &= \{ (u,p,v,q) \in S^\eps_s: \, \, q = + \mathbf{T}(v;\eps) \}.
\end{split}
\end{equation}
They are shown (red curves) in the $(q,v)$ plane in Fig.~\ref{fig:takeoff+touchdown}. 
We note that the coefficients $\tfrac{u_s(v_0)-A}{\lambda(v_0)}$ and $\frac{3\lambda(v_0)}{v_0}$ in \eqref{T}
can be $\mathcal{O}(1), \mathcal{O}(\sqrt{\eps})$, or $\mathcal{O}\left( \tfrac{1}{\sqrt{\eps}} \right)$ depending on the proximity of $B$ to $1$ and the proximity of $v_0$ to $0$ or $\tfrac{B}{A}$, as shown in Fig.~\ref{fig:homoclinicmanifold}.

\begin{remark}
The change in $v$ during the fast jump is given to leading order by 
\[ \Delta v = \int_{I_f} \frac{dv}{d\xi} \, d\xi = \sqrt{\eps} \frac{2q_0}{\lambda(v_0)} + \mathcal{O}\left( \eps^{3/2} \right). \]
Since $\Delta q = \mathcal{O}\left( \sqrt{\eps} \right)$ during a fast jump, it follows that $q_0 = \mathcal{O}\left( \sqrt{\eps} \right)$. 
Thus, $\Delta v = \mathcal{O}(\eps)$. 
\end{remark}

\begin{figure}[h!tbp]
    \centering
    \includegraphics[width=\linewidth]{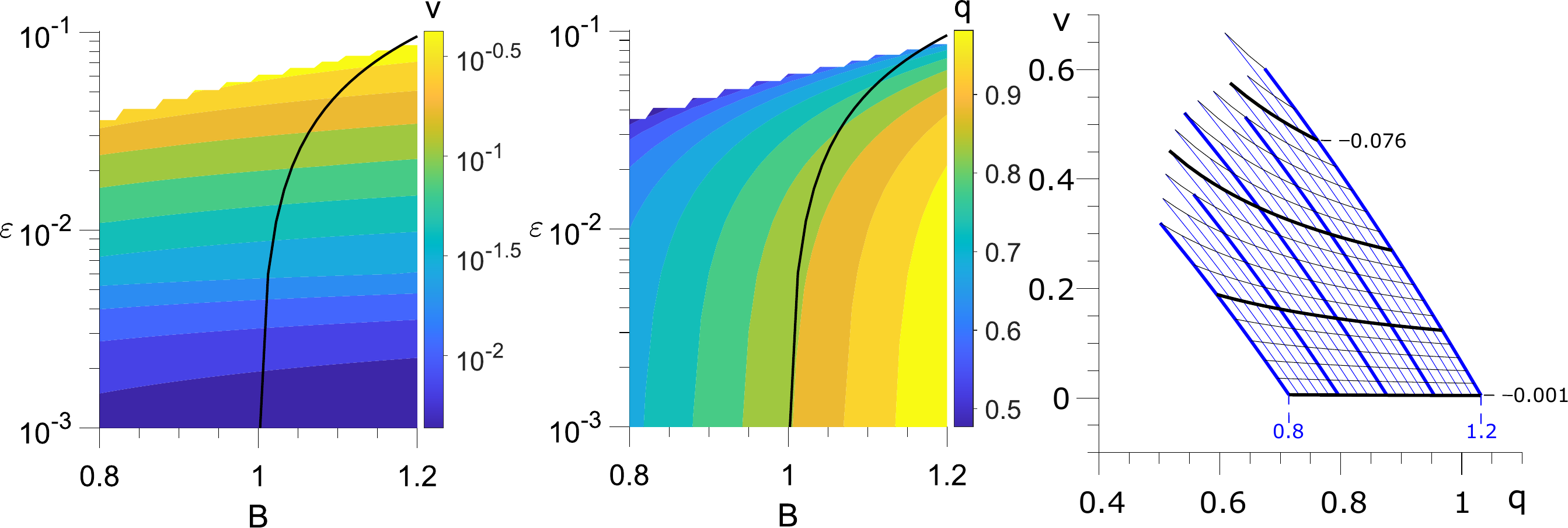}
    \put(-430,155){(a)}
    \put(-270,155){(b)}
    \put(-130,155){(c)}
    \caption{
    (a) The $v$ coordinate (log scale) of the point $T_{\rm touchdown} \cap \Gamma_{\rm true}$ as a function of $(B,\eps)$.
     The \(\varepsilon\) axis is logarithmically scaled, and the black curve indicates the location of the Turing bifurcation.
    (b) The $q$ coordinate (linear scale) of the intersection point.
    (c) The surface of intersection points in the $(q,v)$-plane, ruled with curves of constant \(B\) and \(\varepsilon\).
    Each blue curve (more vertical) increments the \(B\) value by \(0.02\) from \(0.8\) to \(1.2\).
    Each black curve (more horizontal) increments the \(\varepsilon\) values by \(0.005\) from \(0.001\) to \(0.086\).
    For sample points with \(0.091\leq \varepsilon\leq 0.101\), Newton's method did not converge for any test \(B\) value. Here $A=1$.
    }
    \label{fig:placeholder}
\end{figure}

\subsection{Geometric construction of a family of singular spatially-periodic canards}\label{ss:geoconstruction}

\begin{figure}[h!tbp]
\centering
\includegraphics[width=5in]{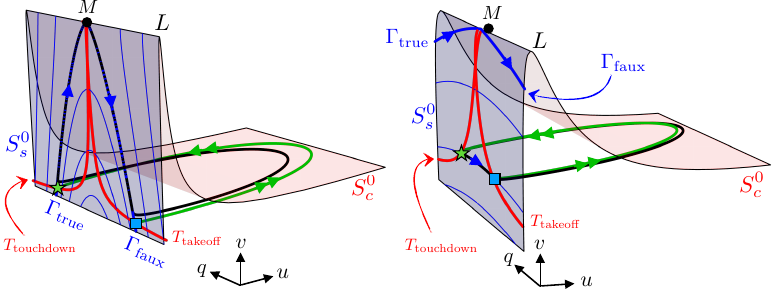}
\put(-364,132){(a)}
\put(-182,132){(b)}
\caption{Pulse construction in the original $(u,q,v)$ phase space, where $u = \hat{u}+u_s(v)$ by \eqref{eq:uuhat}.
(a) A spatially-periodic solution (black curve) with asymptotically large spatial period ($T=1000$) involving canard segments. 
(b) A spatially-periodic solution (black curve) with $T=125$ that does not have canard segments. 
The saddle and center sheets, $S_s^0$ and $S_c^0$ (blue and red surfaces) of the critical manifold meet along the fold set $L$, which contains the RFS point $M$ (black marker). 
$\Gamma_{\rm true}$ and $\Gamma_{\rm faux}$ (thick blue curves) of the RFS point partition the reduced flow (thin blue curves) on $S_s^0$. 
$T_{\rm takeoff}$ and $T_{\rm touchdown}$ are the red curves. 
The homoclinic orbit of the layer problem (green) lies in a plane of constant $v$. 
In (a), the value of $v$ is determined to leading order by where $T_{\rm takeoff} \cap \Gamma_{\rm faux}$ and where $T_{\rm touchdown} \cap \Gamma_{\rm true}$. 
This $v$ is small, and hence the pulse has asymptotically large amplitude (here $u \in [0,70]$). 
Then, for each fixed $v$ above this value, there is a homoclinic orbit connecting $T_{\rm takeoff}$ to $T_{\rm touchdown}$.
The green curve in (b) is one example: $v$ is $\mathcal{O}(1)$, and the pulse has moderate height ($u \in [0,12]$ is shown).
Here, $A=1$, $B=1.03$, and $\eps=0.01$.}
\label{fig:pulseconstruction}
\end{figure}

In this subsection, we construct a family of singular spatially-periodic canards, beginning with those that have asymptotically large spatial periods.
These have three segments:
\begin{itemize}
\item Starting at the RFSN-II point, $M$, the first segment is a slow orbit segment on $S^0_s$ that is close to $\Gamma_{\rm faux}$ of the RFSN point. 
This slow orbit segment terminates at the takeoff point (Fig.~\ref{fig:pulseconstruction}(a), blue square marker), defined by $\Gamma_{\rm faux} \cap T_{\rm takeoff}$, \eqref{takeoff+touchdown-T}. 
\item The second segment is the fast orbit segment (green curve) that connects the takeoff point (blue square marker) to the touchdown point (Fig.~\ref{fig:pulseconstruction}(a), green star marker), defined by $\Gamma_{\rm true} \cap T_{\rm touchdown}$. 
This fast orbit segment lies in $W^u(S_s^{\eps}) \cap W^s (S_s^{\eps})$.
To compute it, we use the asymptotic expansions $\gamma_{\eps}^{u,s}(\xi) = \gamma_0(\xi) + \eps \gamma_1^{u,s}(\xi) + \mathcal{O}(\eps^2)$ to determine the two halves of the orbit. For the orbit $\gamma_{\eps}^{u}(\xi) \in W^u(S_s^{\eps})$, we choose the following endpoint data at $\xi = -\tfrac{1}{\sqrt{\eps} \lambda(v_{\rm totd})}$:
\begin{equation*}
\begin{split}
    \hat{u}_0 &= \hat{u}_{\rm HOM}\left(-\tfrac{1}{\sqrt{\eps} \lambda(v_{\rm totd})};v_{\rm totd}\right),\,\, p_0 = p_{\rm HOM}\left(-\tfrac{1}{\sqrt{\eps} \lambda(v_{\rm totd})};v_{\rm totd}\right), \,\, v_0 = v_{\rm totd}, \,\, q_0 = 0, \\
    \hat{u}_1^u &= 0, \quad p_1^u = 0, \quad v_1^u = 0, \quad q_1^u = -\mathbf{T}(v_{\rm totd};\eps),
\end{split}
\end{equation*} 
where $v_{\rm totd}$ is the $v$-value of the intersection between the faux (true) canard and the takeoff (touchdown) curve (here `totd' denotes takeoff and touchdown).
Similarly, the orbit $\gamma_{\eps}^{s}(\xi) \in W^s(S_s^{\eps})$ can be constructed by choosing the following endpoint data at $\xi = \tfrac{1}{\sqrt{\eps} \lambda(v_{\rm totd})}$:
\begin{equation*}
\begin{split}
    \hat{u}_0 &= \hat{u}_{\rm HOM}\left(\tfrac{1}{\sqrt{\eps} \lambda(v_{\rm totd})};v_{\rm totd}\right), \,\, p_0 = p_{\rm HOM}\left(\tfrac{1}{\sqrt{\eps} \lambda(v_{\rm totd})};v_{\rm totd}\right), \,\, v_0 = v_{\rm totd}, \,\, q_0 = 0, \\
    \hat{u}_1^s &= 0, \quad p_1^s = 0, \quad v_1^s = 0, \quad q_1^s = \mathbf{T}(v_{\rm totd};\eps).
\end{split}
\end{equation*} 
Combining $\gamma_{\eps}^{u}(\xi)$ and $\gamma_{\eps}^{s}(\xi)$, we obtain the following parametrization of the fast orbit segment:
\begin{equation}    \label{eq:fastpulse}
\begin{split}
    \hat{u}(\xi) &= \hat{u}_{\rm HOM}(\xi;v_{\rm totd}) + \mathcal{O} \left( \eps^2 \right), \\
    p(\xi) &= p_{\rm HOM}(\xi;v_{\rm totd}) + \mathcal{O} \left( \eps^2 \right), \\
    v(\xi) &= v_{\rm totd} + \mathcal{O} \left( \eps^2 \right) \\ 
    q(\xi) &= \eps \left( (u_s(v_{\rm totd})-A) \xi + \tfrac{3\lambda(v_{\rm totd})}{v_{\rm totd}} \tanh \left( \tfrac{1}{2} \lambda(v_{\rm totd}) \xi \right) - p_{\rm HOM}(\xi;v_{\rm totd}) \right) + \mathcal{O} \left( \eps^2 \right),
\end{split}
\end{equation}
for $\xi \in I_f$.
\item The third segment is a slow orbit segment on $S^0_s$ that is close to the true canard, $\Gamma_{\rm true}$, from the touchdown point (green star marker) back up to the RFS point at $M$.
\end{itemize}

The key ingredients are:
(i) $W^u(S^\eps_s)\cap W^s(S^\eps_s)$ is transverse (established in Subsec.~\ref{ss:transversality}),
(ii) $\Gamma_{\rm faux} \cap T_{\rm takeoff}$ is transverse (recall Subsec.~\ref{ss:takeofftouchdown} and Figs.~\ref{fig:pulseconstruction}(a) and (b)), and
(iii) $T_{\rm touchdown} \cap \Gamma_{\rm true}$ is transverse.
Hence, standard theory may be used to show that, for $0<\eps\ll 1$, there exists a periodic orbit that is exponentially close to $S^\eps_s$ along almost all of the slow segment near $\Gamma_{\rm ture}$ and $\Gamma_{\rm, faux}$, and that lies close to $W^u(S^\eps_s)\cap W^s(S^\eps_s)$ along the fast segment.

Within this family of canards, those with asymptotically large periods continue into orbits with shorter periods
(see Fig.~\ref{fig:pulseconstruction}(b)). 
The orbits consist of two segments: a slow orbit segment on $S^0_s$ that connects points on $T_{\rm takeoff}$ and $T_{\rm touchdown}$ at the same height in $v$ and that lies between $\Gamma_{\rm true}$ and $\Gamma_{\rm faux}$,
and a fast orbit segment \eqref{eq:fastpulse} that connects these points on $T_{\rm takeoff}$ and $T_{\rm touchdown}$.

We observe that, while the geometric construction has been performed here for $B$ close to $B_T$, and hence the RFSN-II point, $\Gamma_{\rm true}$, and $\Gamma_{\rm faux}$ create the spatial canards, a similar construction shows that spatially-periodic canard orbits exist for each $B > B_T$, with $|B-B_T| = \mathcal{O}(1)$. 
The flow on the saddle sheet of the critical manifold is similar, and there is a general RFS point at $u_M=\tfrac{2A}{B+1}$, which has true and faux canards.
Also, for each $B\ne 1$, there is a two-dimensional surface of orbits homoclinic to $S^0_s$ with similar take-off and touchdown curves.
With these structural ingredients, the spatially-periodic pulses can then be constructed in the same manner.

\subsection{Families of spatially-periodic canards}
\label{ss:families}

\begin{figure}[h!t]
    \centering
    \includegraphics[width=5in]{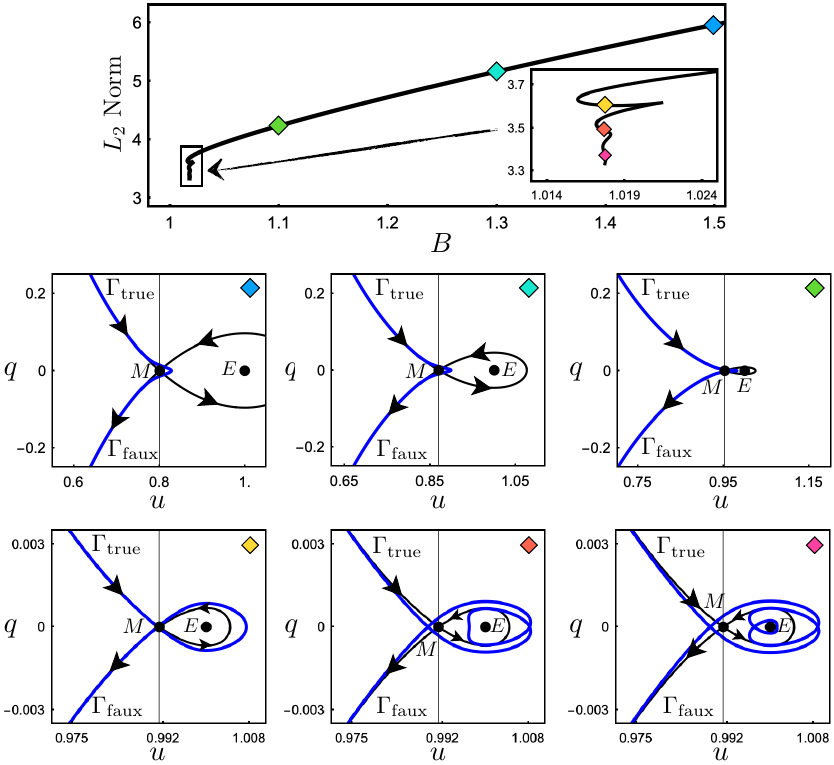}
    \put(-320,322){(a)}
    \put(-361,215){(b)}
    \put(-238,215){(c)}
    \put(-114,215){(d)}
    \put(-361,104){(e)}
    \put(-238,104){(f)}
    \put(-114,104){(g)}
    \caption{(a) Bifurcation diagram of a family of spatially-periodic canards.
    (b)-(g) Projections into the $(u,q)$ plane, into small boxes of different sizes centered on $M_{\rm RFS}$ (black marker), of the solutions (blue curves) for each colored diamond.
    (b) $B=1.5$, (c) $B=1.3$, (d) $B=1.1$, (e) $B=1.01775$ and $T\approx 767.12$, (f) $B=1.01775$ and $T\approx 830.84$, and (g) $B=1.01775$ and $T \approx 898.06$. 
    $\Gamma_{\rm true}$ and $\Gamma_{\rm faux}$ are the black curves.
    The directions of the arrows on $\Gamma_{\rm true}$ and $\Gamma_{\rm faux}$ are from the reduced system \eqref{onS0}, not the desingularized system \eqref{desingularizedreduced}. 
    Hence, the directions along $\Gamma_{\rm true}$ and $\Gamma_{\rm faux}$ on $S^0_c$ are the natural directions
    ({\it i.e.}, the orientation has not been reversed, as it is on $S^0_c$ for \eqref{desingularizedreduced}).
    The inset in (a) shows that, near the Turing bifurcation at $B_T = (1+\eps A)^2=1.0201$, the curve exhibits almost self-similar structure. 
     Here, $A=1$ and $\eps = 0.01$.}
    \label{fig:maxcanards}
\end{figure}
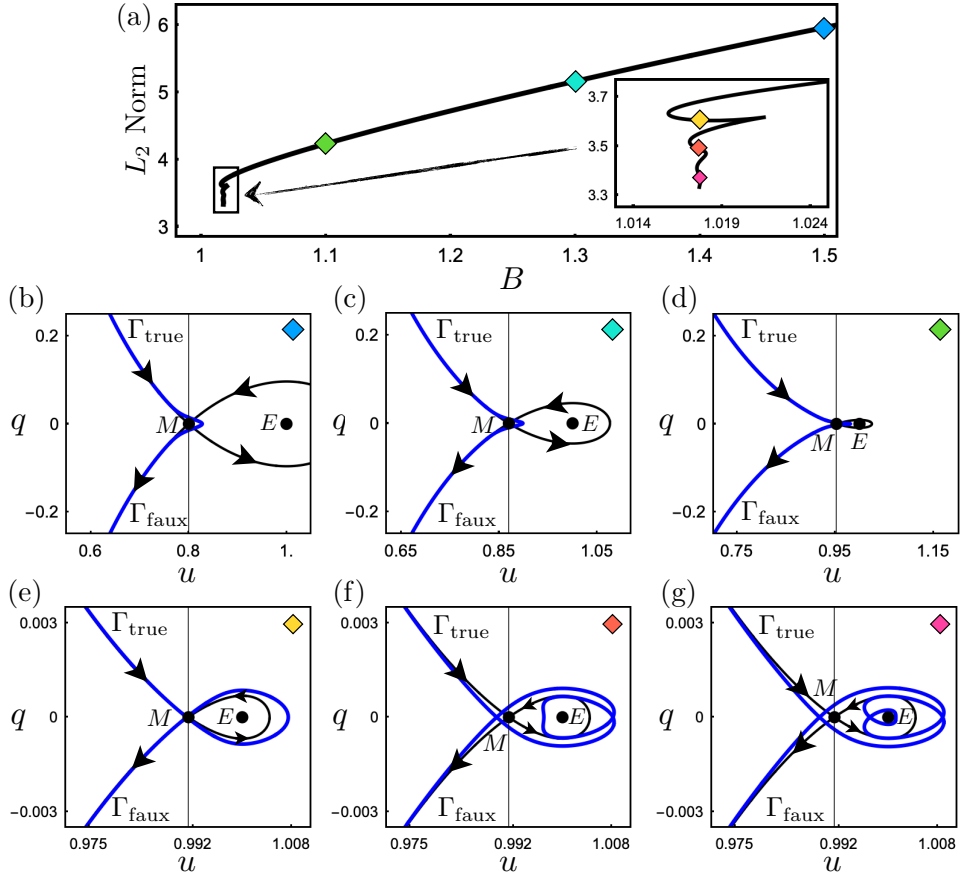

Having shown how the spatially-periodic canards are constructed in Subsec.~\ref{ss:geoconstruction}, we present a typical branch of maximal canards for a sequence of decreasing $B$ values (see Fig.~\ref{fig:maxcanards}(a)). 

For $B$ sufficiently far to the right of $B=1$, the RFS point at $u_M=\tfrac{2A}{B+1}$ lies sufficiently far away from (and to the left of) the equilibrium $E$. 
As shown in Figs.~\ref{fig:maxcanards}(b)-(d), corresponding to the blue, cyan, and green diamonds in Fig.~\ref{fig:maxcanards}(a), the orbits have long segments where they gradually vary in space along the true canard $\Gamma_{\rm true}$ until they reach the level of $u_M$. 
Then, they exhibit a small-amplitude pulse (of which the figure only shows the projection), with $u$ reaching a local maximum just to the right of $u_M$.
Subsequently, they have a long slow segment near $\Gamma_{\rm faux}$, which is symmetric about the $u$-axis to the segment along $\Gamma_{\rm true}$.
Lastly, the connection back up to the neighborhood of $\Gamma_{\rm true}$ consists of a fast pulse of somewhat larger amplitude (not shown).

Then, as we take successively smaller values of $B$ along this branch, the solutions continue to exhibit the same type of long slow segments near $\Gamma_{\rm true}$ and $\Gamma_{\rm faux}$, and they are connected back up (not shown) with fast pulses in which $u$ reaches its maximum.
However, there is a qualitative difference in the short, fast segments near $M$, between the two slow segments.
Namely, for $B$ sufficiently close to $B=1$, the RFS point $M$ is closer to the equilibrium $E$, and in between the long slow segments near $\Gamma_{\rm true}$ and $\Gamma_{\rm faux}$ the solutions exhibit one (or more) small-amplitude, fast oscillations about $E$ (see Figs.~\ref{fig:maxcanards}(e)-(g) for the solutions at the yellow, red, and magenta diamonds).
Also, the bifurcation curve shown in Fig.~\ref{fig:maxcanards}(a) exhibits a certain self-similarity about a critical value at $B \approx 1.01775$, with branches of solutions that are small in the $L^2$ norm (see Sec.~\ref{s:selfsimilarity}).

\subsection{Spatially-periodic canard-mediated bursting} \label{subsec:bursting}

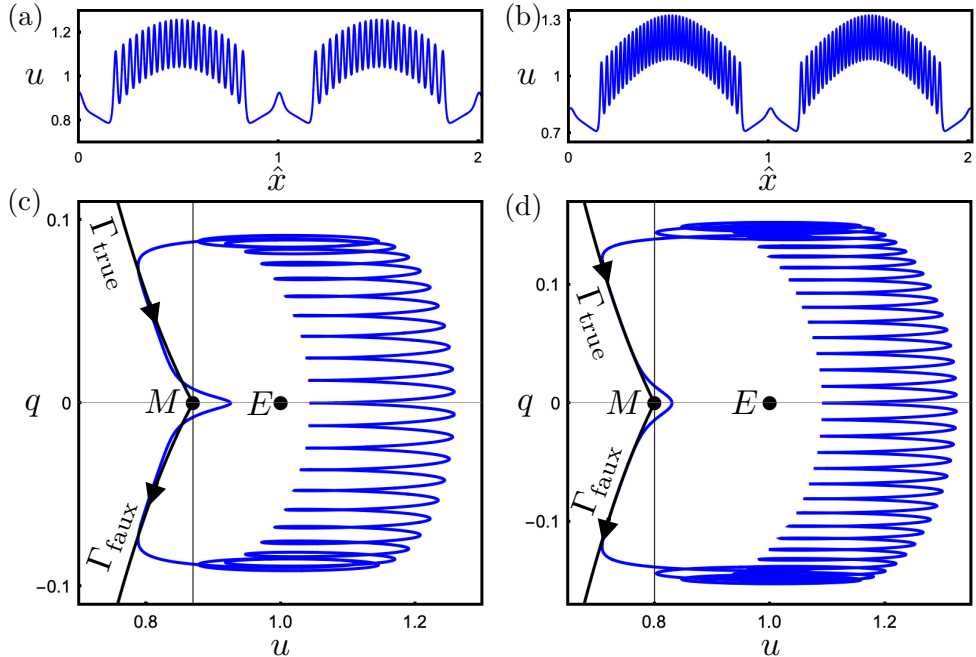
\begin{figure}[h!t]
    \centering
    \includegraphics[width=5in]{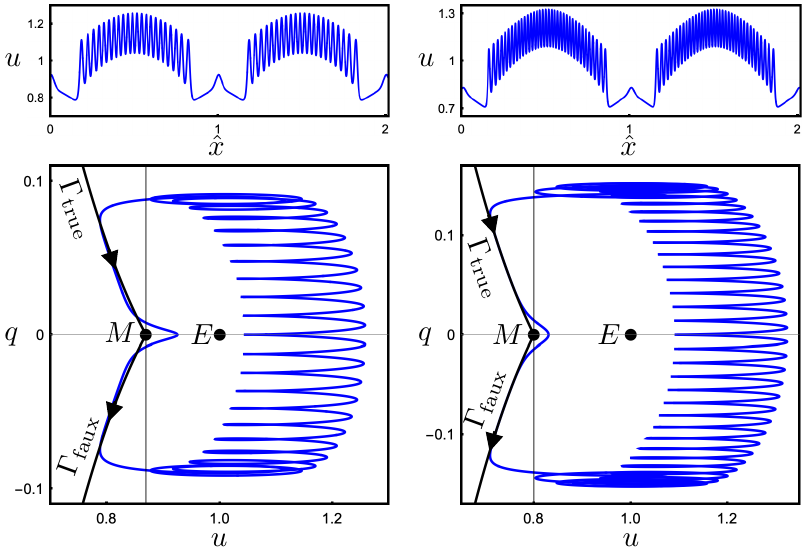}
    \put(-366,240){(a)}
    \put(-366,170){(c)}
    \put(-178,240){(b)}
    \put(-178,168){(d)}
    \caption{Spatially-periodic canard-mediated bursting solutions for $A=1$ and $\eps=0.01$. 
    Left column: $B = 1.3$; right: $B=1.5$. 
    (a) $u(\hat{x})$ with $22$ fast spikes, shown over two periods $T$, with $T=310$. 
    (b) $u(\hat{x})$ with $36$ fast spikes and $T=358$. 
    (c) and (d) Projections of (a) and (b) into the $(u,q)$ plane. 
 The slow segments closely follow $\Gamma_{\rm true}$ and $\Gamma_{\rm faux}$.
 The fast segments consist of spiking oscillations that have gradually increasing spatial frequency (up to the profile center) and then a symmetrical, gradually decreasing frequency. 
    }
    \label{fig:bursting}
\end{figure}

We also find ``bursting'' solutions, which consist of intervals of slow (or gradual) spatial variation in alternation with intervals on which the solutions exhibit a finite sequence of spatial oscillations in rapid succession. 
Two representatives are shown in Fig.~\ref{fig:bursting}. 

The slow `silent' interval corresponds to slow drift along $\Gamma_{\rm true}$ and $\Gamma_{\rm faux}$.
It terminates at a takeoff point along $\Gamma_{\rm faux}$, where the `active' bursting begins. 
On the burst interval, the solution has a finite number of high-frequency spatial oscillations that, on average, slowly drift from $q<0$ to $q>0$. 
The burst interval ends at a touchdown point on $\Gamma_{\rm true}$, and the solution is then back in the silent interval over the rest of the period. 
These spatial canard bursting solutions exist robustly in the Brusselator. 
As $B$ is increased, the number of spikes also increases, see Table~\ref{tab:spikeno}.

\begin{table}[h!t]
    \centering
    \begin{tabular}{ccccccc}
        \hline
        $B$ & $1.01$ & $1.1$ & $1.2$ & $1.3$ & $1.4$ & 1.5 \\
        \hline
        $n_{\rm spikes}$ & $7$ & $11$ & $16$ & $22$ & $29$ & $36$ \\
        \hline
    \end{tabular}
    \caption{Number of spikes during the active phase of a spatially-periodic burst with $A=1$ and $\eps = 0.01$.
    The canard solutions for $B=1.3$ and $B=1.5$ are shown in Fig.~\ref{fig:bursting}.}
    \label{tab:spikeno}
\end{table}

\section{Nearly Self-similar Dynamics of Spatially-Periodic Canards}
\label{s:selfsimilarity}

In this section, we show that the Hamiltonian $H_2$ on the hemisphere $\{r_2=0\}$ in $K_2$ has a self-similarity, and hence that its zero level set has an infinite scale invariance.
In turn, this self-similarity causes some classes of canards to have a nearly-infinite self-similarity for $0<\eps\ll 1$.

\subsection{Self-similarity of the Hamiltonian \texorpdfstring{$H_2$}{Lg}} \label{ss:selfsimilarHamiltonian}

We recall from \eqref{K2 Dynamics} that the dynamics in the central/rescaling chart $K_2$ are governed by
\begin{equation*}
\begin{split}
    U_2' &= P_2 \\
    P_2' &= A\B_2 - A^2 V_2 - \frac{1}{A} U_2^2 - r_2^2 F_2  \\
    V_2' &= Q_2 \\
    Q_2' &= U_2 -r_2^2 \left(A\B_2 -A^2 V_2 - \frac{1}{A} U_2^2 - r_2^2 F_2 \right).
\end{split}
\end{equation*}
Here, $F_2 = -\B_2 U_2 +  2A U_2 V_2 + r_2^2 U_2^2 V_2.$
On the invariant set $\{ r_2=0 \}$,
the unperturbed problem is
\begin{equation*}
\begin{split}
    U_2' &= P_2 \nonumber \\
    P_2' &= A\B_2 - A^2 V_2 - \frac{1}{A} U_2^2 \nonumber \\
    V_2' &= Q_2 \nonumber \\
    Q_2' &= U_2.
\end{split}
\end{equation*}
It is useful to scale the variables via 
\begin{equation*}
 U_2 = A^2 \hat{U}_2, \quad
 P_2 = A^2 \hat{P}_2, \quad
 V_2 = \hat{V}_2, \quad 
 Q_2 = A^2 \hat{Q}_2.
\end{equation*}
Hence, the unperturbed system is
\begin{equation}\label{K2-r=0}
\begin{split}
    {\hat{U}}_2' &= {\hat P}_2 \\
    {{\hat P}}_2' &= -{\hat V}_2 - A {\hat U}_2^2 + \frac{\B_2}{A} \\
    {\hat{V}}_2' &= A^2 {\hat Q}_2 \\
    {{\hat Q}}_2' &= {\hat U}_2.
\end{split}
\end{equation}

Let 
\begin{equation} \label{H2}
H_2(\hat{U}_2, \hat{P}_2, \hat{V}_2, \hat{Q}_2)
= \frac{1}{2} \hat{P}_2^2 - \frac{1}{2} A^2 \hat{Q}_2^2 -  \hat{U}_2 \hat{V}_2 
   - \frac{A}{3} \hat{U}_2^3 + \frac{\B_2}{A} \hat{U}_2.
\end{equation}
Then, we see that the unperturbed system is Hamiltonian,
\begin{equation*}
{\hat{\bf{U}}}_2' 
= J \nabla_{(\hat{Q}_2, \hat{V}_2, \hat{P}_2, \hat{U}_2)} H_2, 
\end{equation*} 
where 
${\hat{\bf U}}_2 = [ \hat{U}_2, \hat{P}_2, \hat{V}_2, \hat{Q}_2 ]^T$ and
$J = \begin{bmatrix} 
0 & {\bf I} \\ -{\bf I} & 0 
\end{bmatrix}$.
The Hamiltonian has the following self-similarity:
\begin{equation}\label{H2selfsimilar}
H_2(\lambda^2 U_2, \lambda^3 P_2, \lambda^4 V_2, \lambda^3 Q_2; \lambda^4 \B_2)
= \lambda^6
H_2(U_2, P_2, V_2, Q_2; \B_2).
\end{equation}
Hence, the zero level set is invariant under this scaling.
(We add that the powers of $\lambda$ are the same as the powers of $r$ used in the geometric desingularization \eqref{powersofr} in Sec.~\ref{s:desingularization}.)

\subsection{Nearly self-similar dynamics of classes of spatially-periodic canard solutions}
\label{ss:nearlyselfsimilar}

\begin{figure}[h!t]
  \centering
  \includegraphics[width=5in]{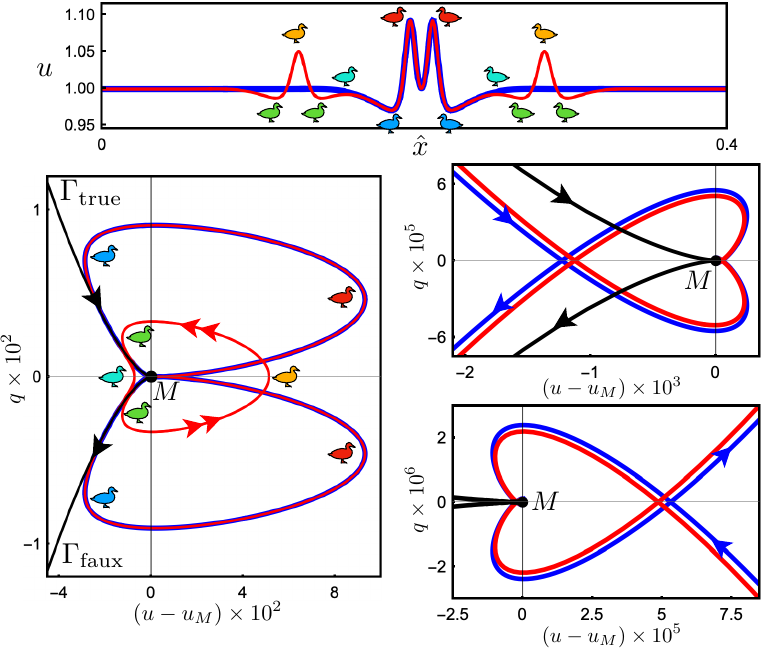}
  \put(-346,298){(a)}
  \put(-360,216){(b)}
  \put(-174,222){(c)}
  \put(-174,106){(d)}
  \caption{Two spatially-periodic solutions (red and blue), both with period $T=1500$ (and $k \approx 0.419$), exhibiting nearly self-similar dynamics for $A=1, B=1$, and $\eps = 0.01$.
  (a) $u(\hat{x})$. 
  (b) The projections into the $(u,q)$ plane.
  $M$ is a RFSN-II point, where the RFS and $E$ coincide. 
  (c) and (d) Magnifications of successively smaller neighborhoods of $M$ reveal several levels of the nearly self-similar cusp-like structure centered around $M$.
  (For $\hat{x} \in [0.4,1]$ (not shown), $u(\hat{x})$ is near the homogeneous state $u=1$.)
  }
  \label{fig:selfsimilarB1}
\end{figure}
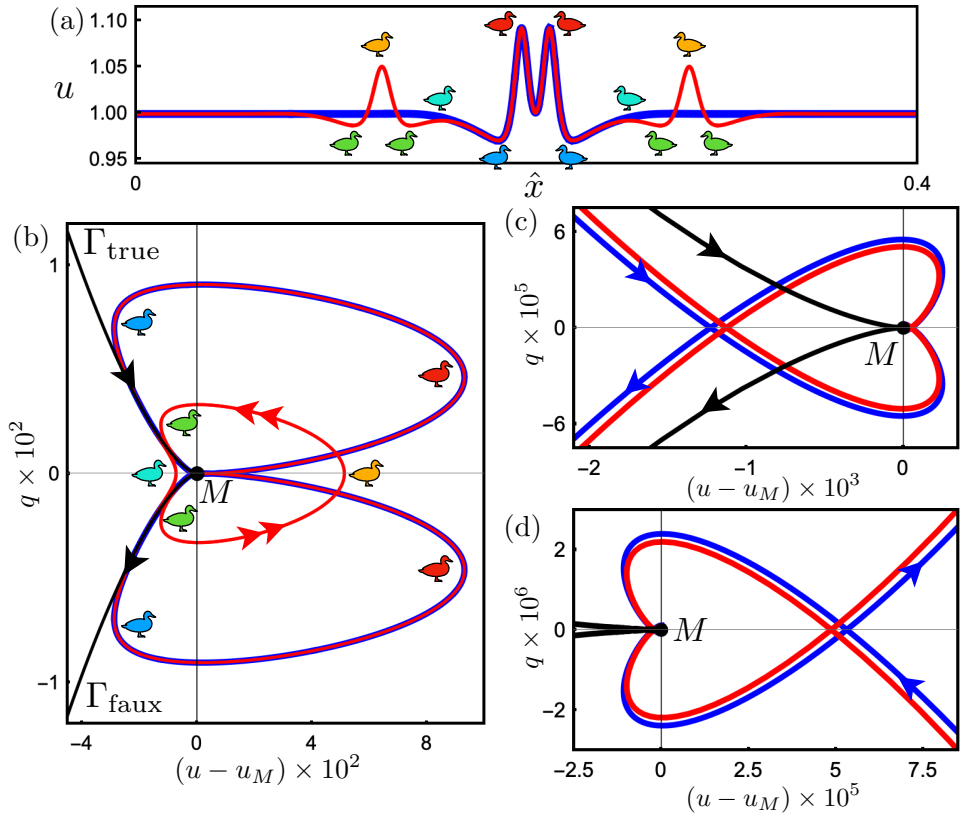

\begin{figure}[h!t]
  \centering
  \includegraphics[width=5in]{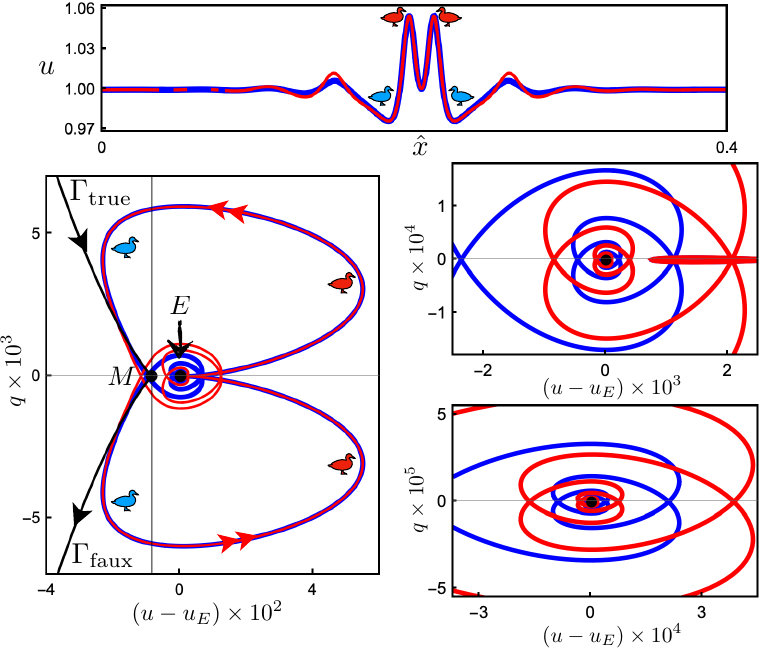}
  \put(-346,298){(a)}
  \put(-360,216){(b)}
  \put(-174,222){(c)}
  \put(-174,106){(d)}
  \caption{Two spatially-periodic solutions (red and blue), both with period $T=1500$, exhibiting nearly self-similar dynamics for $A=1, B=1.0175$, and $\eps=0.01$. 
  Since $B>1$, $M$ is a RFS and $E$ is a center (recall Fig.~\ref{fig:desingularized}(c)).
  (a) $u(\hat{x})$.
  Both patterns have pulses on the wings of the two central pulses, with spatial oscillations there.
  (b)-(d) Projections into the $(u,q)$ plane. Magnifications of successively smaller neighborhoods of $E$ reveal several levels of the nearly self-similar winding around $E$.
  These solutions belong to a different branch than those in Fig.~\ref{fig:maxcanards}.
  In (c), the nearly horizontal segment of the red orbit essentially on the horizontal axis (which enters and exits on the right) is a piece of the two fast pulses near the local minimum.
  }
  \label{fig:selfsimilarB1.0175}
\end{figure}
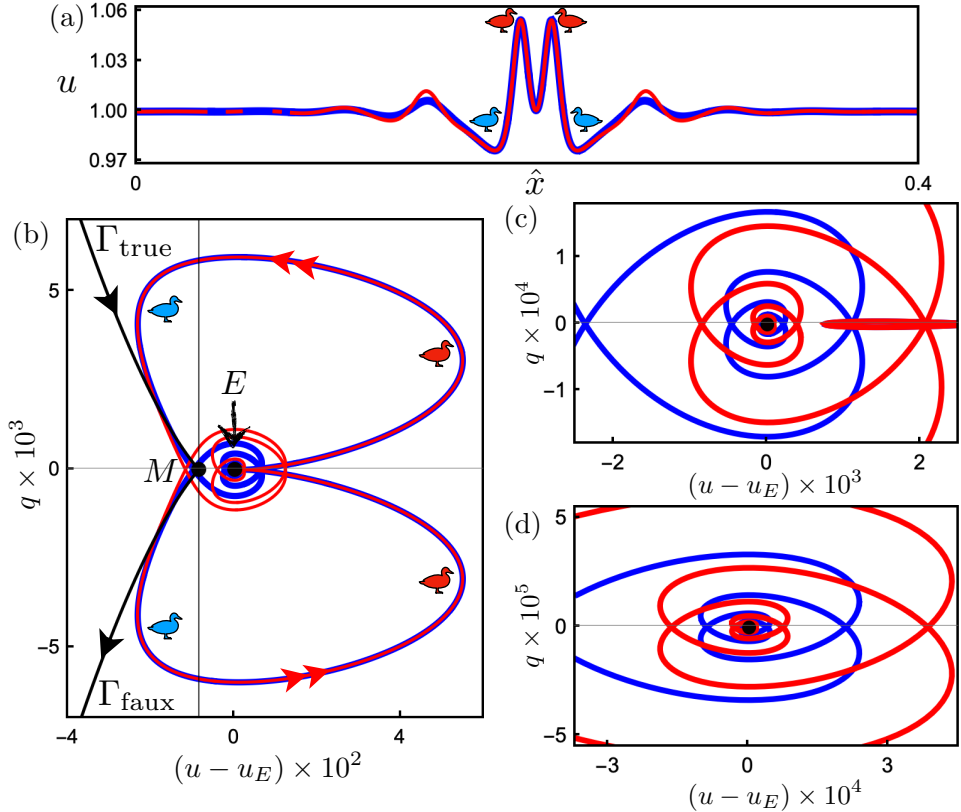

For $0<r_2\ll 1$, the scale invariance manifests (numerically) as nearly self-similar dynamics. 
Furthermore, for $B=1$ or sufficiently close to one, the organizing center for the self-similar dynamics is $M_{\rm RFSN-II}$ (Fig.~\ref{fig:selfsimilarB1}); and, for $B> 1 +\mathcal{O}(\eps)$ the organizing center is the equilibrium $E$  (Fig.~\ref{fig:selfsimilarB1.0175}).

In Fig.~\ref{fig:selfsimilarB1}(a), we show two spatially-periodic solutions (red and blue) of \eqref{Brusselator} with nearly self-similar dynamics close to $M_{\rm RFSN-II}$. 
The red solution has 14 segments, of which we describe the first 7. 
Symmetry gives the remaining 7 segments. 
As shown in (b), the segments are:
\begin{enumerate}[(i)]
\setlength{\itemsep}{0pt}
    \item A slow segment along the faux canard $\Gamma_{\rm faux}$ from the RFSN-II point $M$ to the (takeoff) point near the local $u$-minimum (first, lower green ducky).   
    \item A fast $\mathcal{O}(\eps)$-amplitude pulse from the takeoff point to the local $u$-maximum (yellow ducky). 
    \item The symmetric image of the previous segment (to leading order) takes the solution from the yellow ducky back to the (touchdown) point (second, upper green ducky). 
    \item A slow segment then carries the solution to the cyan ducky at $q=0$. 
    \item A long, slow segment near $\Gamma_{\rm faux}$ from the cyan ducky to the (takeoff) point (lower blue ducky).
    \item A fast $\mathcal{O}(\eps)$-amplitude pulse from the blue ducky to the $u$-maximum (red ducky). 
    \item A fast segment from the red ducky to the neighborhood of the RFSN-II point $M$. 
\end{enumerate}
The nearly self-similar dynamics shown in Figs.~\ref{fig:selfsimilarB1}(c) and (d) are confined near the RFSN-II point in phase space, and are observed along the (seemingly) flat plateaus in the spatial profile. 

The blue solution in Fig.~\ref{fig:selfsimilarB1} has a similar (but simpler) symmetric deconstruction:  
\begin{enumerate}[(i)]
\setlength{\itemsep}{0pt}
    \item A slow segment along $\Gamma_{\rm faux}$ from $M$ to the takeoff point (lower blue ducky). 
    \item A fast $\mathcal{O}(\eps)$-amplitude pulse from the blue ducky to the $u$-maximum (red ducky). 
    \item Another fast segment from the red ducky to the neighborhood of the RFSN-II point $M$. 
\end{enumerate}
As with the red solution, the nearly self-similar dynamics of the blue solution are restricted to a small neighborhood of $M$ and are observed along the flat plateaus in the spatial profile. 

Fig.~\ref{fig:selfsimilarB1.0175} shows two spatially-periodic canards of a different type, with self-similar type oscillations about $E$ (instead of $M=M_{\rm RFS}$). 
These are typical of the patterns observed for $B-1$ sufficiently small and positive.
They both have long slow segments near $\Gamma_{\rm faux}$ from near $M$ down to the takeoff point (lower blue ducky).
Also, both exhibit fast pulses on which $u$ reaches maximal values (red ducky) and then local minima near $M$.
Then, the patterns exhibit the reflected images to the right.
Lastly, after the long slow segment near the true canard from the upper blue ducky into a neighborhood of $M$, both orbits exhibit several levels of self-similar oscillations about $E$ (see Fig.~\ref{fig:selfsimilarB1.0175}(c) and (d)).
Some are visible in Fig.~\ref{fig:selfsimilarB1.0175}(a), occurring before $\hat{x}$ reaches 0.4.
Then, the segments of the solutions along which $u$ oscillates ever closer to $E$, which correspond to deeper levels of the nearly self-similar dynamics, lie to the right (in the portion of the period not shown).

\section{Direct Numerical Simulations of \texorpdfstring{\eqref{Brusselator}}{Lg} on Periodic Domains}
\label{s:PDE-numerics}

In this section, we report on the results of direct numerical simulations of \eqref{Brusselator} to study the stability of spatially-periodic solutions of~\eqref{Brusselator} on intervals with periodic boundary conditions.
Then, in the next section,
we present simulations on large intervals with Neumann boundary conditions.

Direct simulations were performed with Mathematica's NDSolveValue-package with the ``MethodOfLines" and ``SpatialDiscretization" prescribed by 
\{``TensorProductGrid",``MinPoints"\(\to\)10000\} on
intervals \([\tfrac{-L}{2},\tfrac{L}{2}]\), subject to periodic boundary conditions.
The wavenumber is \(k=\tfrac{2\pi}{L}\).

\begin{figure}[h!t]
    \centering
    \includegraphics[width=.45\linewidth]{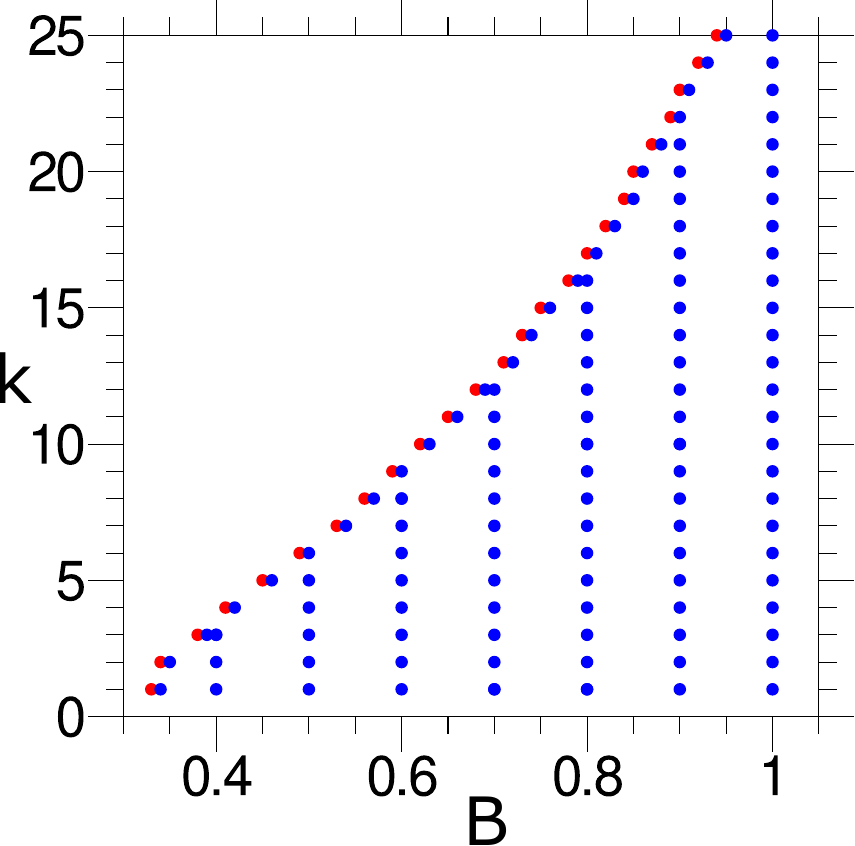}
    \caption{Numerical approximation of a large portion ($B < B_T$) of the regime of stable, spatially-periodic solutions of the Brusselator PDE~\eqref{Brusselator} obtained from direction simulations on $[-\tfrac{L}{2},\tfrac{L}{2}]$ with periodic boundary conditions. 
    At blue markers, the data converges to numerically stable, spatially-periodic solutions.
    At red markers, solutions converge to the homogeneous state. 
    Here, \(A=1\) and \(\eps=0.01\), so that \(k_{T}=10\), \(B_{T}=1.0201\), and the homogeneous state is linearly stable. 
    The Fourier mode, \(k\), is controlled by imposing periodic boundary conditions on \(\tfrac{2\pi}{k}\) and verifying that the steady state contains only a single period.
    }
    \label{fig: 2D stability regime}
\end{figure}

The initial data were set equal to the homogeneous equilibrium value, and a small perturbation term was added to $u$, consisting
of a large positive value on a narrow interval near the origin and a smaller negative value outside of the interval.
The exact widths and amplitudes were varied to find possible basins of attraction, and typical values are \(+10\) on an interval of length \(0.01*L\) and \(-0.5\) outside.
Generally, for these initial data  to be in the basins of attraction, larger amplitude perturbations are required for larger $B$ values and for solutions with smaller wavenumbers.
A large portion of the stability regime is presented in Fig.~\ref{fig: 2D stability regime}, for \(A=1\) and \(\eps=0.01\), so that \(B_{T}=1.0201\).
Simulations were performed for \(k\): $1 \to 25$ using steps of size 1 and for \(B\): $0.34 \to 1$ using steps of size 0.01 near the boundary of the stability regime and of size 0.1 in the interior.
In the bottom left corner (small \(B\) and \(k\)), the boundary appears to be concave down.
Near \((B,k)=(0.5,6)\), there appears to be an inflection point.
Lastly, there appears to be a second inflection point near \((B,k)=(0.88,21)\).
Several features are similar to those in Fig.~\ref{fig:existenceballoon}.
\begin{figure}[h!t]
    \centering
    \includegraphics[width=0.85\linewidth]{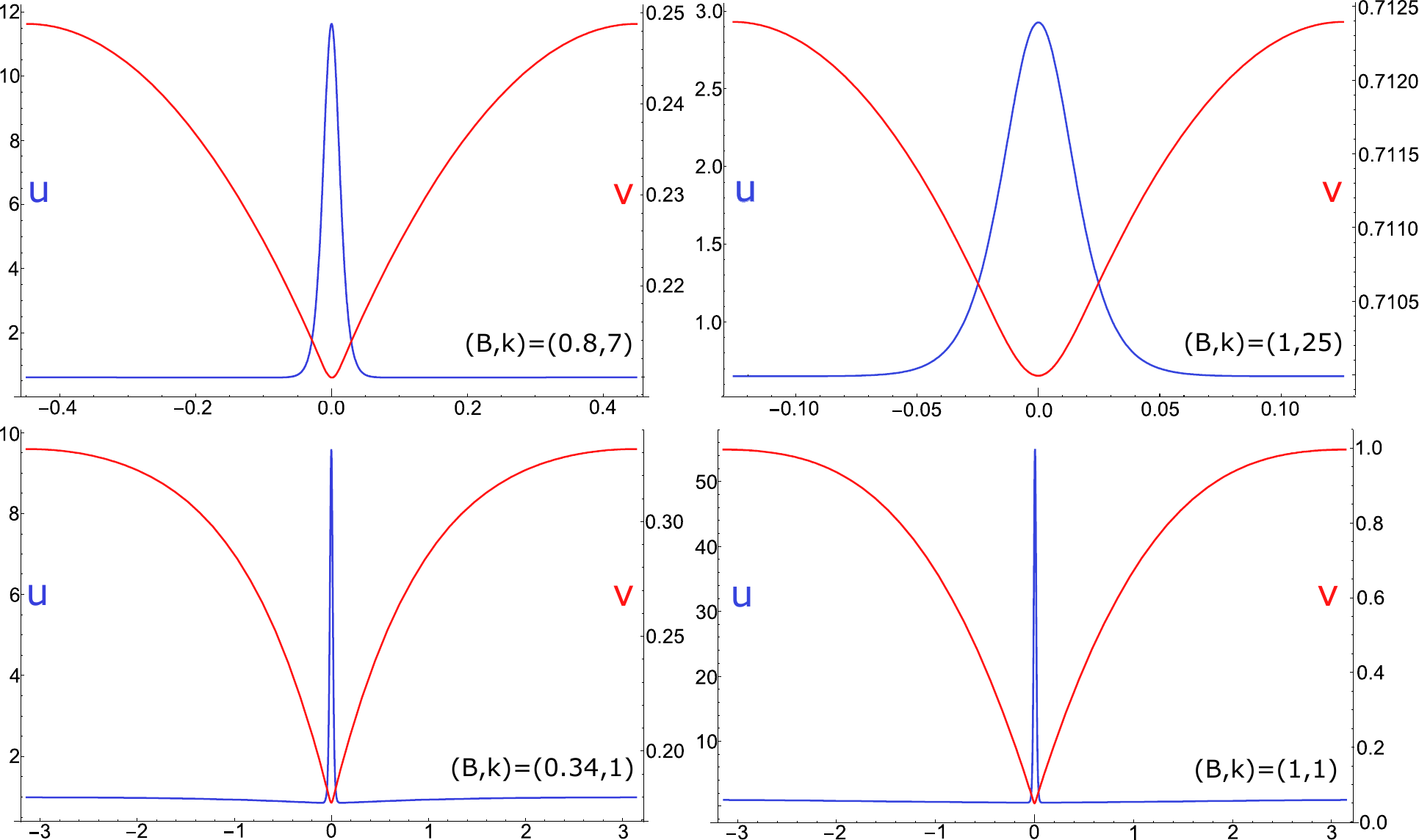}
    \put(-420,235){(a)}
    \put(-210,236){(b)}
    \put(-420,112){(c)}
    \put(-213,112){(d)}
    \caption{Stable spatially-periodic states observed in simulations of \eqref{Brusselator}.
    \(u\) and \(v\) vary over different scales on the left and right axes, respectively.
    The parameters represent the (a) middle, (b) top right, (c) bottom left, and (d) bottom right regions of the stability regime in Fig.~\ref{fig: 2D stability regime}.
    Here, \(A=1\) and \(\eps=0.01\).
    }
    \label{fig: Steady States}
\end{figure}

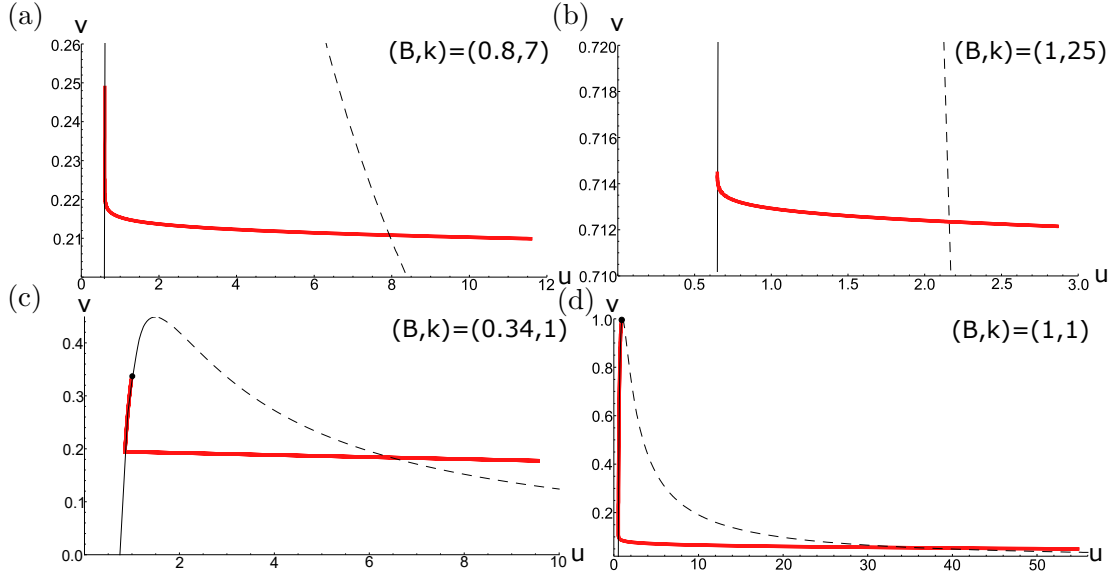
\begin{figure}[h!t]
    \centering
    \includegraphics[width=0.85\linewidth]{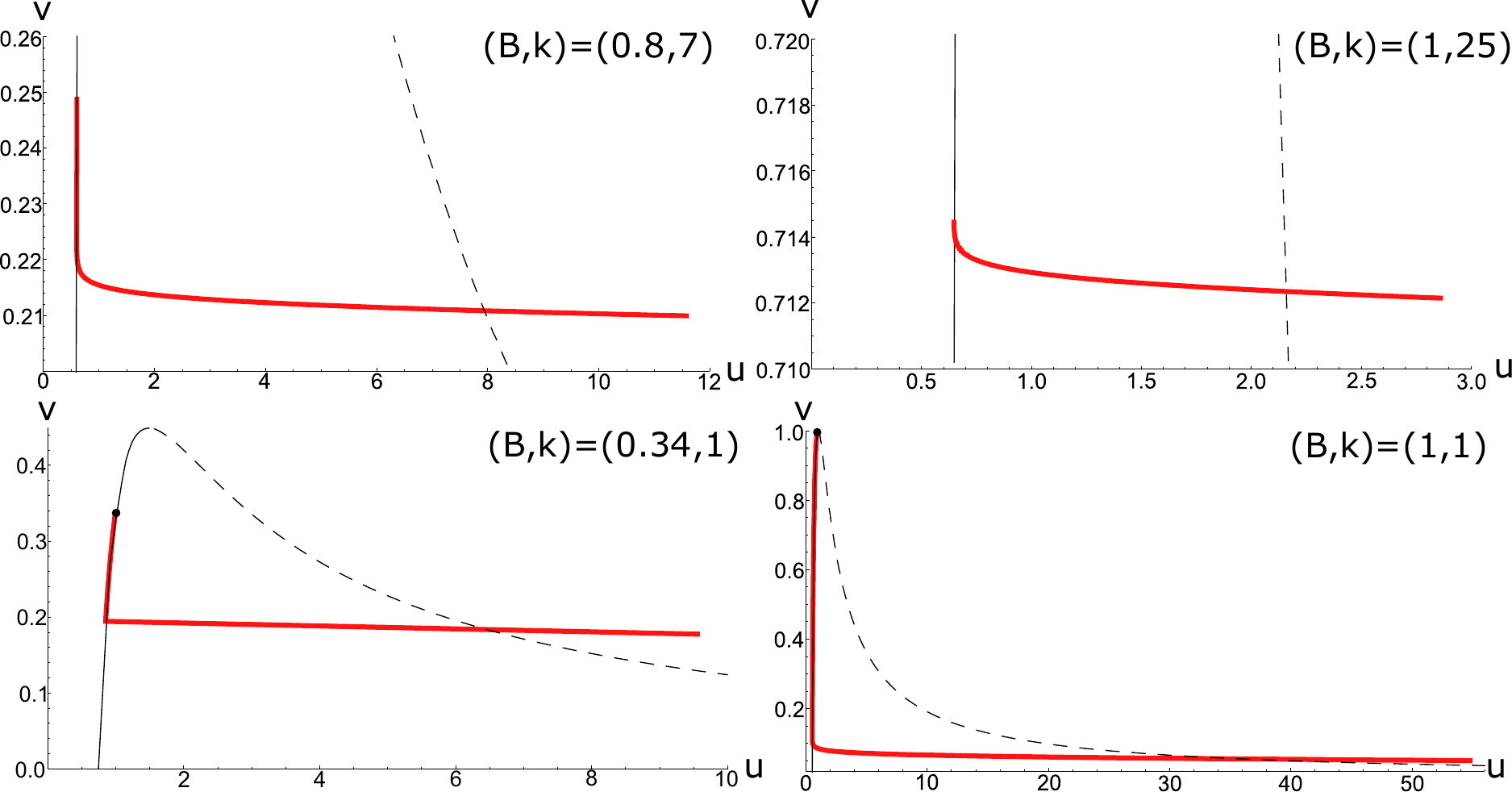}
    \put(-416,207){(a)}
    \put(-210,207){(b)}
    \put(-416,100){(c)}
    \put(-208,99){(d)}
    \caption{Projections into the $(u,v)$-plane of the stable spatially-periodic solutions in Fig.~\ref{fig: Steady States}.
    In each period, the trajectories (red curves) have a nearly-horizontal segment on which $u$ has its large-amplitude pulse and $v$ is close to its minimum, followed by a long slow segment near the left branch of the nullcline (black curve) along which $u$ slowly varies and $v$ gradually varies to its maximum and back.
    The $u$-pulse is an oscillation about the right branch of the nullcline (dashed black).
    }
    \label{fig: Steady State Projections}
\end{figure}

Fig.~\ref{fig: Steady States} presents the $u$ and $v$ profiles of four stable, spatially-periodic solutions observed at parameter values selected from various parts of the stability regime.
$u$ has a large-amplitude pulse on a narrow interval, and it varies gradually outside this interval.
$v$ varies gradually between minimum and maximum values, which are $\mathcal{O}(\eps)$ close. 
$v$ reaches a minimum at the center of the pulse interval.
(Note the differences in the scales on the vertical axes for $u$ and $v$.)

For all solutions here, the absolute pulse widths are approximately the same ($\mathcal{O}(\eps)$).
However, relative to the period lengths, the pulse widths decrease with $k$.
Furthermore, the amplitudes of the pulses in $u$ generally increase as $k$ decreases (for each fixed $B$).
Similarly, for fixed $k$, the amplitudes increase as the value of $B$ is taken to be further from the stability boundary.

Fig.~\ref{fig: Steady State Projections} shows the projections into the $(u,v)$-plane of the stable, spatially-periodic states of Fig.~\ref{fig: Steady States}.
The nearly-horizontal segments ($v$ near its minimum), correspond to the narrow, large-amplitude $u$ pulses that make one oscillation about the right branch of the nullcline.
On the complementary intervals, the trajectories are close to the left branch of the $u$-nullcline, where $v$ varies gradually.
The smaller $k$, the closer the trajectory gets to the local maximum of the nullcline.
These solutions are quantitatively similar to those in Fig.~\ref{fig:pulseconstruction}.

\begin{figure}[h!tbp]
    \centering
    \includegraphics[width=0.9\linewidth]{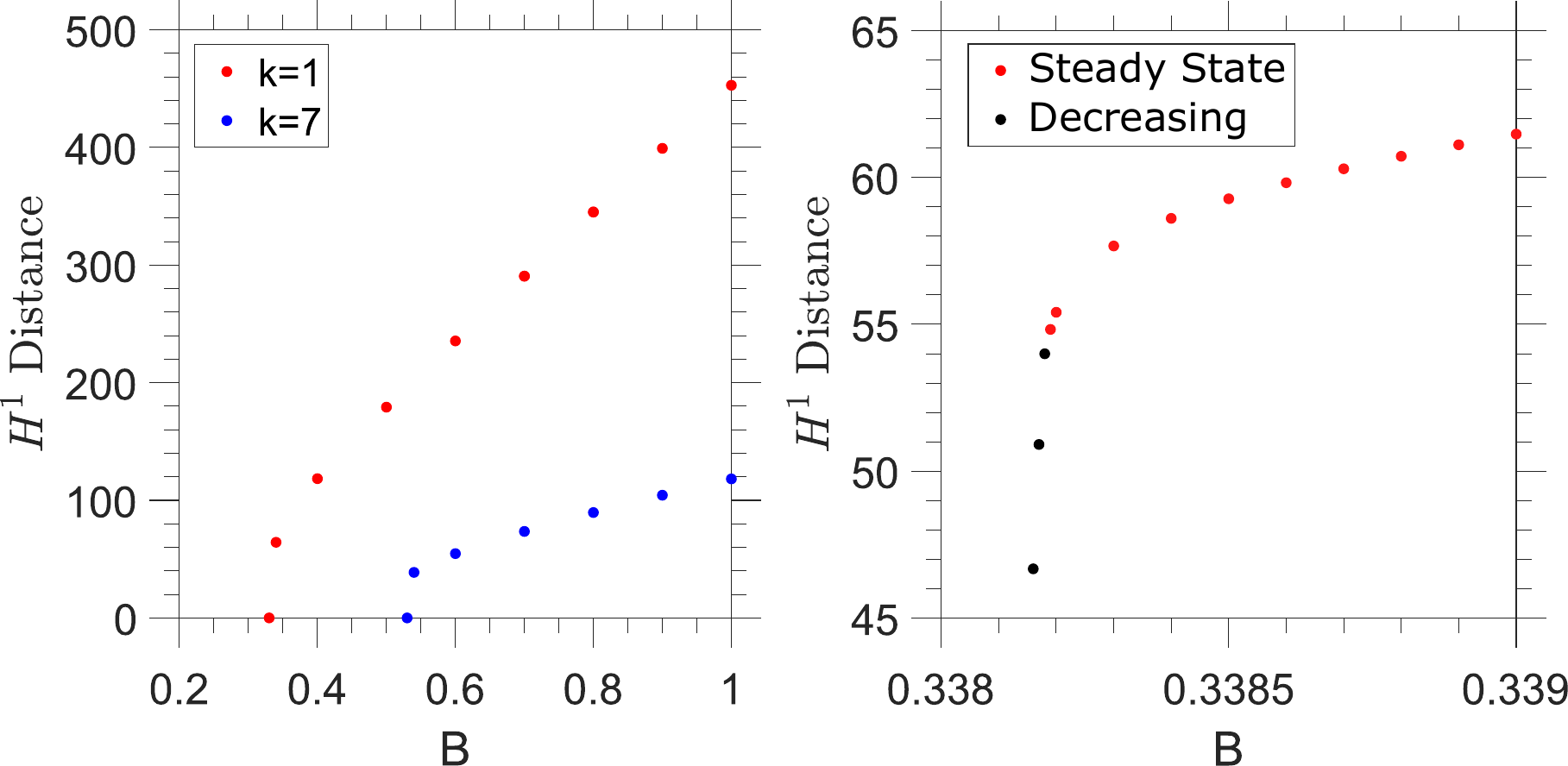}
    \put(-430,200){(a)}
    \put(-210,200){(b)}
    \caption{(a) The distance in the \(H^{1}\) norm between stable spatially-periodic states and the homogeneous state as a function of $B$ for the solutions in Fig.~\ref{fig: 2D stability regime} with \(k=1\) (red) and \(k=7\) (blue).
    The leftmost red and blue markers occur at \(B = 0.34 \) and \(B = 0.53 \), respectively, where the $H^1$ norm is zero to machine precision.
    (b) The $H^1$ distance for solutions with $k=1$ (simulated until time \(t_f=150\)) on a finer scale. 
        }
    \label{fig: H1VsB} 
\end{figure}
The existence of branches of stable, spatially-periodic states suggests that for each \(k\) there is a fold at some \(B_{F}(k)<B_{T}\).
The folds connect branches of unstable periodic solutions which bifurcate off the homogeneous state in the subcritical Turing bifurcation to the branches of (numerically) stable solutions.
To investigate, we plot the distance in the \(H^{1}\) norm (on the interval 
\([\tfrac{-L}{2},\tfrac{L}{2}]\)) between the steady states and the homogeneous equilibrium.

Fig.~\ref{fig: H1VsB}(a) shows the $H^1$ distance as a function of $B$ for constant $k$ ($k=1$ and $k=7$) in the empirical stability regime.  
The (saddle-node) bifurcations $B_F(1) \approx 0.34$ and $B_F(7) \approx 0.53$ correspond to the saddle-node pairs $(k, B_F(k))$ in the spatial ODE.
(Recall Fig.~\ref{fig:bifurcation3d}, where 
the boundaries of the blue surface show that the bifurcation is at $k \approx 1.0$ for $B  = 0.34$ and at $k\approx 7.2$ for $B = 0.53$.)
Near the boundary, the $H^1$ distance grows steeply as $B$ is taken further from critical.
Then, away from the boundary, the $H^1$ distance increases approximately linearly with $B$.
\begin{figure}[h!tbp]
    \centering
    \includegraphics[width=0.75\linewidth]{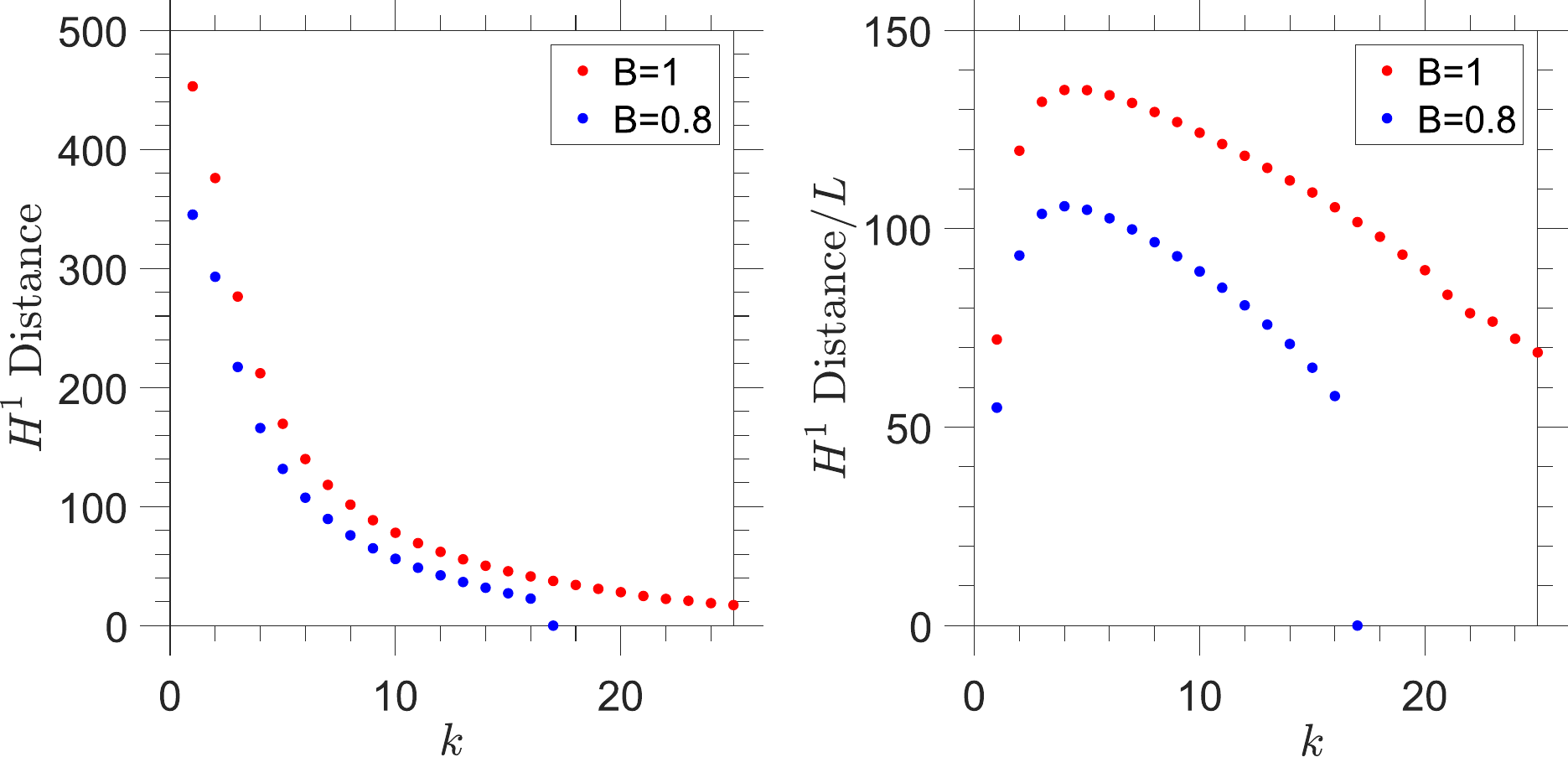}
    \put(-355,165){(a)}
    \put(-172,165){(b)}
    \caption{(a) The distance in the \(H^{1}\) norm between stable periodic steady states and the homogeneous state as a function of $k$. 
    (b) The \(H^{1}\) distance normalized by domain length.
    For each fixed \(B\), the \(H^{1}\) distance decreases to approximately \(30\) as \(k\) approaches the boundary value.}
    \label{fig: H1VsK}
\end{figure}

Fig.~\ref{fig: H1VsB}(b) shows the \(H^{1}\) distance for a finer sampling of \(B\) values near the empirical value $B_F(1) \approx 0.338185$.
For \(B\geq 0.33819\) (red), solutions have reached steady state, and the $H^1$ distance scales in proportion to \(\sqrt{B-B_{F}(1)}\).
In contrast, for \(B\leq 0.33818\) (black), solutions stay near slowly-varying spatially-periodic states for long times before decreasing rapidly to the homogeneous state.
For \(B\) further from \(B_{F}(1)\), the time spent near the periodic solutions decreases, as is common in simulations near saddle-node bifurcations for fixed integration times.

Fig.~\ref{fig: H1VsK}(a) shows the \(H^{1}\) distance between the stable, spatially-periodic solutions and the homogeneous state as a function of $k$ for $B=0.8$ and $B=1$.
The distance is strictly decreasing and concave up in \(k\), except near the boundary of the stability regime. 
Panel (b) shows the distance normalized by the domain length.
It is concave down in \(k\), with an absolute maximum near \(k=4\).

These plots and the trends in Figs.~\ref{fig: Steady States} and \ref{fig: Steady State Projections} reveal the relative strengths of two competing effects.
As \(k\) decreases, the pulse amplitude and sharpness increase, while the relative width of the pulse interval decreases. 
The former causes an increase in the $H^1$ distance, while the latter causes a decrease because the solution is close to the homogeneous state over a longer portion of the interval. 
Hence, Fig.~\ref{fig: H1VsB} shows that the latter dominates the former as $k$ decreases (below $k\approx 4$). 
The narrowing of the pulse interval has a stronger effect than an amplitude increase.


\section{Numerical simulations on large domains}
\label{s:PDE-largedomains}

\begin{figure}[h!t]
        \includegraphics[width=0.4\textwidth]{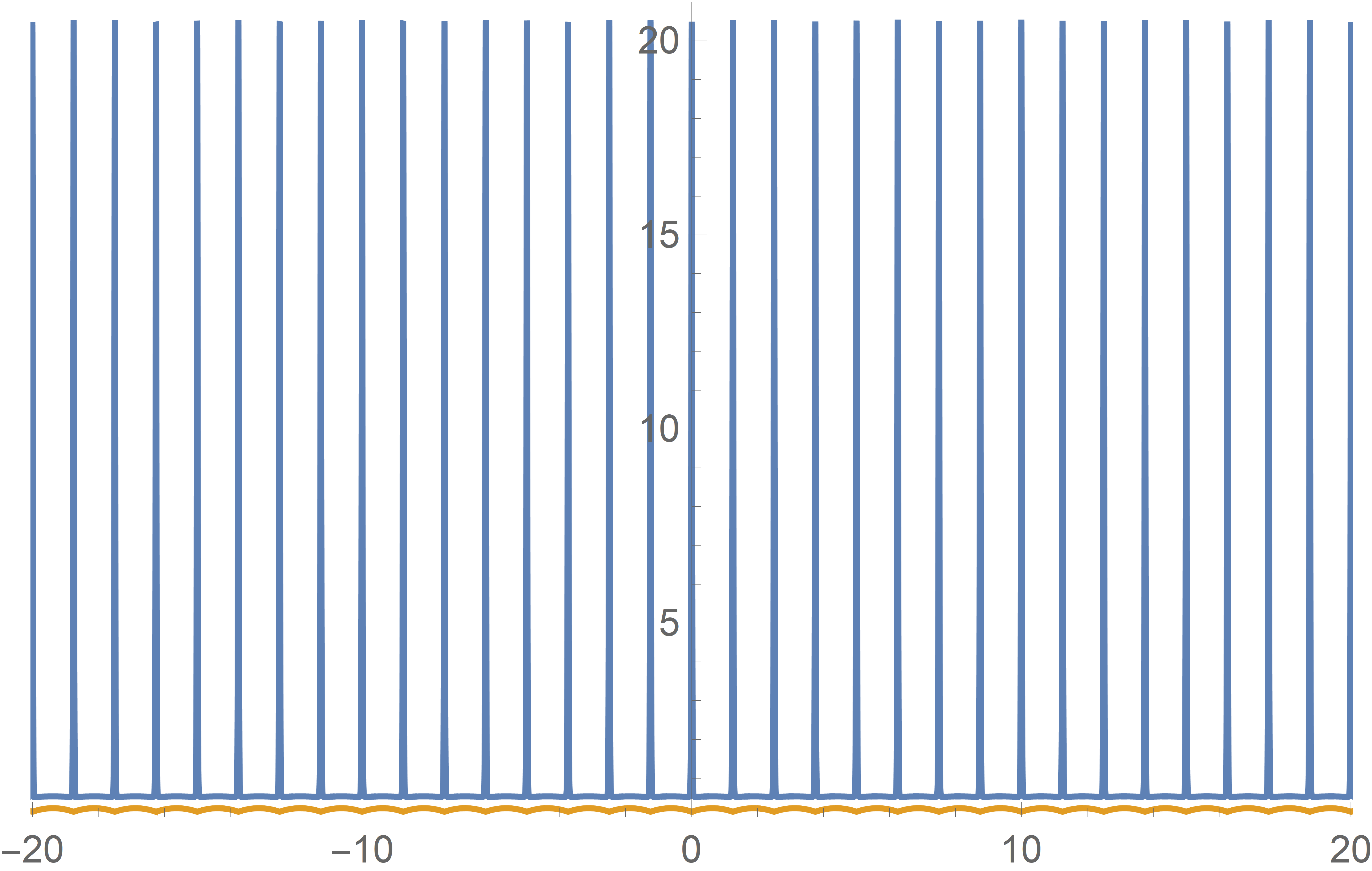}
        \qquad 
        \includegraphics[width=0.4\textwidth]{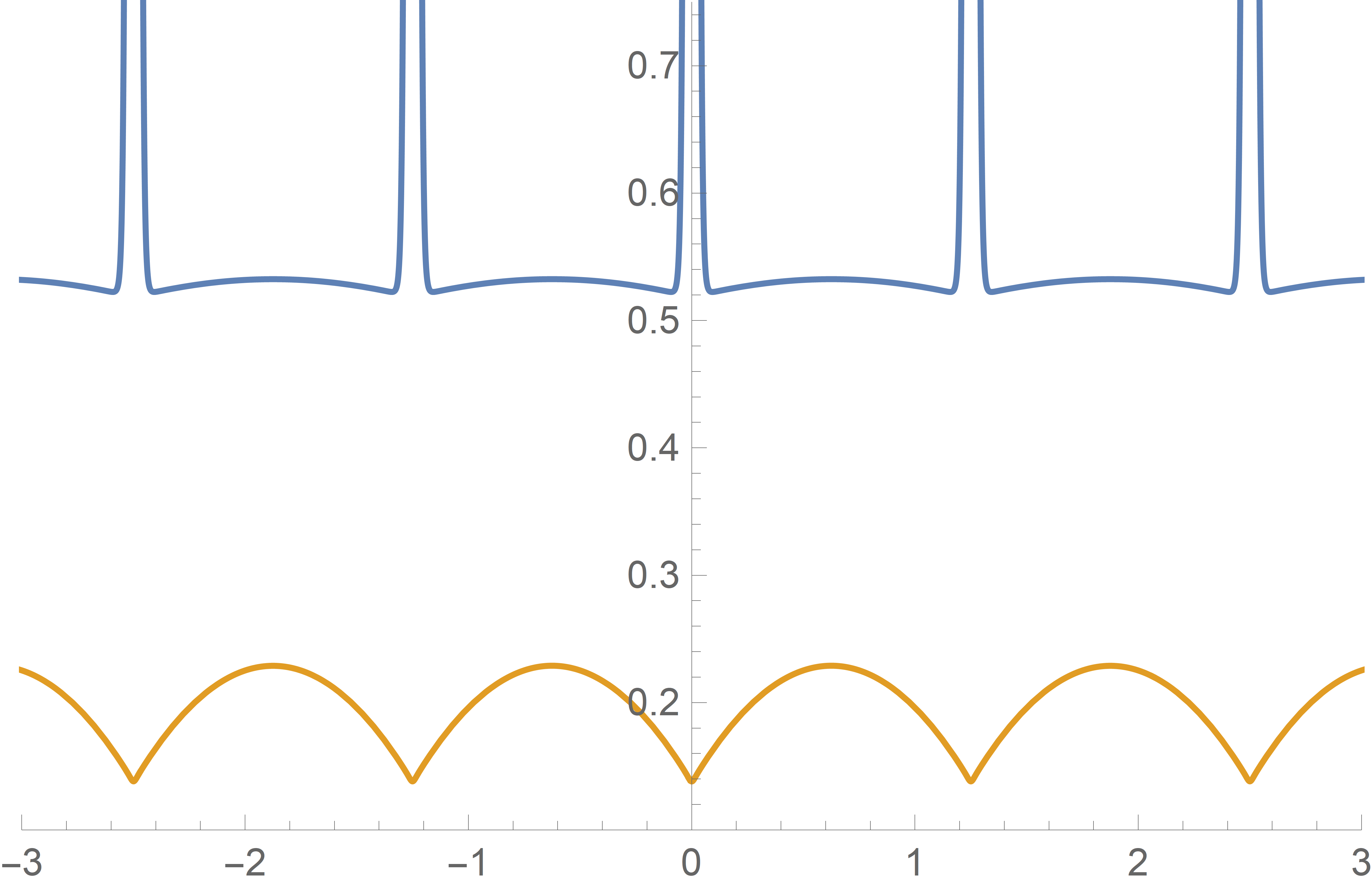}
\put(-420,110){(a)}
\put(-200,110){(b)} \\
\includegraphics[width=0.4\textwidth]{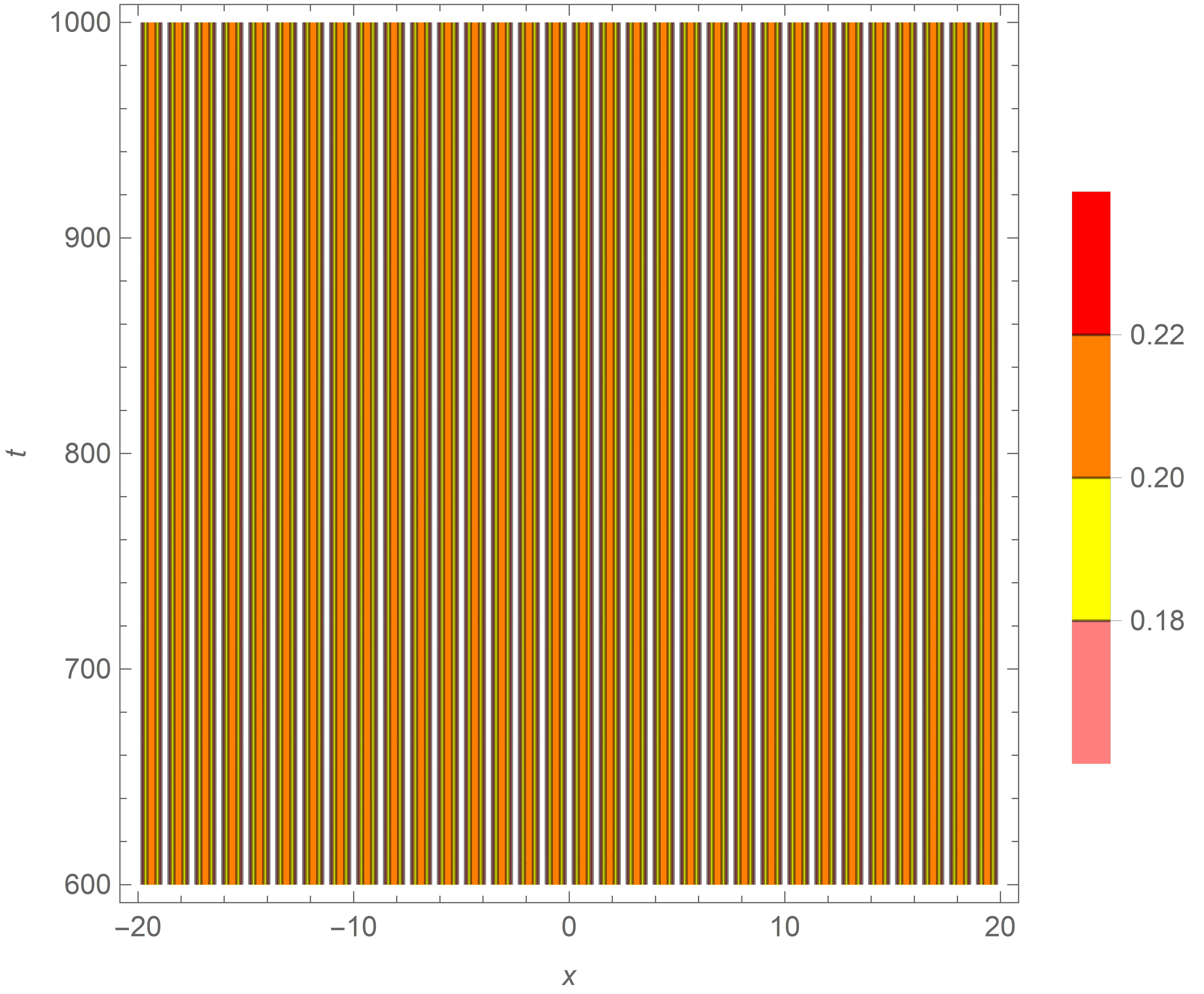}
\qquad 
\includegraphics[width=0.4\textwidth]{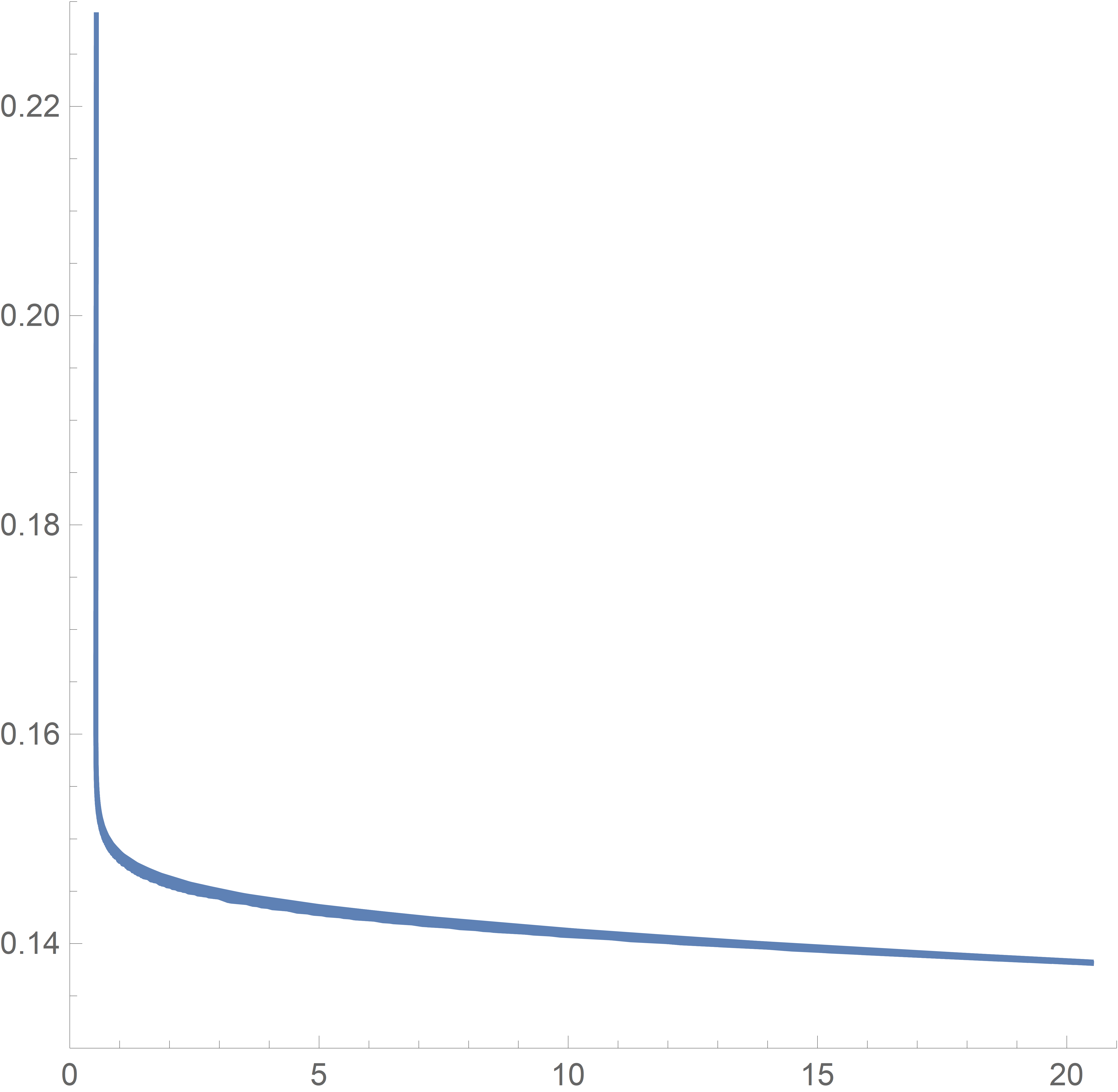}
\put(-420,145){(c)}
\put(-200,170){(d)} 
\caption{A stable spatially-periodic canard pattern observed in direct numerical simulations of \eqref{Brusselator} with wavelength (spatial period) $T = 5/4$ and wavenumber $k = 8 \pi/5$ on the domain $[-30,30]$ with homogeneous Neumann boundary conditions starting from the initial conditions $(u(x,0), v(x,0)) = (1 + 0.1 [(\cos 4 \pi x/5)]^{8}, 1 + 0.1 [\sin(4\pi x/5)]^2)$. Both the initial and the final patterns contain $48$ full periods; however, only the portion on $[-20,20]$ is shown. 
(a) $u(x,1000)$ (blue) and $v(x,1000)$ (orange).
(b) A zoom of (a).
(c) A contour plot of $v(x,t)$.
(d) $(u(x,1000),v(x,1000))$ for $x \in [20,20]$.
The initial conditions $u(x,0)$ and $v(x,0)$ are sinusoidal. 
The parameters are $A=1$, $B=1$, and $\eps=0.01$.}
\label{fig:sim1}
\end{figure}

We performed simulations on large domains with homogeneous Neumann boundary conditions, as approximations for the PDE \eqref{Brusselator} on the real line.
First, we present a stable spatially-periodic pattern that evolves from small-amplitude periodic initial data with the same wavenumber inside the Busse balloon (see Fig.~\ref{fig:sim1}).

Consistent with weakly nonlinear analysis, the initial data are $\mathcal{O}(\sqrt{| B_{\rm T} - B |})$ close to the state $(A, B/A)$ (recall Sec.~\ref{s:TH}), and the spatial variation about $(A, B/A)$ has a sinusoidal structure of 48 full periods on $[-30,30]$.
The initial data rapidly (before the time $t=30$) becomes an $\mathcal{O}(1)$ spatially-periodic pattern with fast-slow structure, consisting of localized pulses in $u$ and small-amplitude, scalloped, spatial oscillations in $v$. 
$T=1.25$, with wavenumber $k=8\pi/5$, and the number of pulses stays the same during the simulation.
The attractor has the form of the patterns constructed in Fig.~\ref{fig:pulseconstruction}(b), and is similar to those in Fig.~\ref{fig: Steady State Projections}.

\begin{figure}[h!t]
\centering
        \includegraphics[width=0.46\textwidth]{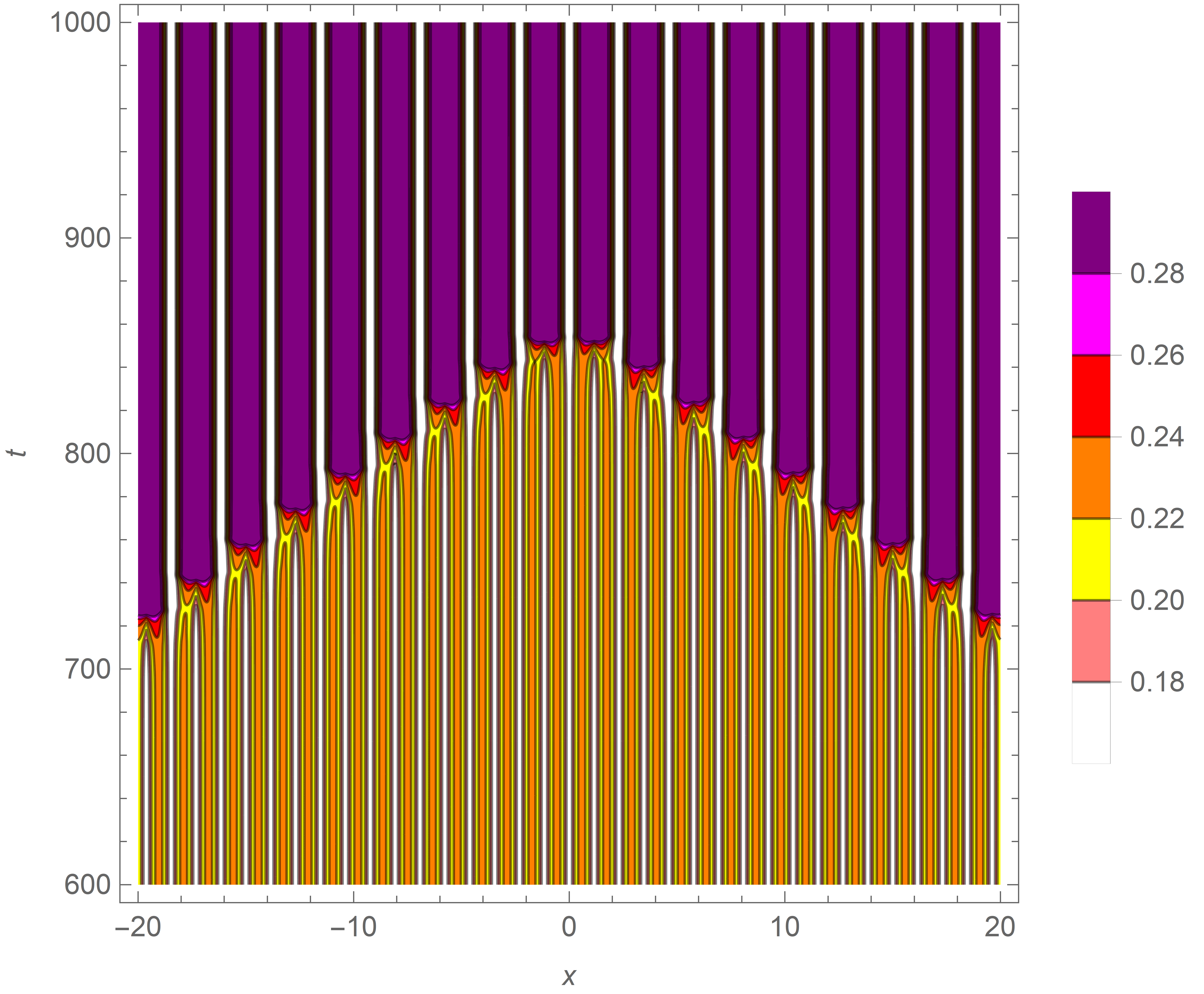}
        \includegraphics[width=0.46\textwidth]{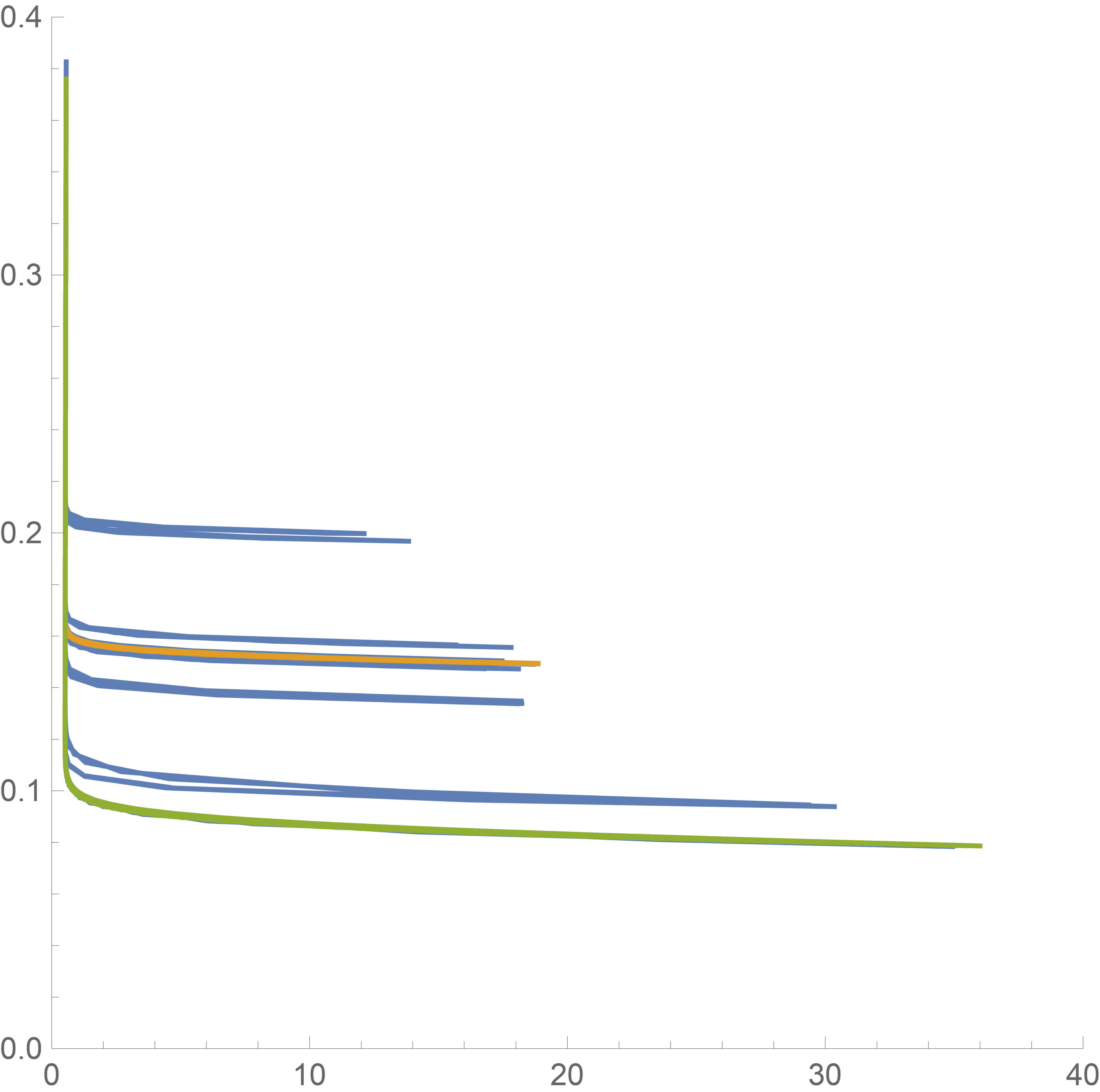}
        \put(-445,170){(a)}
\put(-230,170){(b)} 
\\ 
\vskip0.1truein 
        \includegraphics[width=0.31\textwidth]{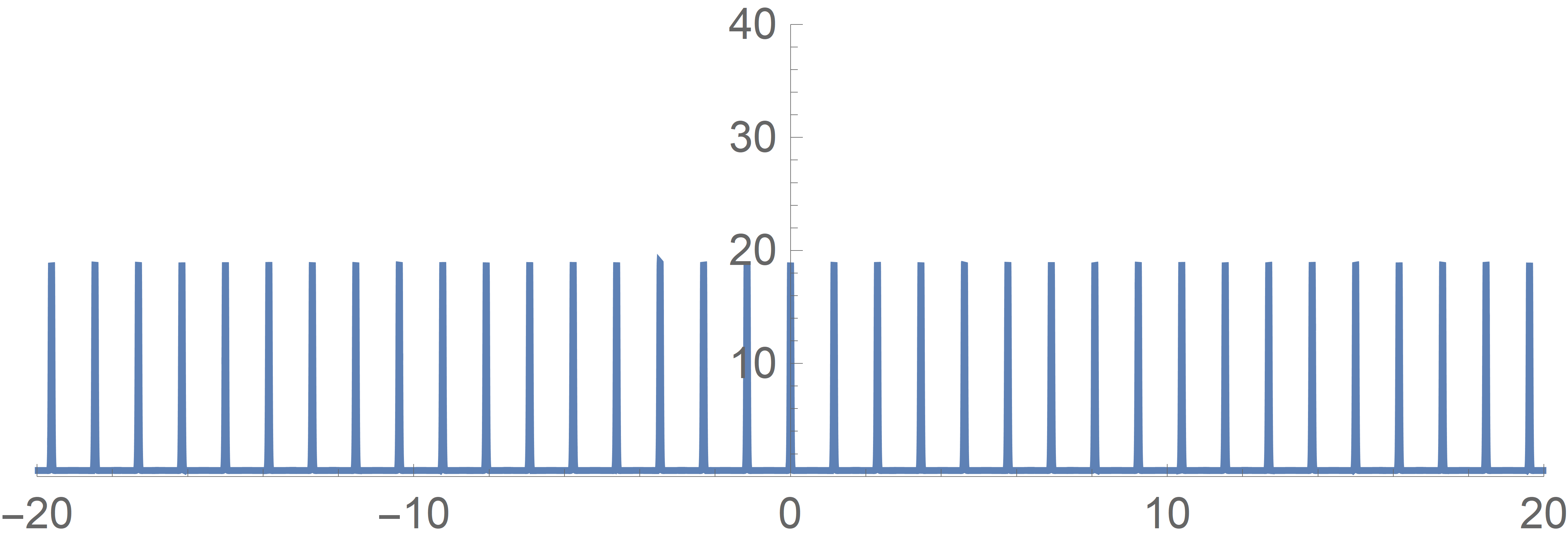}
        \includegraphics[width=0.31\textwidth]{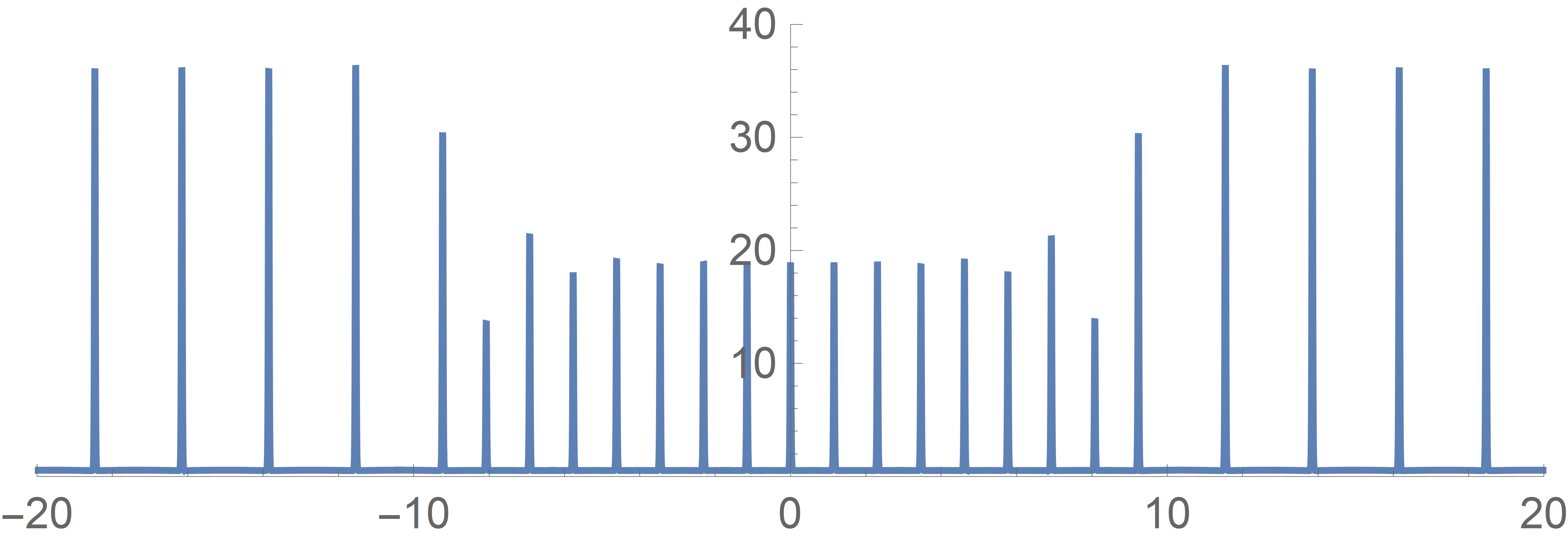}
        \includegraphics[width=0.31\textwidth]{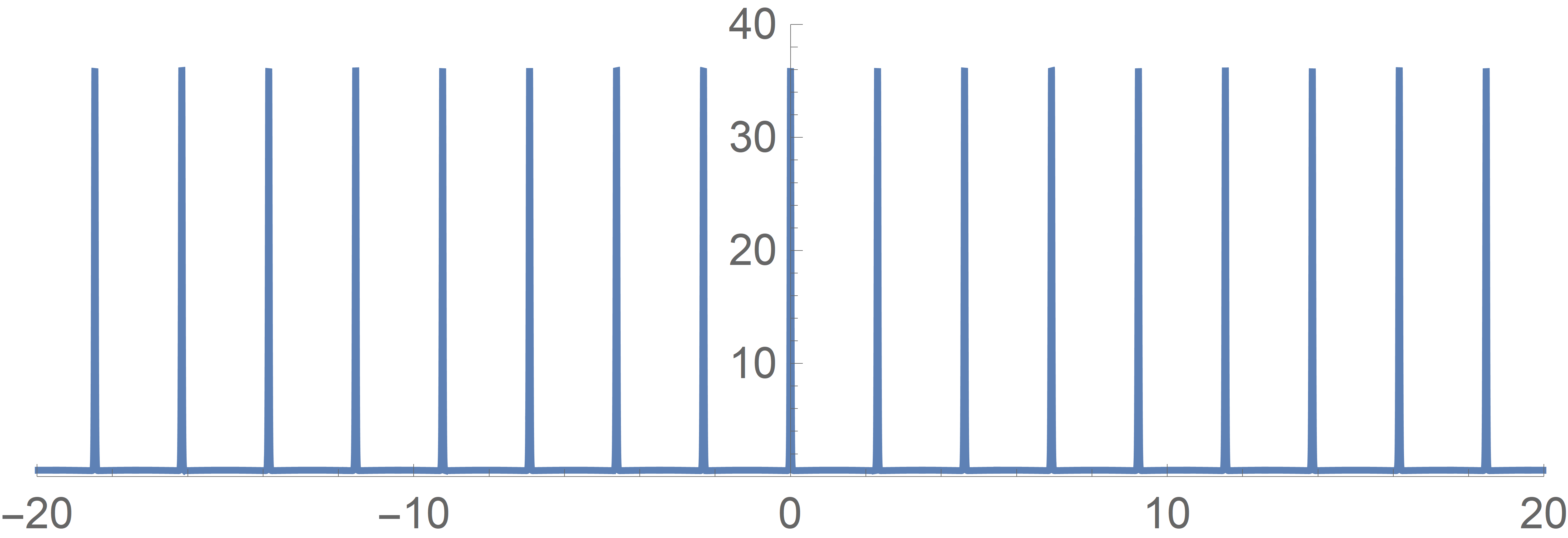}
        \put(-455,44){(c)}
        \put(-302,44){(d)}
        \put(-151,44){(e)}
        \\
\vskip0.22truein 
        \includegraphics[width=0.31\textwidth]{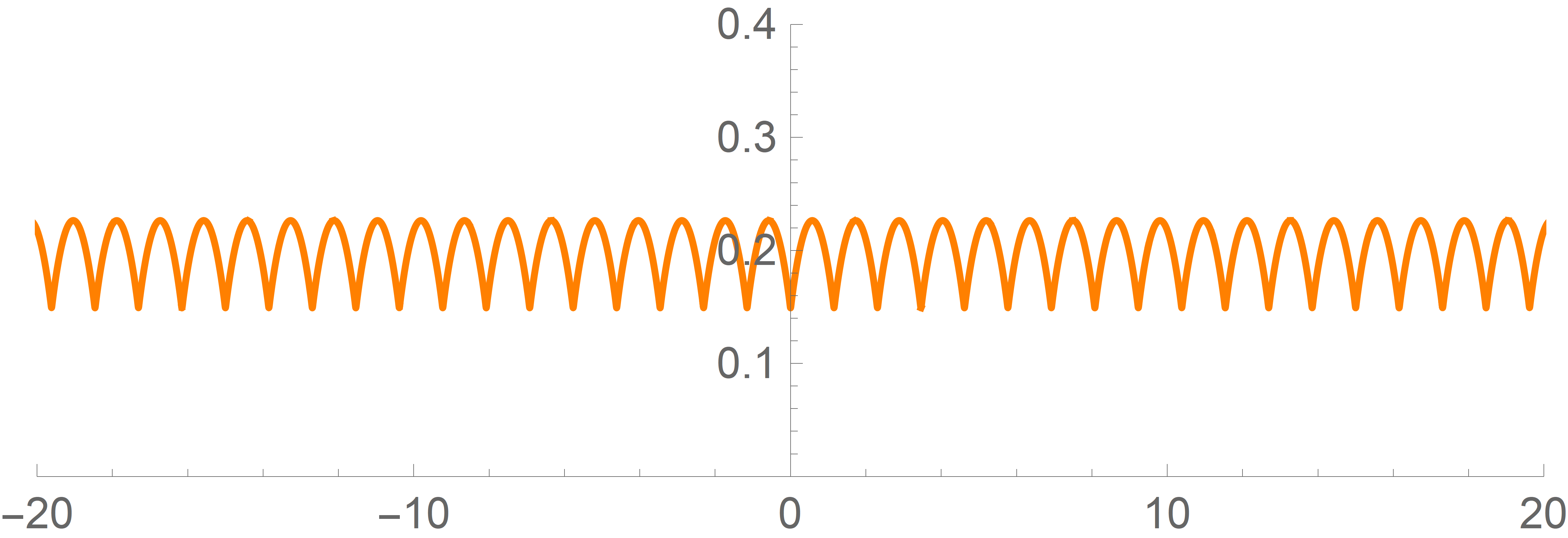}
        \includegraphics[width=0.31\textwidth]{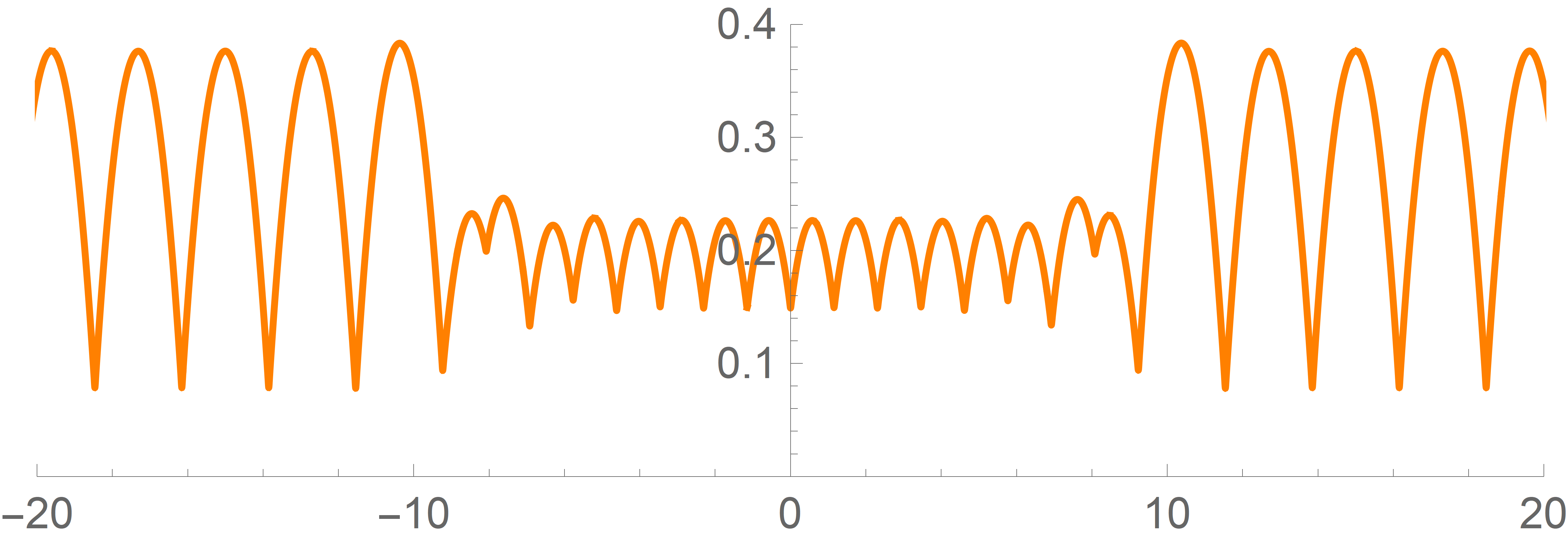}
        \includegraphics[width=0.31\textwidth]{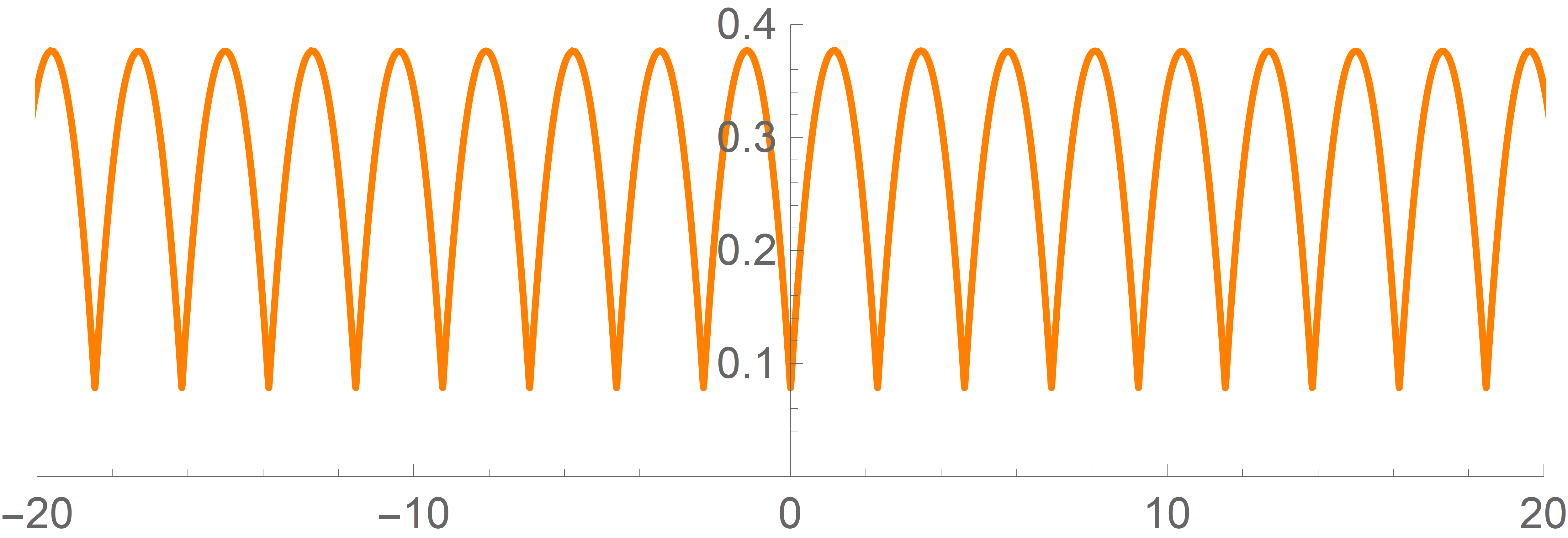}
            \put(-455,50){(f)}
        \put(-302,50){(g)}
        \put(-151,50){(h)}
\caption{The evolution to a stable spatially-periodic canard pattern with $T = 30/13$ ($k = 13 \pi/15$). 
(a) Contour plot of $v(x,t)$.
The (sinusoidal) initial data with $T_{\rm i} = 15/13$ ($k = 26 \pi/13$), {\it i.e.,} with 52 periods on $[-30,30]$, evolves rapidly (before time $t=30$) to a marginally unstable canard pattern with the same period/wavelength. (The subscript $i$ refers to initial.)
This marginally unstable pattern lies just outside the Busse balloon.
Just after $t=700$, the pattern begins to fall back into the Busse balloon by a side band/period doubling `bifurcation wave' that `invades' the original pattern (cf. \cite{SD2025}).
The invasion is complete before $t=860$.
(b) $(u(x,\bar{t}),v(x,\bar{t}))$ for $\bar{t} = 600$ (orange), $800$ (blue), $1000$ (green).
(c) $u(x,600)$.
(d) $u(x,800)$.
(e) $u(x,1000)$.
(f) $v(x,600)$.
(g) $v(x,800)$.
(h) $v(x,1000)$.
The parameters are $A=1$, $B=1$, and $\eps=0.01$.
We only plotted the sub-interval $[-20,20]$.
}
\label{fig:sim2}
\end{figure}

Next, we show a stable periodic state that evolves from initial data with a wavenumber (just) outside the Busse balloon (see Fig.~\ref{fig:sim2}).
In particular, for the same parameter values and on the same domain as in Fig.~\ref{fig:sim1}, the periodic initial data in Fig.~\ref{fig:sim2} has 52 pulses. 
This marginally unstable pattern evolves back into the Busse balloon by an invasion wave that eliminates half of the pulses (Fig.~\ref{fig:sim2}(c)-(h)). 
Although the destabilization process has the nature of a (co-dimension 1) period doubling instability \cite{RS2007}, it may be expected that the boundary of the Busse balloon is determined by a more standard sideband instability.
Indeed, the sideband and period doubling curves (in wavenumber-parameter space) can be close to each other for singularly perturbed R-D systems of this type, as for the Busse balloon in the Klausmeier-Gray-Scott model that is similar to \eqref{Brusselator}, see \cite{SSERDR014}. 
Simulations with other $B$ and initial data with multiple pulses and wavenumbers outside the Busse balloon show similar results. 

We also performed simulations in which the parameter $B$ is slowly varied in time.
These simulations provide insight into the structure of the Busse balloon and into the role played by the multi-pulse canard patterns, where new spikes are generated by the RFSN-II and RFS singularities and their canards, as shown for example in Figs.~\ref{fig:kbelowkT} and \ref{fig:selfsimilarB1.0175}.

\begin{figure}[h!t]
        \includegraphics[width=0.44\textwidth]{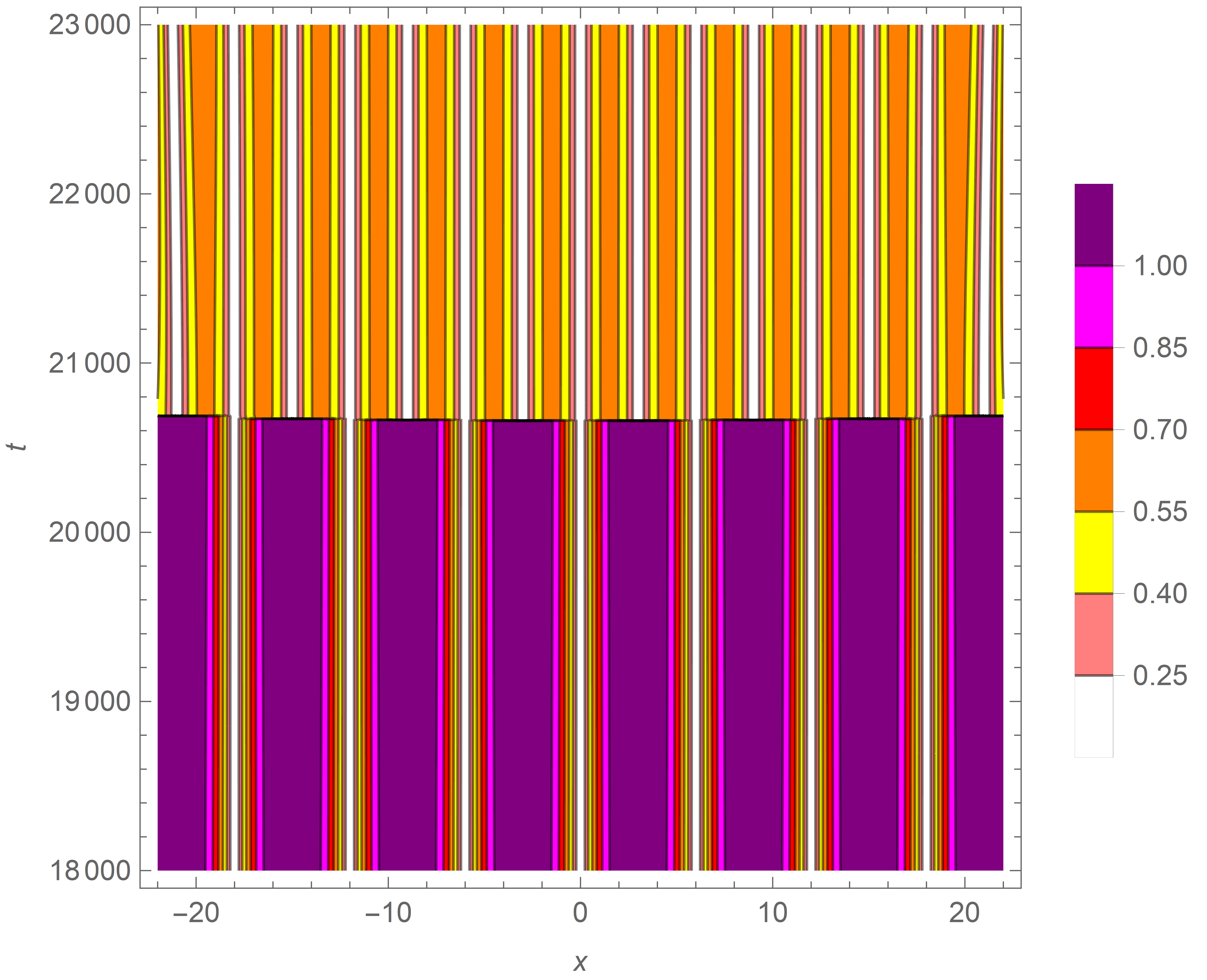} 
        \hspace{0.2truein} 
        \includegraphics[width=0.44\textwidth]{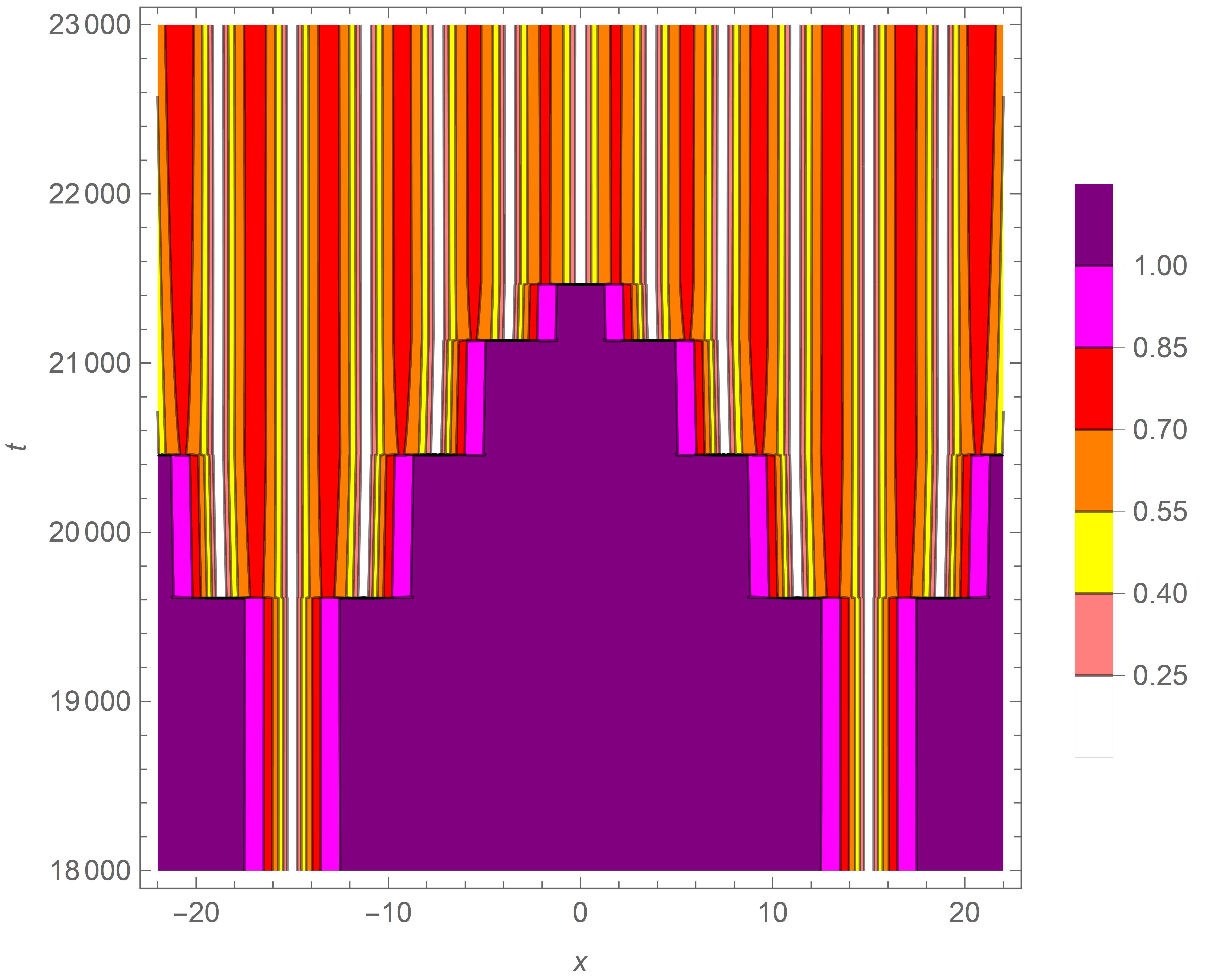}
        \put(-444,158){(a)}
        \put(-220,158){(b)}
\\
        \includegraphics[width=0.42\textwidth]{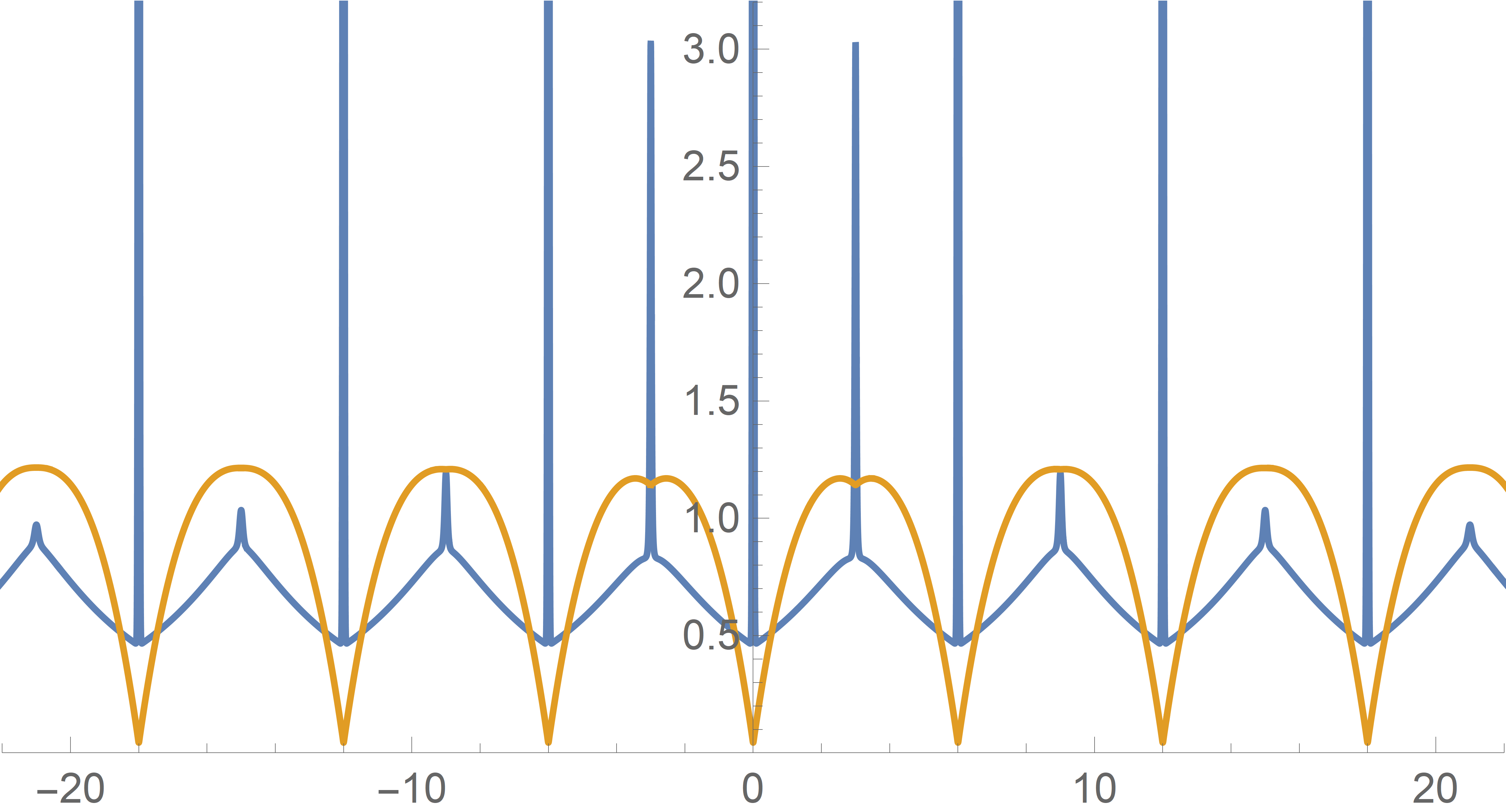}
\hspace{0.25truein}
        \includegraphics[width=0.42\textwidth]{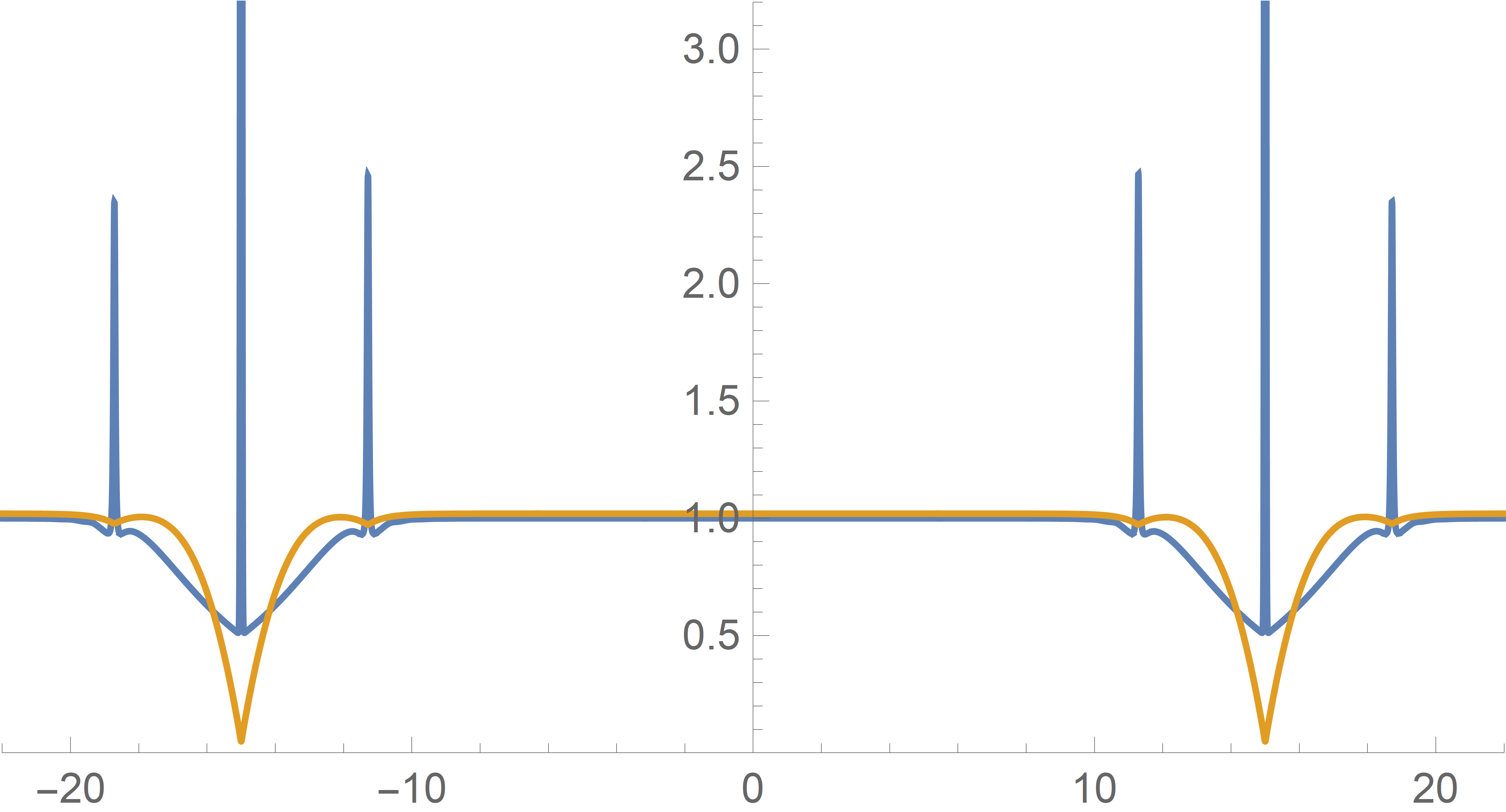}
        \put(-428,100){(c)}
        \put(-205,100){(d)}
\caption{The evolution of spatially-periodic canard patterns of \eqref{Brusselator} with a slowly increasing function $B(t)=1+10^{-5}t$.
(a) Contour plot of $v(x,t)$ from initial data with $T_{\rm i}=6$. 
The initial 10-pulse pattern (on $[-30,30]$) passes through the boundary of the Busse balloon near $t = 20660$, {\it i.e.,} at $B \approx 1.2066$.
A period-halving bifurcation occurs across the domain.
(b) Contour plot of $v(x,t)$ from initial data with $T_{\rm i}=30$.
The initial two-pulse data gains pulses in an almost stepwise manner after it passes through the boundary of the Busse balloon near $B \approx 1.1961$.
Further steps occur near $B \approx 1.2045$, $B \approx 1.2113$, and $B \approx 1.2146$. 
(c) $u(x,20659.5)$ (blue) and $v(x,20659.5)$ (orange) for the first simulation.
(d) $u(x,19610.5)$ (blue) and $v(x,19610.5)$ (orange) for the second simulation, just at the first step. 
Here, $A=1$ and $\eps=0.01$.
Homogeneous Neumann boundary conditions on $[-30,30]$, plotting only $[-22,22]$.
}
\label{fig:sim4}
\end{figure}

We begin with simulations in which $B$ increases slowly and linearly in time (see Fig.~\ref{fig:sim4}). 
We take $B(0)=1$ and $B'(t) \equiv 10^{-5}$.
In the first column, the initial data for the simulation are sinusoidal with spatial period $T_{\rm i}=6$ and are similar to those in the simulation shown in Fig.~\ref{fig:sim1}. 
In the second column, the initial data has $T_{\rm i}=30$.
It is inside the domain of attraction of the long wavelength canard pattern. 
It is a (sharply) peaked combination of $\tanh$ and $|\sin|$ functions: $u(x,0) = 
\{(1/2 - 15 [\tanh(50(x + 15) - 1/40)) - \tanh(50(x + 15) + 1/40)]\}/2 +
\{(1/2 - 15 [\tanh(50(x - 15) - 1/40)) - \tanh(50(x - 15) + 1/40)]\}/2
+ \{|\sin (\pi(x + 15)/50)|^{1/4} + |\sin (\pi(x - 15)/50)|^{1/4}]\}/4$.
The contour plots of $v(x,t)$ are shown in Figs.~\ref{fig:sim4}(a) and (b).
In both, the transient patterns and the final attractors are spatially-periodic canard patterns of the type constructed in Sec.~\ref{s:geoconstruction}.

Fig.~\ref{fig:sim4} also shows that the evolution to the attractors can take place via two routes.
In (a), period-halving occurs across the domain in a brief interval of time near $t\approx 20660$, and then the steady stare attractor is reached. 
Alternatively, as illustrated in (b), patterns gain pulses in a step-wise fashion, after $B$ passes through the boundary of the Busse balloon, in this case near $B \approx 1.1961$.
These added pulses invade the domain until the steady state attractor is reached.

\begin{figure}[h!tbp]
\centering
        \includegraphics[width=0.31\textwidth]{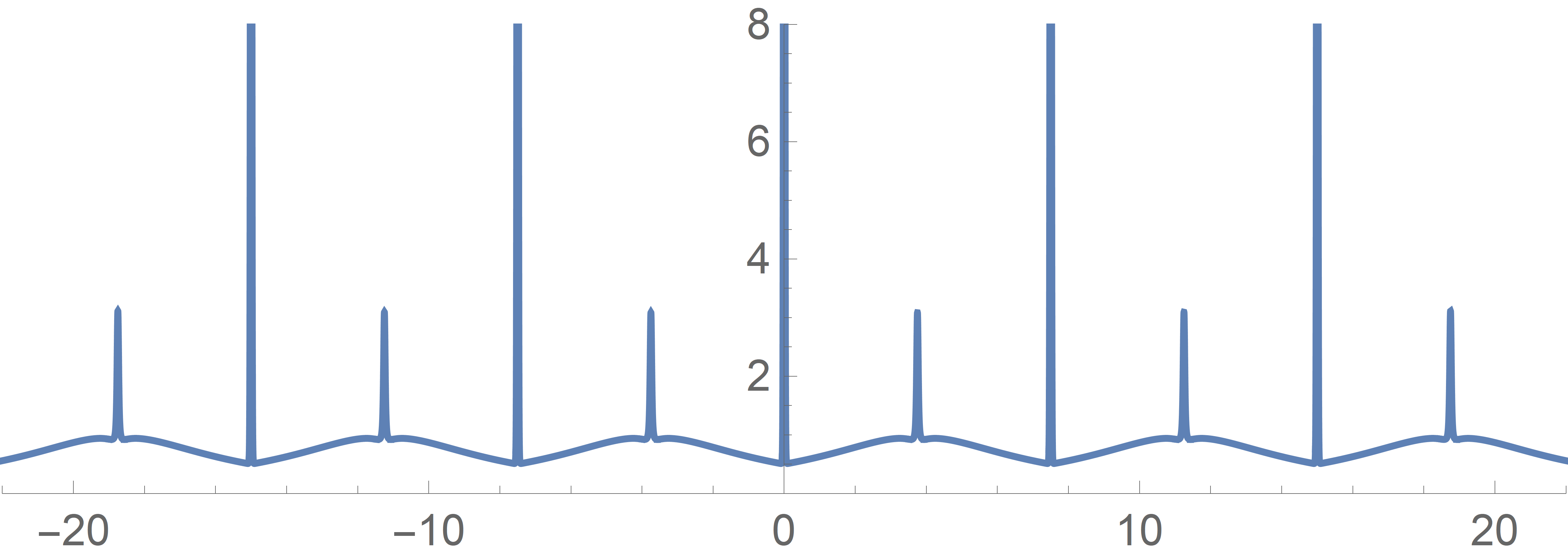} \, \, 
        \includegraphics[width=0.31\textwidth]{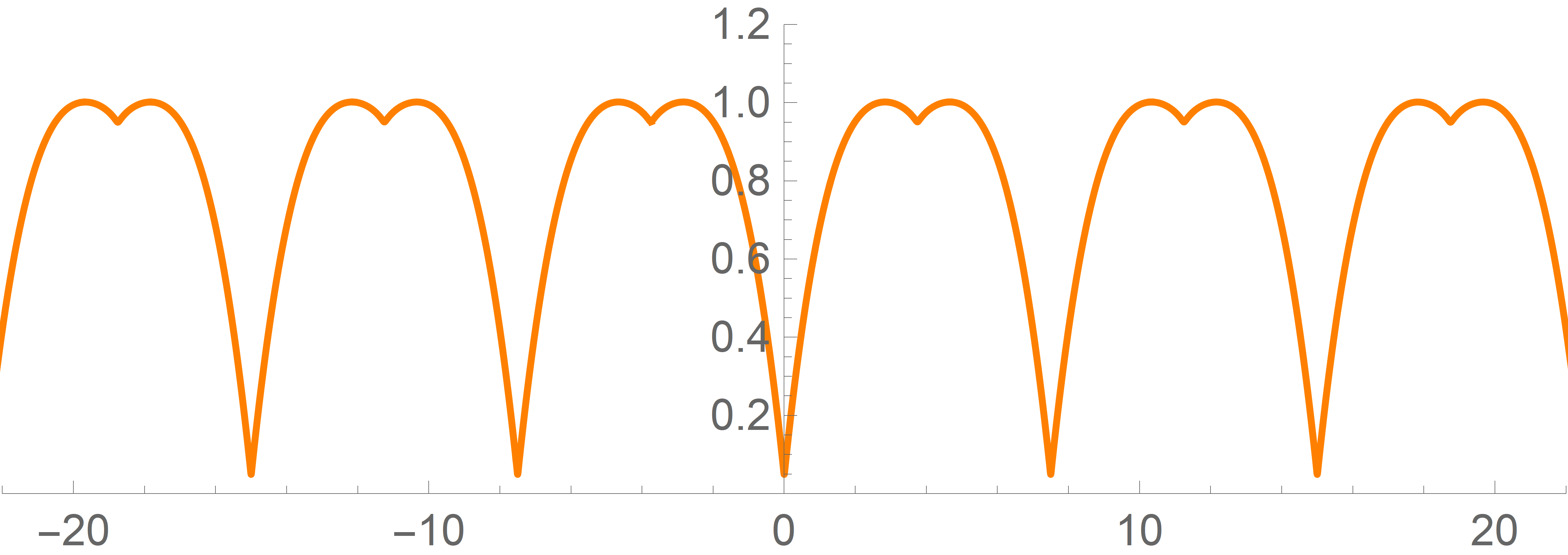} \, \, 
        \includegraphics[width=0.31\textwidth]{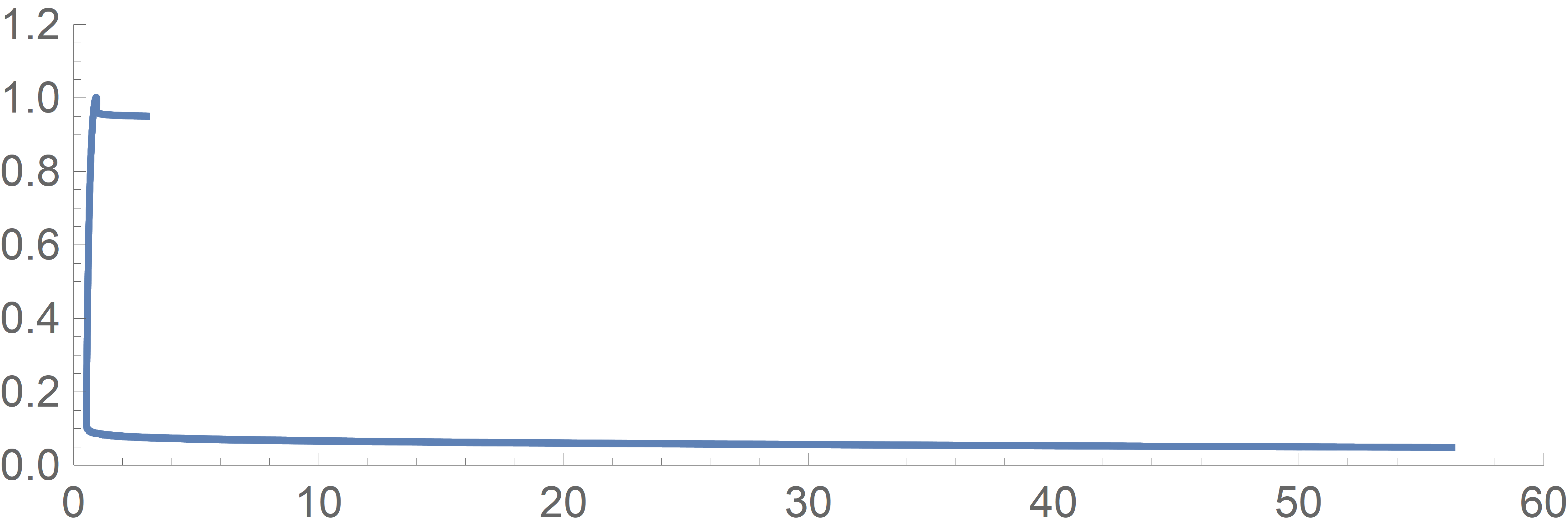}
        \put (-475,50){(a)}
        \put(-312,50){(b)}
        \put(-162,50){(c)}
        \\
\centering
        \includegraphics[width=0.31\textwidth]{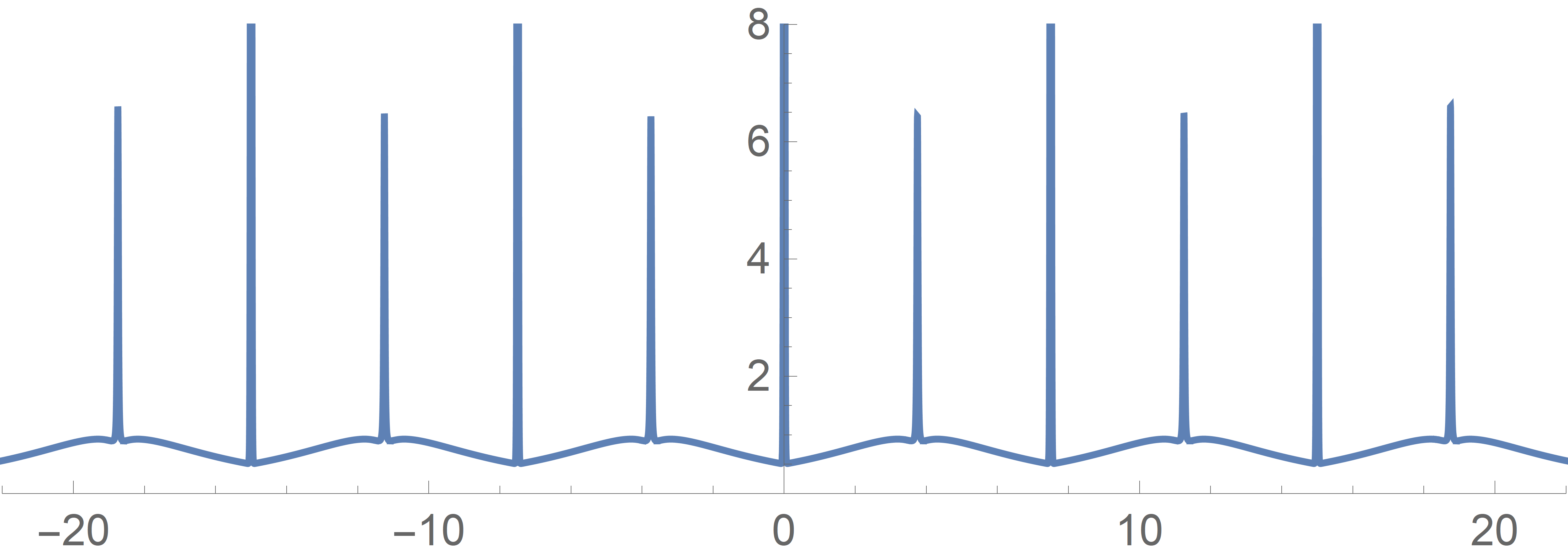} \, \, 
        \includegraphics[width=0.31\textwidth]{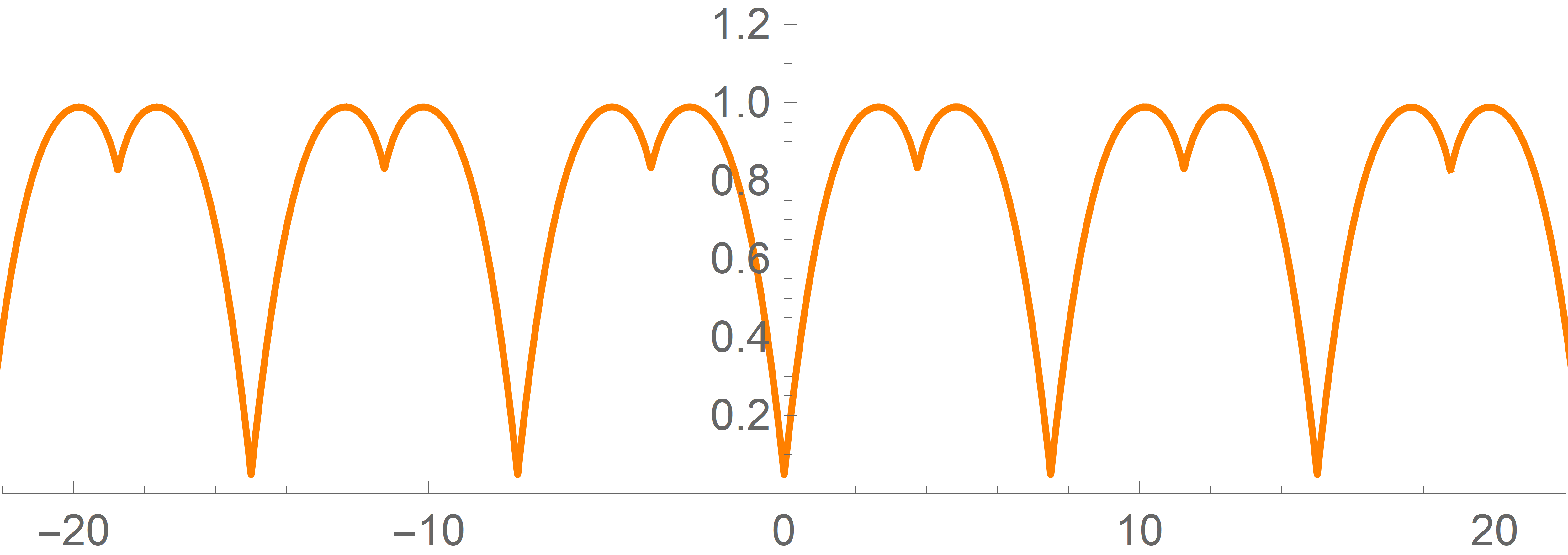} \, \, 
        \includegraphics[width=0.31\textwidth]{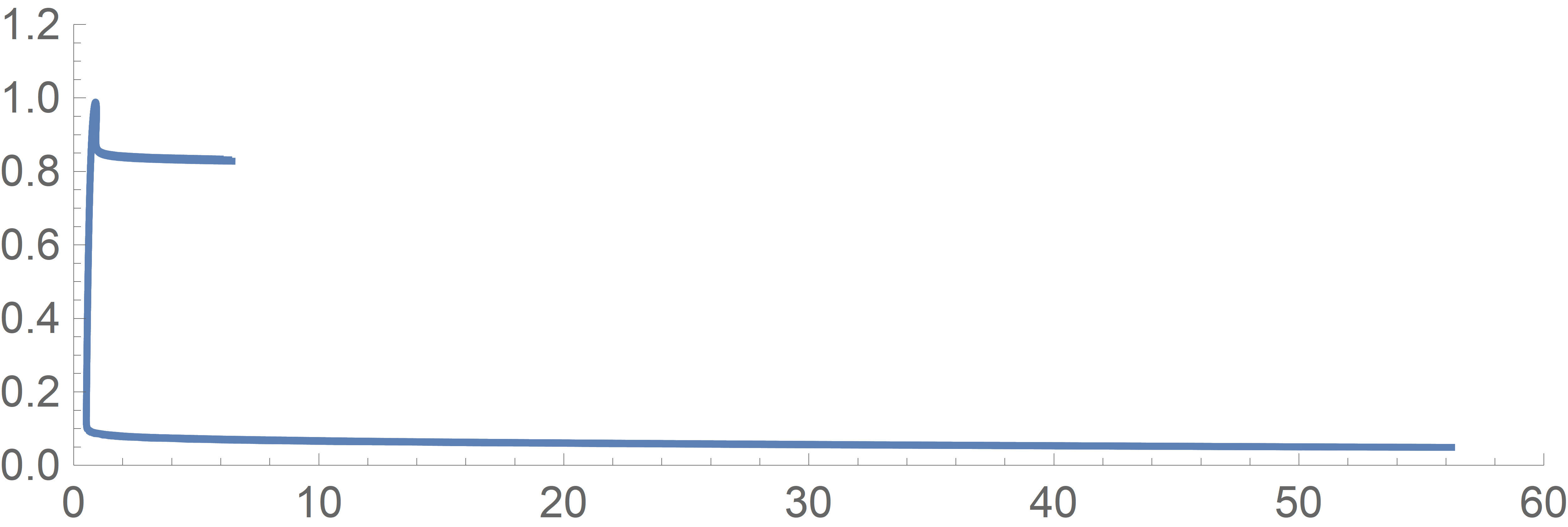}
        \put(-475,50){(d)}
        \put(-312,50){(e)}
        \put(-162,50){(f)}
        \\
\centering 
        \includegraphics[width=0.31\textwidth]{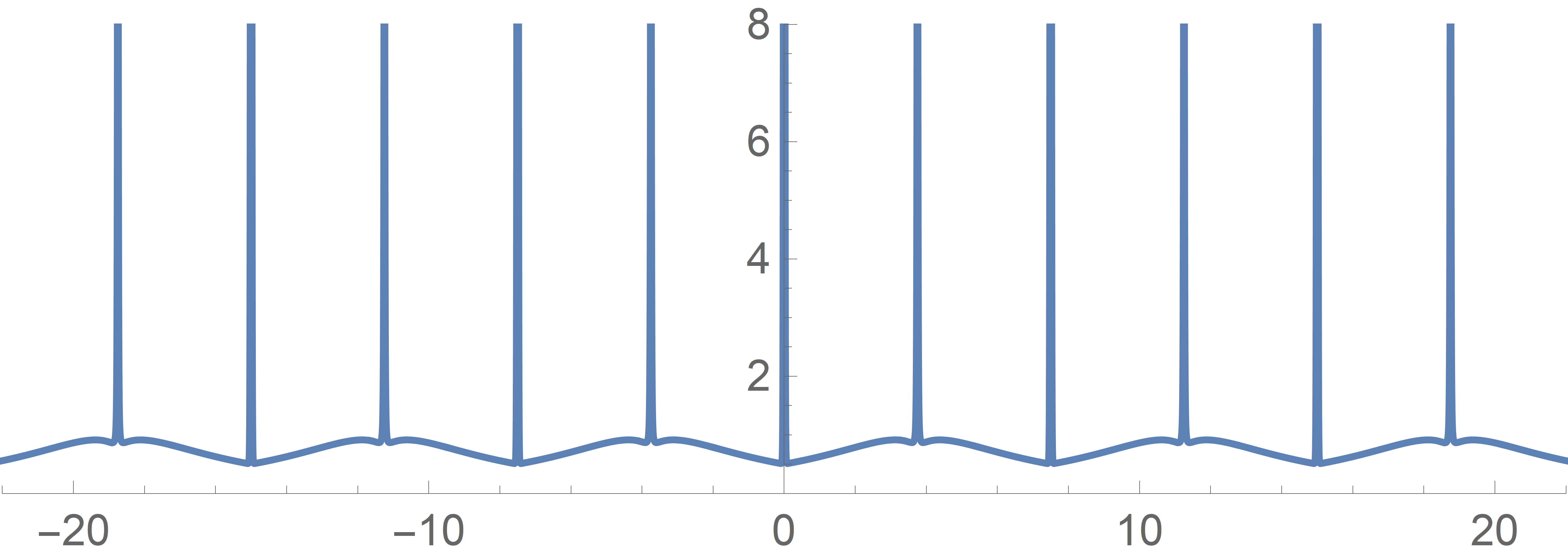} \, \, 
        \includegraphics[width=0.31\textwidth]{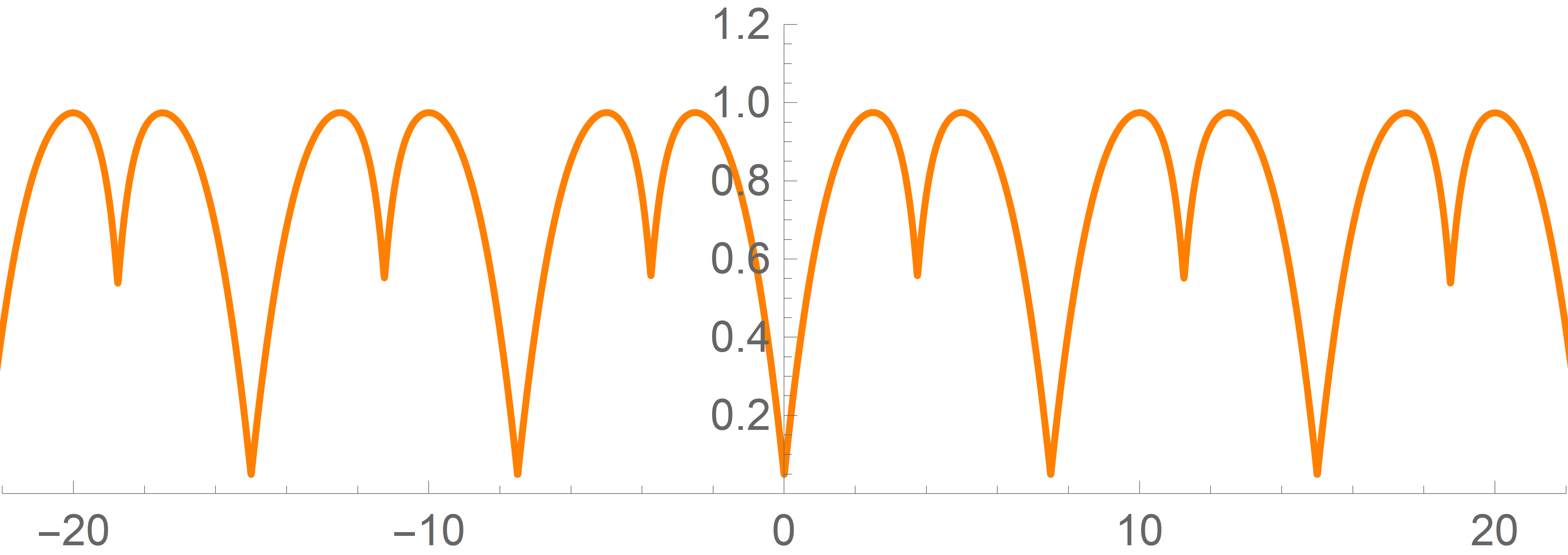} \, \, 
        \includegraphics[width=0.31\textwidth]{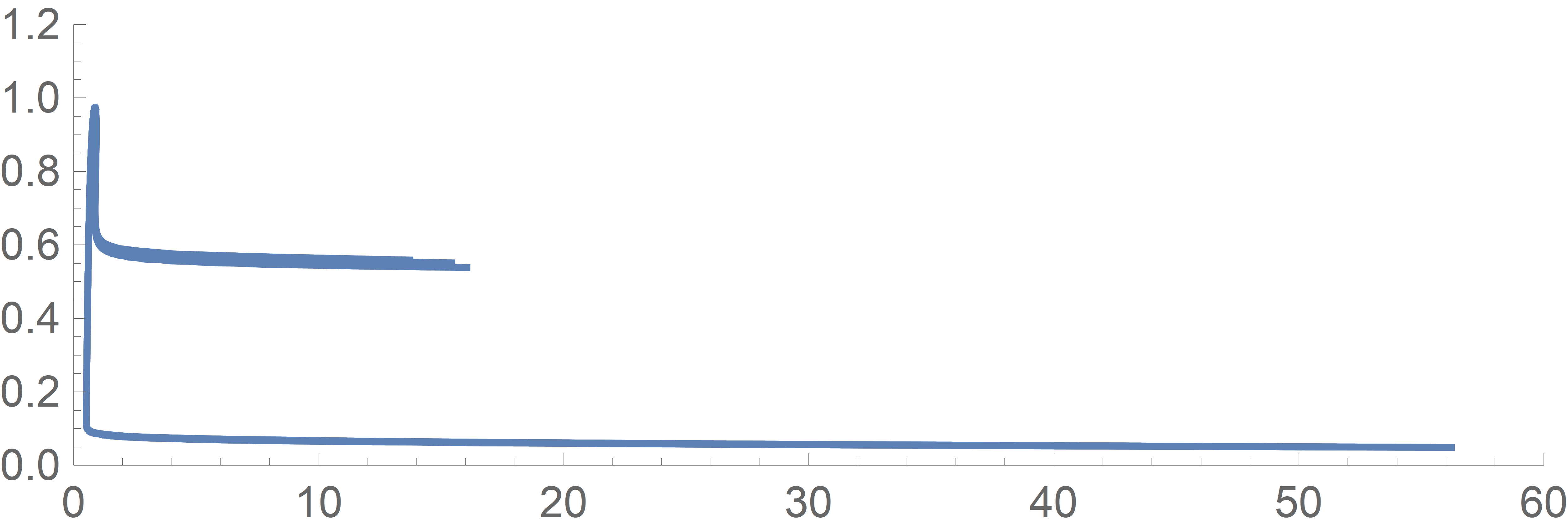}
         \put(-475,50){(g)}
        \put(-312,50){(h)}
        \put(-162,50){(i)}
\caption{The (fast) growth of spikes in \eqref{Brusselator} near $B \approx 1.0218$ at the transition from an 8-pulse pattern to a 16-pulse pattern on the interval $[-30,30]$, obtained with  $B(t) = 1 + 10^{-5}t$ and homogeneous Neumann boundary conditions.
(a) $u(x,2017.4)$; 
(b) $v(x,2017.4)$;
(c) $(u(x,2017.4),v(x,2017.4))$.
(d) $u(x,2017.8)$;
(e) $v(x,2017.8)$;
(f) $(u(x,2017.8),v(x,2017.8))$.
(g) $u(x,2018.0)$;
(h) $v(x,2018.0)$;
(i) $(u(x,2018.0),v(x, 2018.0))$). 
Here, $A=1$ and $\eps=0.01$, and we plot the interval $x \in [-22,22]$. 
}
\label{fig:sim5}
\end{figure}

For the evolution through period-halving (Figs.~\ref{fig:sim4}(c) and \ref{fig:sim5}), the formation of new pulses is guided by the true and faux canards of the folded singularity as $B$ passes through the Turing point.
The new pulses are formed exactly at the locations of the local maxima of $u(x,t)$ between adjacent existing pulses.
The tiny initial spikes are of the form created by the canards of the RFSN-II and RFS singularities, and observed, for example, in the green pattern shown in Fig.~\ref{fig:kbelowkT}. 
As time increases, the solution of the PDE evolves through a family of these (apparently marginally unstable) canard patterns, with the magnitudes of the new pulses slowly increasing. 

\begin{figure}[h!tbp]
        \includegraphics[width=0.44\textwidth]{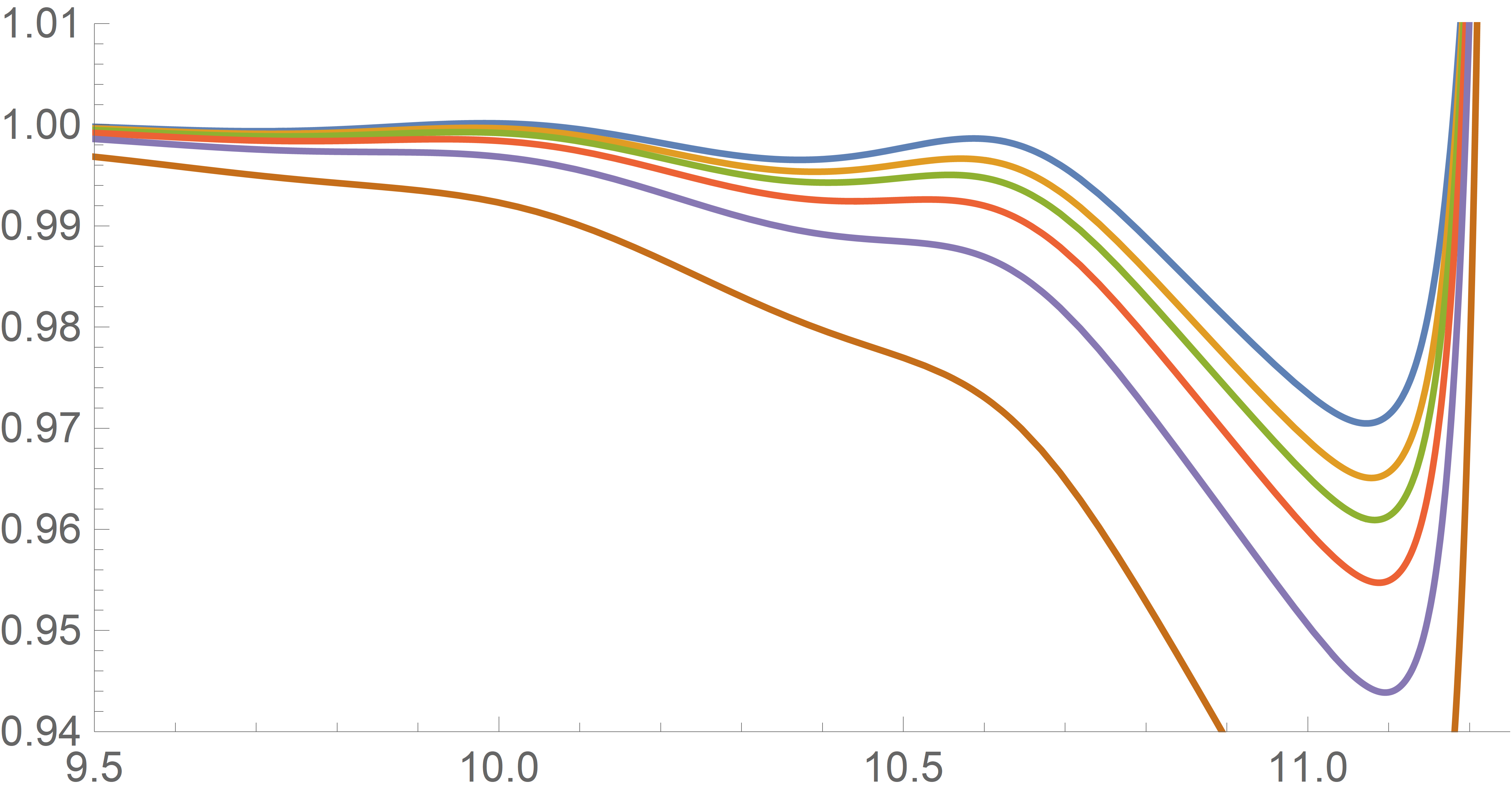} 
        \hspace{0.2truein} 
        \includegraphics[width=0.44\textwidth]{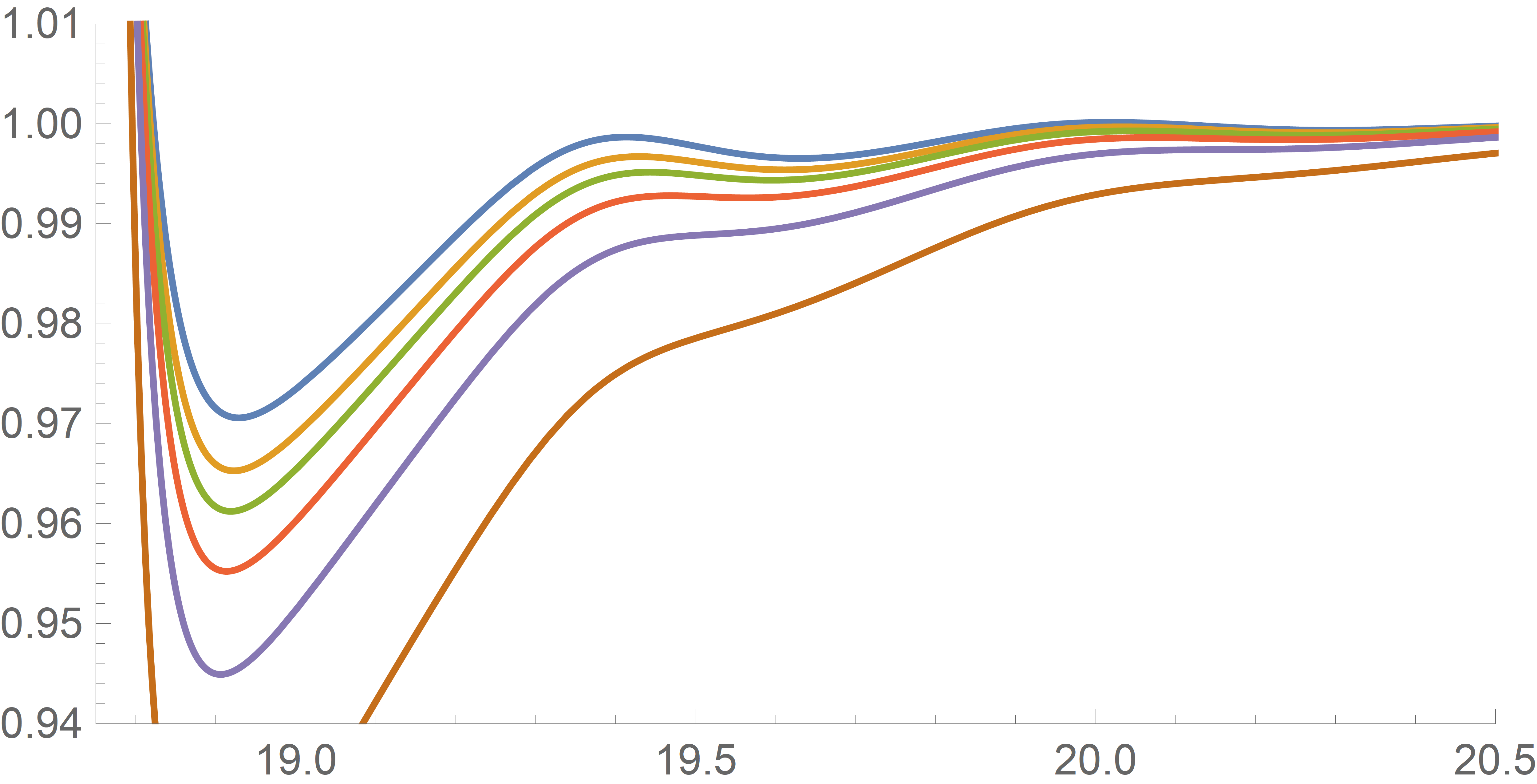}
        \put(-450,98){(a)}
        \put(-225,98){(b)}
\caption{The fast evolution of (spatial) oscillations on the tails of the new spikes appearing on the wings of the original pulse at $x = 15$ during the first fast transition in the simulation shown in Figs. \ref{fig:sim4}(b).
(a) zoom of the $u$-component of the left new spike near $x \approx 11.27$.
(b) zoom of the right new spike near $x \approx 18.71$, at seven equidistant moments in  time ($t \in [19605, 19611] \supset 19610.5$ of Fig. \ref{fig:sim4}(b)). 
The `pulse adding' mechanism shown in Figs. \ref{fig:sim4}(b) and (d) is associated to the appearance of spatial oscillations on the `tails' of the new spikes -- similar to the spatial patterns shown in Fig.~\ref{fig:selfsimilarB1.0175}.
}
\label{fig:sim4b}
\end{figure}

For the spatially-periodic attractors created by stepwise pulse addition (Fig.~\ref{fig:sim4} (b) and (d)), the new spikes appear on the `wings' of the original pulses, once the distance between these pulses is sufficiently large. 
Details of the fast evolution of (spatial) oscillations on the wings are shown in Fig.~\ref{fig:sim4b}.
The new spikes are nucleated where $v$ has a local minimum on either side of the overshoot, {\it i.e.,} near $v=1$. 
As the transients develop, the patterns have  increasingly longer segments along the canards of the folded singularity.
They are of the type constructed in Sec.~\ref{s:geoconstruction}.

\begin{figure}[h!tbp]
        \includegraphics[width=0.325\textwidth]{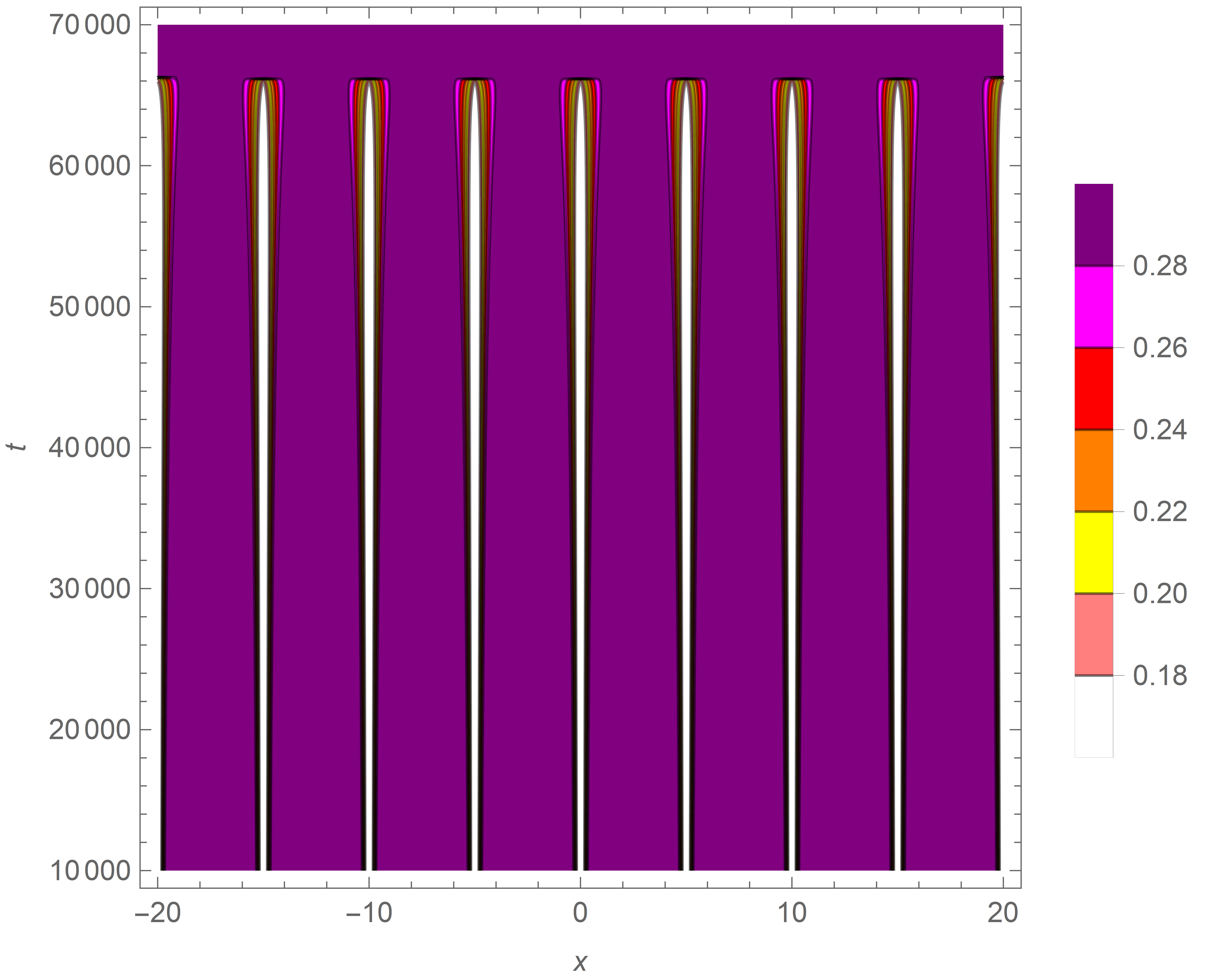}
        \includegraphics[width=0.325\textwidth]{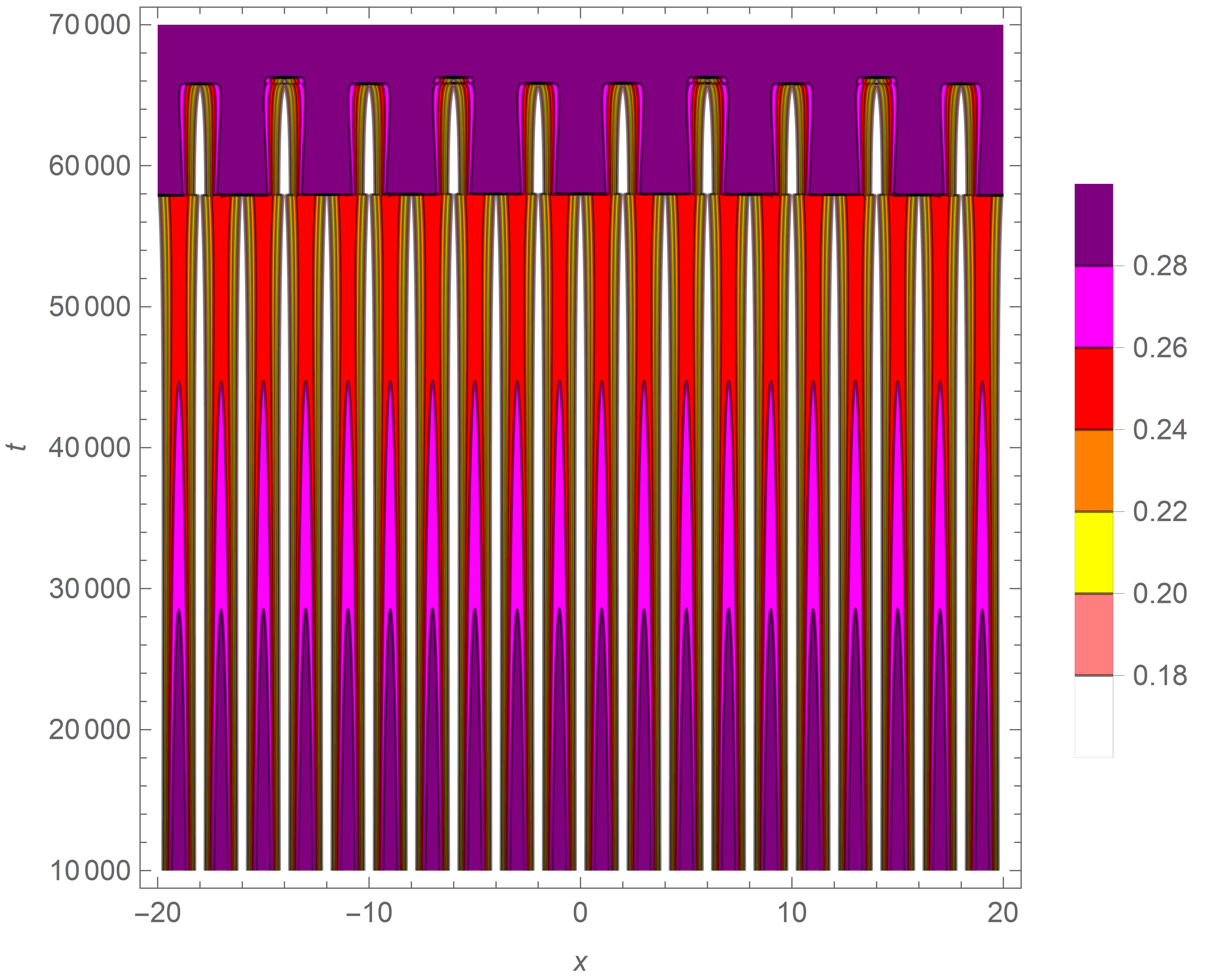}
        \includegraphics[width=0.325\textwidth]{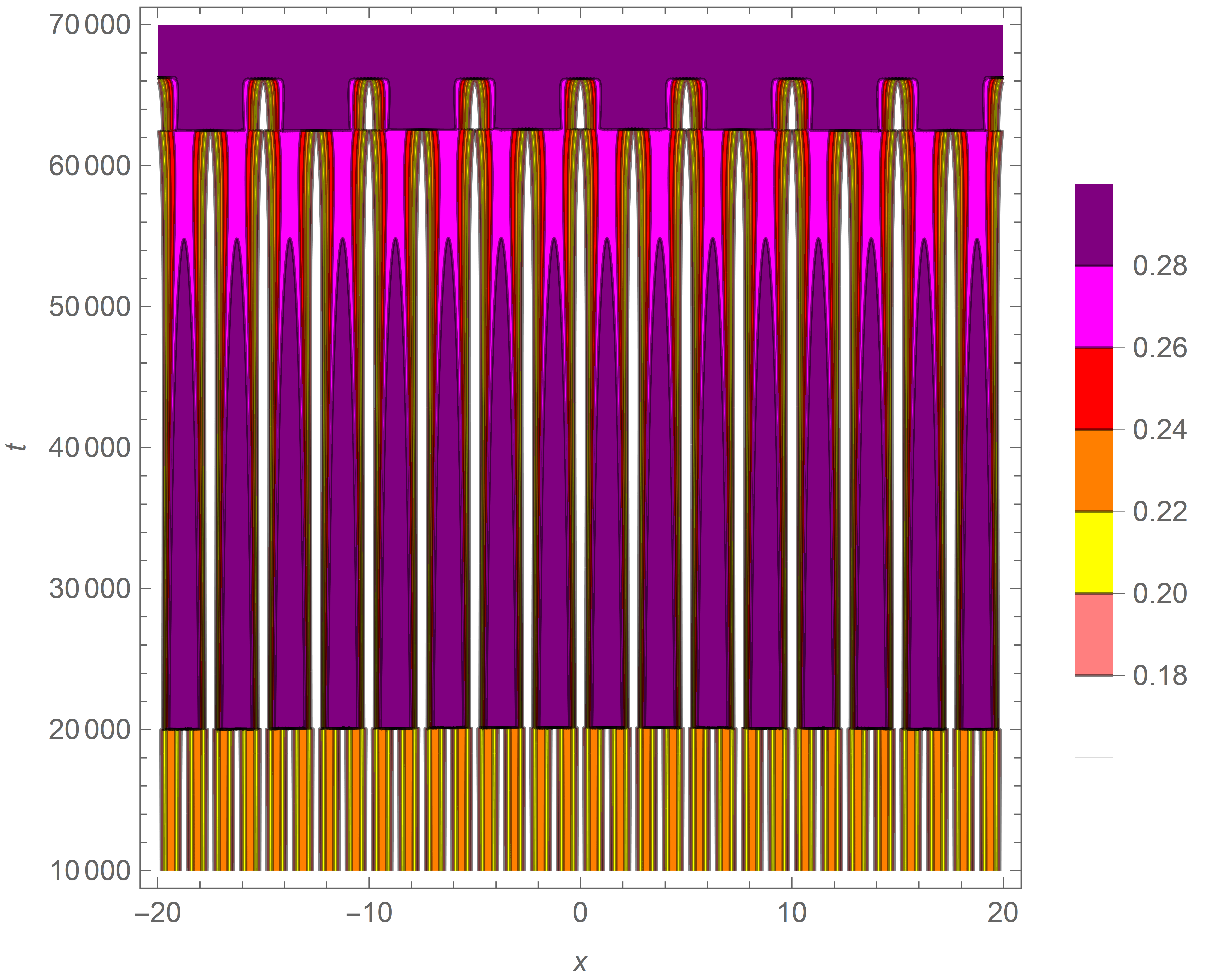}
        \put(-465,110){(a)}
        \put(-309,110){(b)}
        \put(-151,110){(c)}
\caption{Contour plots of $v(x,t)$ for the case of  $B$ a slowly, linearly decreasing function of time.
Simulations for three different spatially-periodic initial data, with $T_{\rm i}=5$ in (a), $T_{\rm i}=2$ in (b), and $T_{\rm i}=5/4$ in (c). 
Each pattern samples part of the Busse balloon, and some undergo period-doubling.
Then, for some $t$ near $66,000$, $B(t)$ crosses the saddle-node bifurcation at $B \approx 0.34$ at which the periodic orbits disappear, so that only the homogeneous state is stable.
$B = B(t)$ with $B'(t) \equiv -10^{-5}$ and $B(0)=1$. 
Here, $A=1$, $\eps=0.01$, and the domain is $[-30,30]$ with homogeneous Neumann boundary conditions. 
}
\label{fig:sim3}
\end{figure}

To complement the study with $B(t)$ increasing, we also consider the case in which the parameter $B$ is a slowly and linearly {\it decreasing} function of time.
We take $B(t) = 1 - 10^{-5}t$.
In Fig.~\ref{fig:sim3}, we show results for solutions generated by three different initial data of the same type as used in Fig.~\ref{fig:sim1}, with spatial periods $T_{\rm i} = 5, 3$, and $5/4$, respectively. 

The analysis in Sec.~\ref{s:TH} shows that substantial portions of the branches of spatially periodic patterns terminate in saddle-node (s-n) bifurcations
(recall Fig.~\ref{fig:bifurcation3d}). 
In the three simulations shown in Fig.~\ref{fig:sim3}, the s-n bifurcation occurs near $t = 66,000$, {\it i.e.,} at $B \approx 0.34$, after which the solutions become homogeneous. 
This agrees well with the critical value $B_F(1) \approx 0.34$ that is obtained both from the surface (plotted numerically in Fig.~\ref{fig:bifurcation3d}) of spatially-periodic canard patterns in the spatial ODE system \eqref{spatialODE} and from the direct PDE simulations on domains of length $L=2\pi/k$ in Fig.~\ref{fig: H1VsB}.
For the simulation in Fig.~\ref{fig:sim3}(a), the wavenumber does not change, $k= 2 \pi/5 \approx 1.2$ ({\it i.e.}, there are eight full periods --and seven full inner pulses-- within the interval $[-20,20]$), until the s-n bifurcation occurs. 
The same holds for the simulation in Fig.~\ref{fig:sim3}(c), since also here the wavenumber is $k \approx 1.2$ after two period-doublings and until just before the s-n bifurcation. 

Although it is more complex, the pattern in Fig.~\ref{fig:sim3}(b) also involves collapse of the pulses via a s-n bifurcation.
There is an interplay between a period-doubling bifurcation near $t = 65,800$, from which four (non-equidistant) pulses appear (out of the 10 that persisted through the preceding period doubling near $t = 58,000$) and the final s-n collapse near $t=66,300$.
The period-doublings exhibited by the patterns in Figs.~\ref{fig:sim3}(b) and (c) are typical for spatially-periodic patterns driven by a slowly varying parameter (see \cite{AACS2025, SSERDR014, SD2025} and the references therein). 

\begin{figure}[h!t]
\centering
        \includegraphics[width=0.31\textwidth]{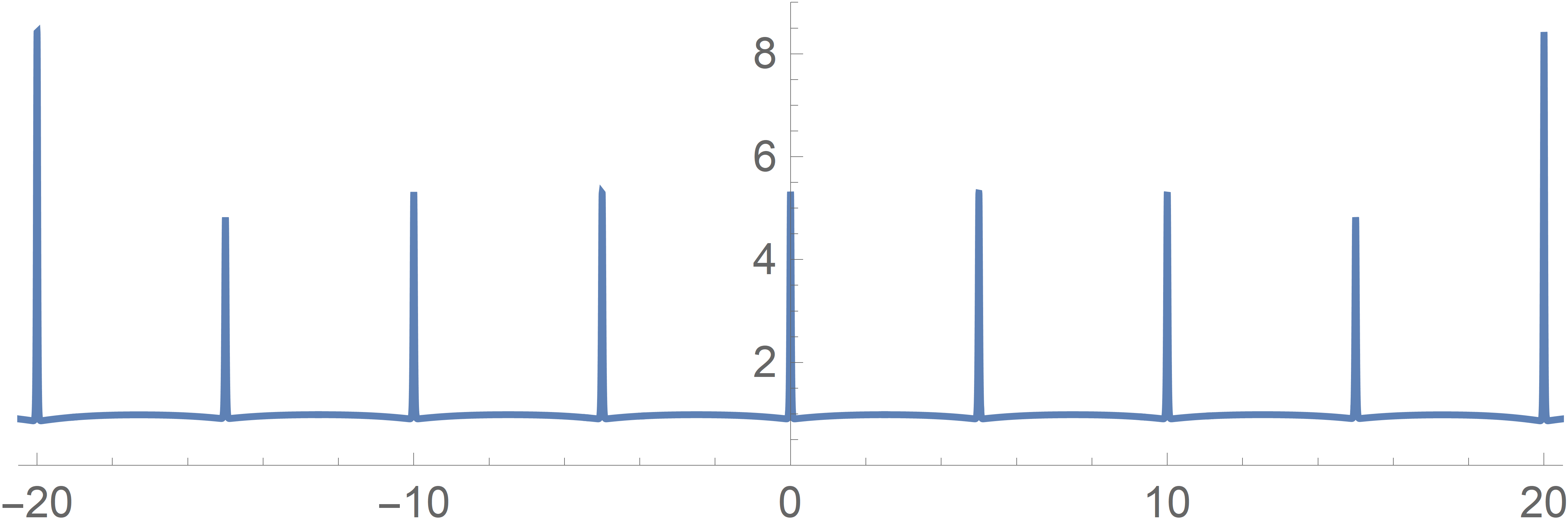} \, \, 
        \includegraphics[width=0.31\textwidth]{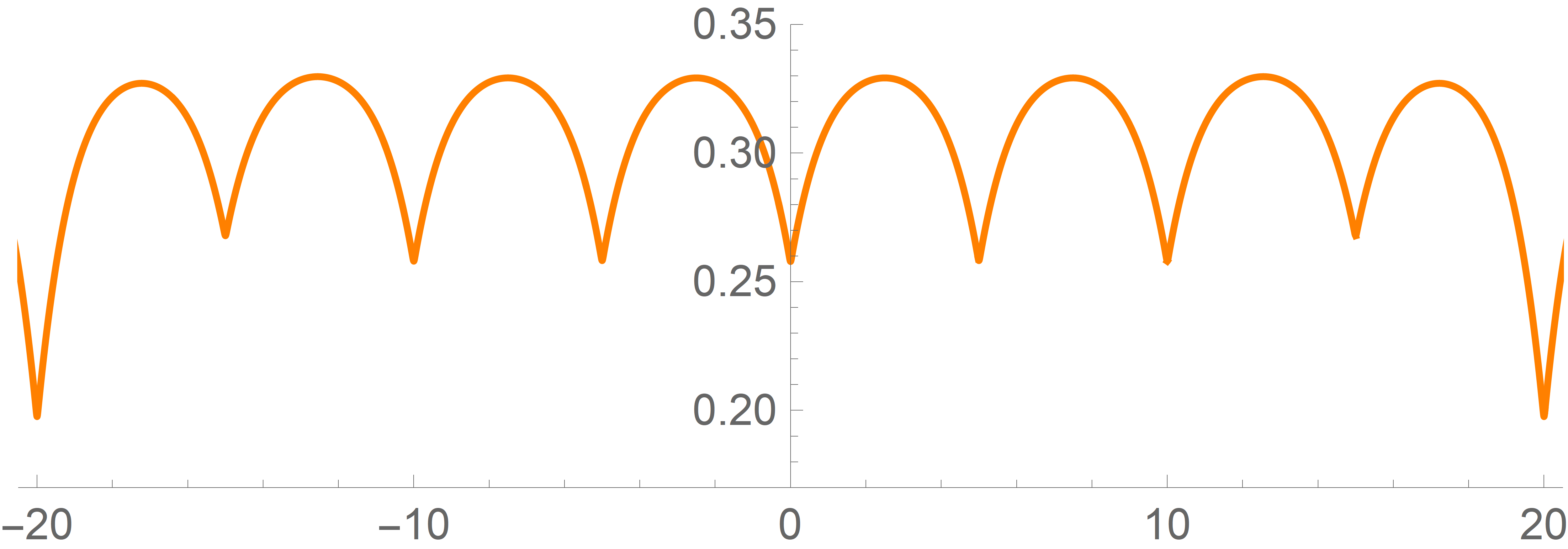} \, \, 
        \includegraphics[width=0.31\textwidth]{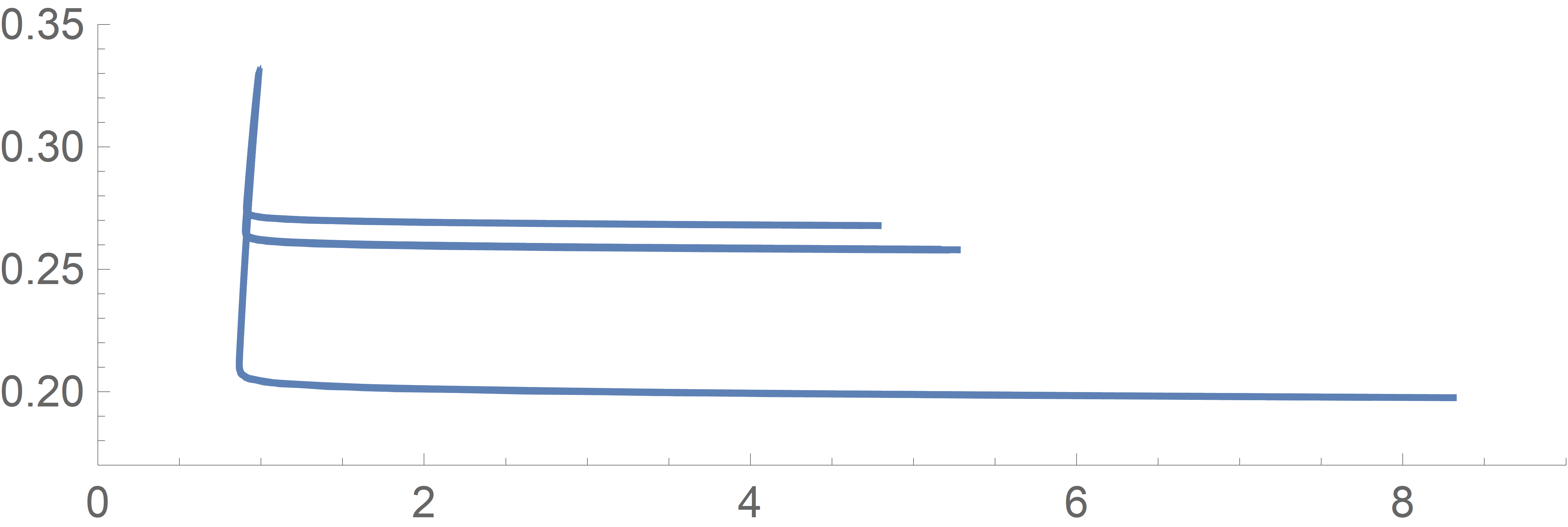}
        \put (-482,40){(a)}
        \put(-315,40){(b)}
        \put(-161,40){(c)}
        \\
\centering
        \includegraphics[width=0.31\textwidth]{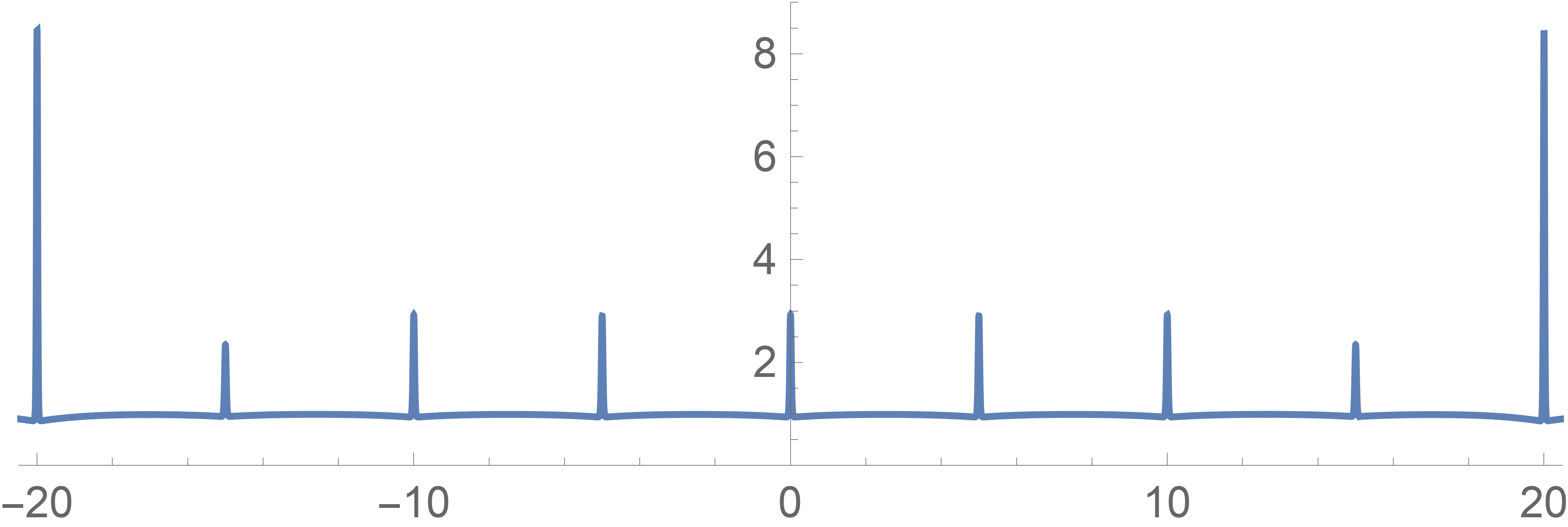} \, \, 
        \includegraphics[width=0.31\textwidth]{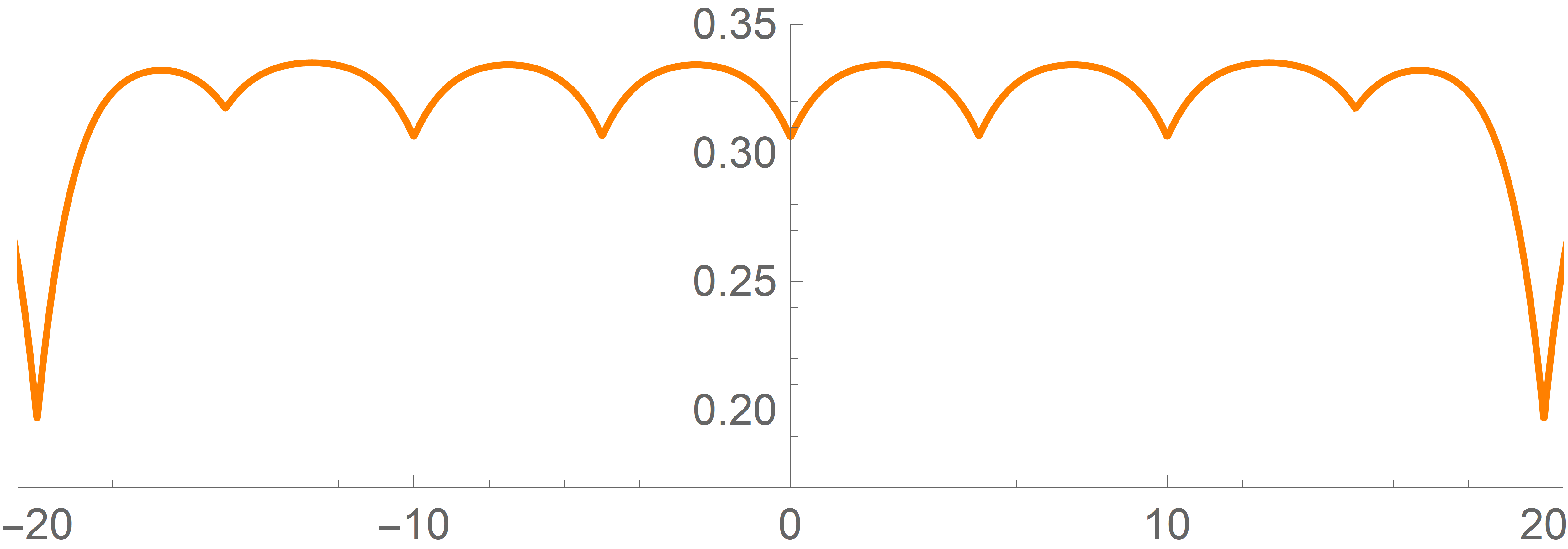} \, \, 
        \includegraphics[width=0.31\textwidth]{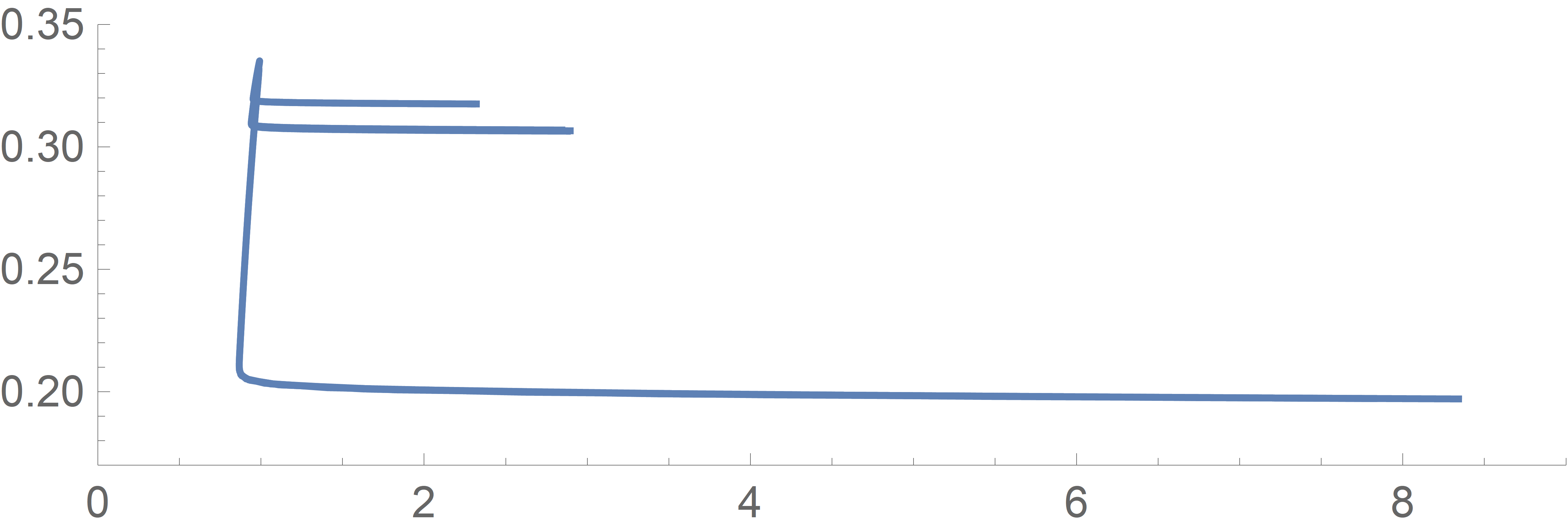}
        \put(-482,40){(d)}
        \put(-315,40){(e)}
        \put(-161,40){(f)}
        \\
\centering 
        \includegraphics[width=0.31\textwidth]{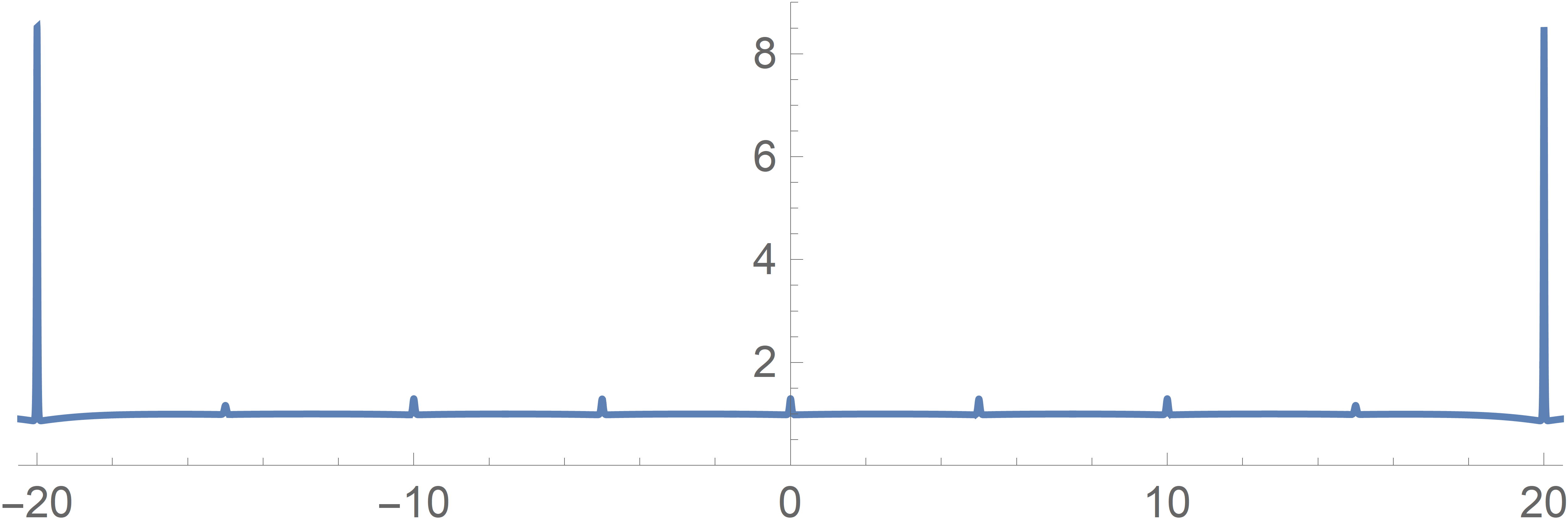} \, \, 
        \includegraphics[width=0.31\textwidth]{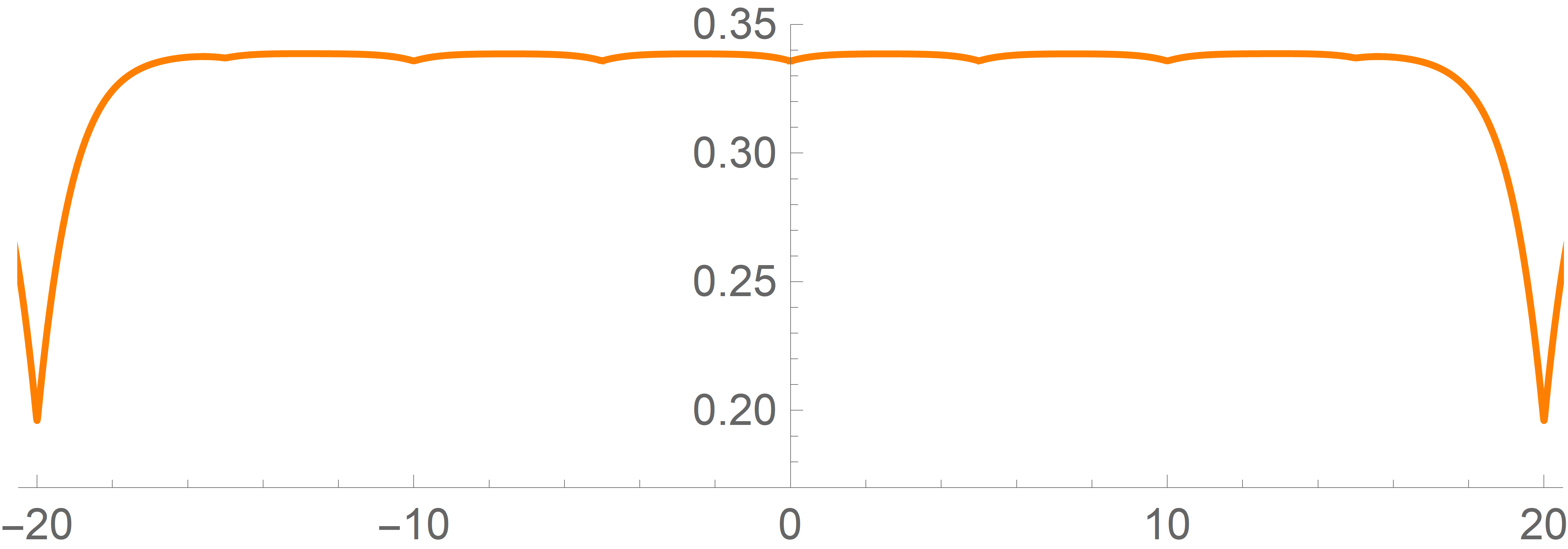} \, \, 
        \includegraphics[width=0.31\textwidth]{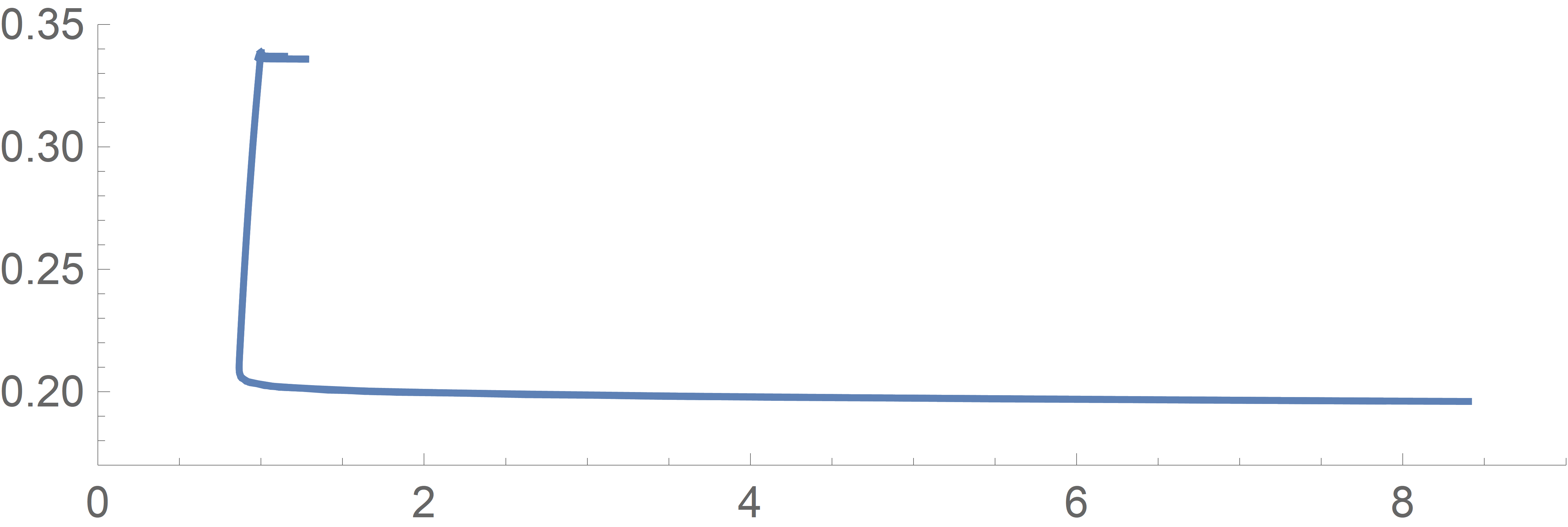}
         \put(-482,40){(g)}
        \put(-315,40){(h)}
        \put(-161,40){(i)}
\caption{
Transient patterns observed during the final stage of the collapse of the periodic state to the homogeneous state for the pattern with $T_{\rm i} = 5$ shown in Fig.~\ref{fig:sim3}(a). 
In the first, second, and third rows, we show snapshots of $u(x,\bar{t})$, $v(x,\bar{t})$, and $(u(x,\bar{t}),v(x,\bar{t}))$ for $\bar{t} = 66159$, $66162$, and $66165$, respectively. 
The transient patterns have the form of the spatially-periodic canards constructed in Sec.~\ref{s:geoconstruction}.
Similar dynamics occur for the transitions in Figs.~\ref{fig:sim3}(b) and (c).
}
\label{fig:sim3b}
\end{figure}

Hence, for $B(t)$ slowly decreasing, the dynamics are guided by multi-pulse canard patterns of the types constructed in Sec.~\ref{s:geoconstruction}, just as for the opposite case of slowly increasing $B(t)$. 
We observe that, as $B(t)$ slowly decreases and the patterns evolve to steady state, they are close to a sequence of periodic canards with medium-amplitude and then small-amplitude spikes.
See also Fig. \ref{fig:sim3b}.

Overall, these simulations of the dynamics on large intervals complement the PDE simulations of Section \ref{s:PDE-numerics} on intervals one-period in length.
They provide a fuller picture of the Busse balloon (as expected).
Indeed, the s-n curves shown in Figs.~\ref{fig:bifurcation3d} and \ref{fig: H1VsB} are outside the Busse balloon for patterns with wavenumber $k > k^*$, for a certain $k^* > 1.0$. 
As Figs.~\ref{fig:sim3}(b) and (c) show, patterns with initial wavenumber $k_{\rm i} > k^*$ evolve --via some steps in which they cross through a non-saddle-node part of the Busse balloon boundary and `bounce' back into it-- until the wavenumbers satisfy $k < k^*$.

\section{Conclusions, Discussion, and Open Problems} \label{sec:conclusions}
In this article, we reported on the discovery and geometric construction of families of spatially-periodic canard solutions of the Brusselator PDE \eqref{Brusselator}. 
They are created in subcritical Turing bifurcations at $B_T=(1+\eps A)^2$,  \eqref{kTBT}, and they lie on stable and unstable branches that meet at saddle-node bifurcations of periodic patterns (Fig.~\ref{fig:bifurcation3d}).
On the upper branches, many of these patterns are observed numerically to be attractors for the PDE, both on periodic domains and on large domains.
Furthermore, the growth of new spikes/pulses and the evolution toward stable periodic patterns were observed numerically to occur through families of unstable multi-pulse canard patterns, which lie along the lower branches.

\subsection{Conclusions}
\label{ss:concl}

The spatially-periodic canards have multiscale structure.
They consist of segments along which the solution varies gradually, interspersed with one or more brief intervals on which the $u$ component has one or more steep pulses and the gradient of the $v$ component has a small jump.
Periodic canards with single pulses have been shown in Figs.~\ref{fig:representative-SAO}, \ref{fig:representative}, \ref{fig:kabovekT}, \ref{fig:kbelowkT}, \ref{fig: Steady States}, \ref{fig:sim1}, and \ref{fig:sim2}. 
The pulse amplitudes range from small, through intermediate, up to large.
The former (of $\mathcal{O}(\eps)$ amplitude or smaller) are found close to the Turing bifurcation, as shown for example in Fig.~\ref{fig:representative-SAO} and the lavender pattern in Fig.~\ref{fig:kabovekT}(a).
Solutions with pulses of intermediate ($\mathcal{O}(1)$) height lie slightly further from the Turing bifurcation, such as the cyan and blue solutions in Fig.~\ref{fig:kabovekT} and the attracting state in Fig.~\ref{fig: Steady States} (b).
There are also canard solutions that have asymptotically large ($\mathcal{O}(1/\sqrt{\eps})$ or $\mathcal{O}(1/\eps)$) amplitude pulses (and asymptotically small wavenumber, {\it i.e.,} large spatial period), recall Figs.~\ref{fig:representative}, \ref{fig:kbelowkT}, \ref{fig: Steady States}(a), (c), (d), \ref{fig:sim1}, and \ref{fig:sim2}.
Furthermore, we found families of canards with two or three pulses per period (Figs. \ref{fig:selfsimilarB1} and \ref{fig:selfsimilarB1.0175}), as well as canards with many pulses per period (Fig.~\ref{fig:bursting}). 

These families of pulsatile spatially-periodic canards are observed in a large regime in the $(B,k)$ plane, recall Fig.~\ref{fig:existenceballoon}, for a broad interval of $A$ values, and two or more canards exist at all points $(A,B,k)$ that we sampled in this regime.
Along the initial segments of the branches that emerge from the subcritical Turing bifurcation, our numerical simulations indicate that these periodic canards are unstable.
Then, for $B<B_T$, the branches turn around where they undergo a saddle-node of periodics along a fold curve $B_F(k)$, recall Fig.~\ref{fig:bifurcation3d}, in which families of (numerically) stable spatially-periodic canards are created. 
Stable patterns with $\mathcal{O}(1)$ or asymptotically large amplitude are observed along the upper branches.
Recall for example Figs.~\ref{fig: Steady States}, \ref{fig:sim1}, and \ref{fig:sim2}. 

The results of numerical solutions -both on periodic domains and on large domains with homogeneous Neumann boundary conditions- reveal rich dynamics of the spatially-periodic canards.
Simulations with constant values of the bifurcation parameter $B$ confirm that the Turing bifurcation is subcritical and that there are saddle-node bifurcations (fold curves) along which the stable and unstable periodics disappear. 
Then, simulations with slowly increasing and decreasing $B(t)$ were used to identify more of the structure of the Busse balloon of stable canard patterns. 
These simulations also revealed how unstable canard solutions appear to guide the evolution during the transient stages, leading to the development of attracting spatially-periodic canards.
New spikes/pulses are observed to originate from the folded singularities.
Moreover, on large domains, we have observed the simultaneous (to leading order) formation of new spikes leading to period halving across the entire domain, as well as the step-wise formation of new spikes on the wings of existing pulses leading to invading fronts. 
Recall Figs.~\ref{fig:sim4}(a) and (b), respectively.

In the analysis of the spatial ODE system \eqref{spatialODE}, we have demonstrated that the families of spatially-periodic canards are created by folded singularities. In particular, for $B$ near $B_T$, the canards are created by a reversible folded saddle-node singularity of type II (RFSN-II), recall \eqref{M-RFSNII}, which coincides with the Turing bifurcation point, asymptotically in the limit as $\eps \to 0$.
The canards with asymptotically large periods have extensive segments that lie along the true and faux canards of the RFSN-II point,
as shown in Figs.~\ref{fig:takeoff+touchdown} and \ref{fig:pulseconstruction}.
(Recall also Figs.~\ref{fig:saddleslowmanifold} and  \ref{fig:persistence} and the analysis in Sec.~\ref{s:K1} for more about the properties and persistence of the true and faux canards of the RFSN-II point.) 
Then, for values of the bifurcation parameter $B > B_T$, with $|B-B_T| =  \mathcal{O}(1)$ ({\it i.e.,} away from the Turing point), the system has a reversible folded saddle (RFS) point, recall \eqref{M-RFS}.
In this large regime in parameter space, the spatially-periodic canards are created by the true and faux canards of the RFS point. (See Fig.~\ref{fig:sim4}, for example.)

For all of these canards, the projections into the $(u,v)$ plane illustrate that the slow segments lie close to the left branch of the kinetics nullcline and that the pulses are excursions around the right branch. 
See Figs.~\ref{fig:kabovekT}(b), \ref{fig:kbelowkT}(b), \ref{fig: Steady State Projections}, \ref{fig:sim1}(d), \ref{fig:sim2}(b), \ref{fig:sim5}(c), (f), (i), and \ref{fig:sim3b}(c), (f), (i).
The left branch of the nullcline is a two-dimensional saddle sheet $S^0_s$ \eqref{criticalmanifold} in the singular ($\eps=0$) limit, and it persists for $0<\eps\ll 1$ as a two-dimensional saddle slow manifold $S^\eps_s$ in the four-dimensional phase space of the spatial ODE \eqref{spatialODE} (recall Figs.~\ref{fig:criticalmanifold} and \ref{fig:saddleslowmanifold}), just as is the case for the van der Pol PDE in \cite{VDK2025}.
The slow segments lie exponentially close to $S^\eps_s$, and for the periodic patterns with large-amplitude pulses, the slow segments lie exponentially close to the true and faux canards of the RFSN-II and RFS points on $S^\eps_s$.

Then, as shown in the geometric construction (recall Sec.~\ref{s:geoconstruction}), the patterns diverge from $S^\eps_s$ at points on the takeoff curve $T_{\rm takeoff}$ to leading order and make fast excursions (pulses in the $u$ component) about the right branch of the kinetics nullcline, which is a two-dimensional manifold of center points $S^0_c$.
The $v$-component is constant to leading order in $\mathcal{O}(\eps)$ during the fast jump.
There is an inverse relation between this value of $v$ and the size of the pulse in $u$: the smaller the $v$ value during the jump, the larger the pulse amplitude (recall also Fig.~\ref{fig:persistence}).
Upon completing the fast excursions, the patterns return to the neighborhood of $S^\eps_s$ near points on the touchdown curve $T_{\rm touchdown}$.
Overall, the large $\mathcal{O}(1/\eps)$ amplitude pulses start near the transverse intersection of $T_{\rm takeoff}$ with $\Gamma_{\rm faux}$ and terminate near the transverse intersection of $T_{\rm touchdown}$ with $\Gamma_{\rm true}$.
The periodic canards with $\mathcal{O}(1)$ and $\mathcal{O}(\eps)$ amplitude pulses in $u$ connect points that are higher up (in $v$) on $T_{\rm takeoff}$ and $T_{\rm touchdown}$.
(See Fig.~\ref{fig:pulseconstruction}(a) and (b), respectively.)

In addition, the projections into the $(u,q)$ plane and the $(u,p)$ plane reveal important properties of the canard patterns. 
For example, the $(u,q)$ projections show that the RFSN-II point is a cusp point (Fig.~\ref{fig:desingularized}(b)), and that the true and faux canards are the single-branched stable and unstable manifolds of the cusp.
In turn, this leads to the formation of patterns that have distinct cusp-like components midway along the slow segments, including patterns with $B=1 +\mathcal{O}(\eps)$, recall Figs.~\ref{fig:kbelowkT}(d), \ref{fig:maxcanards}, \ref{fig:selfsimilarB1}, and \ref{fig:selfsimilarB1.0175}.

Lastly, nearly self-similar dynamics was observed for a class of spatially-periodic canards in the Brusselator PDE \eqref{Brusselator}. 
In the singular limit, the system in the rescaling chart is Hamiltonian \eqref{H2}, and the key level set of this Hamiltonian is scale invariant \eqref{H2selfsimilar}. 
Then, for $0<\eps\ll 1$, even though the system is no longer Hamiltonian, the solutions can exhibit nearly self-similar dynamics.
For $B=1$ and for $B$ close to one,
the organizing center about which the nearly self-similar oscillations take place is the RFSN-II point \eqref{M-RFSNII}.
Then, for $B > B_T$, with $|B-B_T| = \mathcal{O}(1)$, the nearly self-similar oscillations are about the equilibrium point $E$.
(Recall Figs.~\ref{fig:selfsimilarB1} and \ref{fig:selfsimilarB1.0175}, respectively.)

Overall, the Brusselator PDE \eqref{Brusselator} constitutes a useful second example of ``les canards de Turing". 
The results here show that Turing's ducks occur in systems with monostable reaction kinetics, and not just in systems with bistable kinetics as in the van der Pol PDE in \cite{VDK2025}.
Also, the Brusselator results show that it is not necessary for the full fourth-order spatial ODE system to have a conserved quantity, as is the case for the van der Pol system.
Furthermore, in the Brusselator model, 
the Hopf and Turing modes are separated by an O(1) amount, which facilitates studying the spatial canards separately, whereas they are asymptotically close in the van der Pol PDE.
Therefore, the results here further support the conjecture, made in Section 12.2 of \cite{VDK2025}, that these families of spatially-periodic canards exist in broad classes of reaction-diffusion systems.

\subsection{Discussion} 
\label{s:discussion}

In this subsection, we discuss the similarities and differences between the spatially-periodic canards here in the Brusselator PDE and the classical temporal limit cycle canards \cite{BCDD1981,D1984,E1983} known to exist in fast-slow systems of ordinary differential equations, such as the classic van der Pol ODE in the limit of relaxation oscillations. 
In particular, as was first reported in \cite{VDK2025} for ``les canards de Turing" ({\it i.e.}, for the spatially-periodic canards) in the van der Pol PDE, these canard solutions are the spatial analogs of the temporal canards.

Both types of systems have Hopf bifurcations: the spatial ODEs have reversible 1:1 resonant singular Hopf bifurcations (aka singular Turing bifurcations) (see \cite{HI2011,IP1993} in general, and \eqref{kTBT} and Fig.~\ref{f-quartet} here), and the fast-slow systems undergo singular Hopf bifurcations in time.
Immediately beyond the relevant Hopf bifurcations, the systems exhibit small-scale oscillations: about the homogeneous state here just beyond the singular Turing bifurcation in the case of spatial canards, and about the local extremum just beyond the singular temporal Hopf in the case of temporal canards.
There are critical values of the bifurcation parameter, which are asymptotically close to the Hopf points, at which the systems exhibit a canard explosion.
Furthermore, both types of canards are created by folded singularities in the corresponding ODEs: the spatial canards by RFSN-II and RFS points (\eqref{M-RFSNII} and \eqref{M-RFS}), and the temporal canards by FSN-II and FS points.
Lastly, the canards of these folded singularities are both represented by key algebraic solutions of the ODEs in the rescaling charts 
(recall Lemma~\ref{lem-gamma0} here and the solution $\gamma_{c,2}$ in Chart $K_2$ in Section 3.3 of \cite{KS2001} for classical temporal canards in planar fast-slow systems, for example).

There are also some differences that stem from the structures of the underlying spatial and temporal ODE systems. 
The spatial ODE system \eqref{spatialODE} is four-dimensional and has a reversibility symmetry \eqref{reversibility}. 
Hence, the equilibria of the layer system of the spatial ODE can only be saddles, centers, and saddle-node points.
In contrast, the equilibria of the layer systems in fast-slow ODEs with temporal canards are attracting or repelling.
This translates directly to the main difference between the critical manifolds, which here are two-dimensional saddle and normally elliptic manifolds, $S_s^\eps$ and $S_c^\eps$, in the four-dimensional phase space, whereas they are typically one- or two-dimensional, normally attracting or repelling manifolds in the case of temporal canards.
Furthermore, while the drift here along the saddle slow manifolds for spatial canards is similar to the drift along the normally attracting manifolds for temporal canards, the spatial canards oscillate around the elliptic manifold 
(recall Figs.~\ref{fig:representative-SAO}, \ref{fig:representative}, \ref{fig:kabovekT}, \ref{fig:kbelowkT}, \ref{fig: Steady State Projections}, \ref{fig:sim1}, \ref{fig:sim2}, \ref{fig:sim5}, and \ref{fig:sim3b}), rather than drift slowly along normally repelling manifolds, as the temporal canards do.
Lastly, the maximal spatial canards here occur when the persistent true and faux canards of an RFSN-II or RFS point coincide to all orders (recall Fig.~\ref{fig:persistence}), whereas in the case of temporal ODEs the maximal canards occur when the branches of the normally attracting and repelling manifolds coincide to all orders.

\subsection{Some open problems}
\label{s:openproblems}

In this subsection, we briefly mention some open problems. 
One natural next step is to validate the observations made in the numerical simulations of the PDE and to establish the stability (or instability -- depending on the parameters) of canard patterns with single pulses.
This will help to further understand the impacts of the canard patterns on the dynamics of the Brusselator model (\ref{Brusselator}).
Such an analysis is expected to be possible, given the methods developed in the literature by which the stability can be determined of pulses that pass along non-normally hyperbolic parts of a slow manifold (in the existence problem associated to a system of singularly perturbed reaction-diffusion equations), see for example \cite{BJSW2008, CdRS2016, HDK2013}. 
Although these methods have been developed for homoclinic patterns, they may be extended to spatially periodic patterns, for instance along the lines developed in \cite{dR2018, dRDR2016}. 
These results will form the foundation of the subsequent steps in which the (in)stability and marginal stability of the multi-pulse spatially-periodic patterns is determined -- with a focus the multi-pulse canard patterns with tiny spikes and more generally those on the lower parts of the solution branches, which emerge from the Turing bifurcations and which end at the saddle-node bifurcations. 

These steps, in turn, prepare for the challenging and relevant investigations of the guiding role of these multi-pulse canards in the selection process that takes place during transitions between pulsatile canard patterns. 
As discussed in Secs.~\ref{s:PDE-largedomains} and \ref{ss:concl}, two distinct mechanisms seem to be playing a role here. 
The first mechanism consists of the fast and nearly-simultaneous growth or disappearance of tiny spikes in between pairs of adjacent spikes on intervals filled with periodic patterns, as shown for example in Figs. \ref{fig:sim2}(c)-(h), \ref{fig:sim4}(a), (c) and \ref{fig:sim5}.
The second mechanism, which is illustrated by the patterns in Figs.~\ref{fig:sim4}(b), (d) and \ref{fig:sim4b}, consists of new spikes forming on the wings of existing pulses, and it is reminiscent of (a stationary version of) the process of trace firing in which new pulses subsequently appear at the trailing edges of traveling wave trains  \cite{R2005} (since the present patterns are stationary, the new pulses appear here on the wings on both sides of existing pulses). 
To obtain insight in these processes, novel methods, for instance building on the ideas presented in \cite{CS2018}, need to be developed that also incorporate the interactions between spatial canards, their stability, and their temporal dynamics.

Similar novel methods will need to be developed to obtain an analytical understanding of the numerical simulations of the Brusselator PDE \eqref{Brusselator} that show rich pulse interaction dynamics during the evolution to some of the stable spatially-periodic steady state patterns. 
For example, in simulations with $A=1$, $B=1$, and $\eps=0.01$, we started from initial data with $N$-pulses ($N=2,3,4$) relatively close to each other and observed the pulses slowly drift apart.
These pulses evolve until they reach a separation distance that makes the solution periodic on the domain, {\it i.e.,} with a period that is a factor of $1/N$ times the domain length, with a wavenumber inside the Busse balloon. 
There is a well-developed theory to describe these kinds of interactions for weakly interacting pulses \cite{EMN2002, P2002} and for semi-strong interactions in singularly perturbed systems \cite{DKP2007}.
However, these methods cannot be applied directly to canard patterns (due to their (partly) non-normally hyperbolic character).

Inclusion of advection in reaction-diffusion systems with singular Turing bifurcations (or some other terms with odd-order spatial derivatives) would also be of interest, since such terms break the reversibility symmetry and since such terms arise in applications (see \cite{BCD2019,C2024,Harley2014a,Harley2014b,WechselbergerPettet2010}). 
In such systems, the fixed points of the layer (fast) problems can be nodes and spirals, rather than just saddles, centers, and saddle-node points, as imposed by the reversibility symmetry \eqref{reversibility} here.
Then, as a result, the persistent slow manifolds can be of more general stability types, and also the folded singularities on the fold curves can be more general, including folded nodes and folded spirals. Finally, in turn, the possible canard dynamics induced by these folded singularities on the patterns can be even more complex. 
Examples of canard dynamics in the traveling wave solutions of advection-reaction-diffusion systems can be found in \cite{Harley2014a,Harley2014b,WechselbergerPettet2010}. 
The relation of these solutions to the singular Turing bifurcations studied here would be of interest. 

In current work, we are addressing the conjecture \cite{VDK2025} to show that there is a broad class of coupled reaction-diffusion equations of activator-inhibitor type that exhibit singular Turing bifurcations, 
RFSN-II and RFS points, and families of spatially-periodic canard patterns. 
In this direction, it is also of interest to determine the broad conditions under which PDEs with singular Turing bifurcations have an RFSN-II point asymptotically close to the Turing bifurcation, whether the Turing bifurcations must be subcritical or can also be supercritical, and which of the sufficient conditions are necessary.

In \cite{VDK2025} and this article, we leverage the fact that the elliptic problem for the steady state equation can be rewritten as a spatial dynamics problem, since we work in 1-D.
It would be of interest to study these PDEs on two-dimensional domains, and to determine if there is a geometric property of the elliptic problem that generalizes the spatial canard phenomenon to 2-D.

The addition of pulses or spikes to solutions also occurs in neuroscience, where the phenomenon is often labeled as spike-adding.
See for example \cite{BIPS2021,C2020,CS2018,DKK2013,DKS2016,DK2018}.
In these systems, the (kinetics) nullclines are typically bistable, and temporal canards --especially the true and faux canards of folded saddles-- play central roles in spike-adding.
It would be of interest, on the one hand, to explore what additional information one can obtain for (spatial) pulse formation in the Brusselator and other PDEs based on the known results about spike-adding in temporal ODEs and, on the other hand, what effects inclusion of diffusion (which naturally occurs in ion channel models) has on spike-adding.


Overall, this work is part of the growing larger field of spatial canards and spatio-temporal canards in PDEs. 
We refer the reader to the Introduction of \cite{VDK2025}, the references cited therein,
and the open problems listed there.

\bigskip
\noindent
{\bf Acknowledgments}
The authors thank Margaret Beck, Ryan Goh, Gene Wayne for useful conversations while serving on the thesis committee of R.J., and Hans Kaper for a useful conversation.
The research of A.D. is supported by the ERC-Synergy project RESILIENCE (101071417).
The work of R.J. was supported in part by a research account sponsored by Boston University.

\appendix

\section{The Proof of Proposition 2.1.}

We begin by translating the reversible Hopf bifurcation point of \eqref{spatialODE} to the origin by introducing the new variables
\[ u_1 = u - A, \quad u_2 = p, \quad u_3= v - \tfrac{B}{A}, \quad u_4 = q, \quad {\rm and} \quad \mu =  B - B_T, \]
where $B_T = \left( 1+\eps A \right)^2$. 
Hence, the spatial ODE system \eqref{spatialODE} transforms to
\begin{equation} \label{eq:HI_TaylorSeries}
    \tfrac{d \boldsymbol{u}}{dx} = \mathbb L \boldsymbol{u} + R_{20}(\boldsymbol{u},\boldsymbol{u}) + R_{30}(\boldsymbol{u},\boldsymbol{u},\boldsymbol{u}) + \mu R_{11}(\boldsymbol{u})+\mu R_{21}(\boldsymbol{u},\boldsymbol{u}),
\end{equation}
where the linear operator $\mathbb L$ given by
\[ \mathbb L = \begin{bmatrix} 0 & 1 & 0 & 0 \\ 1-B_T & 0 & -A^2 & 0 \\ 0 & 0 & 0 & \eps \\ \eps B_T & 0 & \eps A^2 & 0 \end{bmatrix}, \]
corresponds to the linearization of \eqref{spatialODE} at the reversible Hopf bifurcation and has eigenvalues $\pm i \omega$ (each of multiplicity 2), where $\omega = \sqrt{\eps A}$ is the critical frequency. The multi-linear functions $R_{ij}$ denote nonlinear terms of order $\boldsymbol{u}^i \mu^j$ where $\boldsymbol{u}^i$ denotes a monomial of the form $u_1^{i_1} u_2^{i_2} u_3^{i_3} u_4^{i_4}$ such that $i_1+i_2+i_3+i_4 = i$. The parameter-independent nonlinear terms are given by 
\begin{equation}
\begin{split}
    R_{20}(\boldsymbol{u},\boldsymbol{v}) &= \begin{bmatrix} 0 \\ -\tfrac{B_T}{A} u_1 v_1 -2A u_1 v_3 \\ 0 \\ \eps \left( \tfrac{B_T}{A} u_1 v_1 + 2A u_1 v_3 \right)   \end{bmatrix}, 
    \quad
    R_{30}(\boldsymbol{u},\boldsymbol{v},\boldsymbol{w}) = \begin{bmatrix} 0 \\ -u_1 v_1 w_3 \\ 0 \\ \eps u_1 v_1 w_3 \end{bmatrix}, 
\end{split}    
\end{equation}
and the parameter-dependent nonlinearities are given by 
\begin{equation}
\begin{split}
    R_{11}(\boldsymbol{u}) &= \begin{bmatrix} 0 \\ -u_1 \\ 0 \\ \eps u_1 \end{bmatrix}, 
    \qquad \text{ and } \qquad 
    R_{12}(\boldsymbol{u},\boldsymbol{v}) = \begin{bmatrix} 0 \\ -\tfrac{1}{A} u_1 v_1 \\ 0 \\ \eps \tfrac{1}{A} u_1 v_1 \end{bmatrix}.
\end{split}    
\end{equation}

Recall from \eqref{reversibility} that the system \eqref{eq:HI_TaylorSeries} has the reversibility symmetry 
\begin{equation}\label{eq:app_reversibility}
    \mathcal{R} {\bf F} ({\bf u}) = - {\bf F}( \mathcal{R}({\bf u})),
\end{equation}
where $\boldsymbol{F}$ denotes the vector field in \eqref{eq:HI_TaylorSeries} and 
\[ \mathcal R = \begin{bmatrix} 1 & 0 & 0 & 0 \\ 0 & -1 & 0 & 0 \\ 0 & 0 & 1 & 0 \\ 0 & 0 & 0 & -1 \end{bmatrix} \]
reflects the fact that the Brusselator PDE has no odd-order spatial derivatives. 

Next, we establish that the system \eqref{eq:HI_TaylorSeries} satisfies the criteria for the normal form theory of \cite{HI2011}. To do this, we examine the eigenvector, $\zeta_0$, and generalized eigenvector, $\zeta_1$, of $\mathbb L$ corresponding to the eigenvalue $i \omega$. Without loss of generality, we fix 
\[ \zeta_0 = \begin{bmatrix} A \\ i A \omega \\ -\eps (1+\omega^2) \\ - i \omega (1+\omega^2)  \end{bmatrix} \quad \text{ and } \quad \zeta_1 = \begin{bmatrix} 0 \\ A \\ -2i\, \tfrac{1}{A}\, \omega \\ 1-\omega^2 \end{bmatrix}. \]
Then, it can be shown (by direct computation) that 
\[ \mathbb L \zeta_0 = i \omega \zeta_0, \quad (\mathbb L - i \omega) \zeta_1 = \zeta_0, \quad \mathcal{R} \zeta_0 = \overline{\zeta_0}, \quad \text{ and } \quad \mathcal{R} \zeta_1 = - \overline{\zeta_1}. \]
Hence, the system \eqref{eq:HI_TaylorSeries} satisfies the criteria for the $(i\omega)^2$-bifurcation from Chapter 4.3.3 of \cite{HI2011}. 
That is, there exists a change of variables of the form 
\begin{equation} \label{eq:app_trans}
\boldsymbol{u} = U \zeta_0 + V \zeta_1 + \overline{U} \overline{\zeta_0}+\overline{V} \overline{\zeta_1} + \Phi(U,V,\overline{U},\overline{V},\mu), 
\end{equation}
such that the system \eqref{eq:HI_TaylorSeries} is transformed into the normal form
\begin{equation} \label{eq:HI_NormalForm}
\begin{split}
    \tfrac{dU}{dx} &= i \omega U + V + i U P + \rho_U \\ 
    \tfrac{dV}{dx} &= i \omega V + i V P + U Q + \rho_V
\end{split}
\end{equation}
where $P$ and $Q$ are real-valued functions given by 
\begin{equation*}
\begin{split}
    P &= \alpha \mu + \beta U \overline{U} + \tfrac{1}{2} i \gamma \left( U \overline{V} - \overline{U} V \right) \\ 
    Q &= a \mu + b U \overline{U} + \tfrac{1}{2} i c \left( U \overline{V} - \overline{U} V \right),
\end{split}
\end{equation*}
and $\rho_U$ and $\rho_V$ encode the higher order terms. 

To compute the normal form coefficients, we seek a vector $\zeta_1^*$ orthogonal to the range of $(\mathbb L - i \omega)$ such that 
\[ 
\langle \zeta_0, \zeta_1^* \rangle = 0, \quad 
\langle \overline{\zeta_0}, \zeta_1^* \rangle = 0, \quad 
\langle \zeta_1, \zeta_1^* \rangle = 1, \quad \text{ and } \quad 
\langle \overline{\zeta_1}, \zeta_1^* \rangle = 0.
\]
With the choices of $\zeta_0$ and $\zeta_1$ given above, we find 
\[ \zeta_1^* = \frac{1}{4} \begin{bmatrix} i \eps \tfrac{1+\omega^2}{\omega} \\ \tfrac{1+\omega^2}{A} \\ i \tfrac{A}{\omega} \\ 1 \end{bmatrix}. \]
Then, by differentiating \eqref{eq:app_trans} with respect to $x$ and using \eqref{eq:HI_TaylorSeries} and \eqref{eq:HI_NormalForm}, we obtain the following invariance equation
\begin{equation}  \label{eq:app_invariance}
\begin{split}
    \mathbb L \Phi + R &= \left( i\omega U + V \right) \tfrac{\partial \Phi}{\partial U} + i \omega V \tfrac{\partial \Phi}{\partial V} + \left( -i \omega \overline{U} + \overline{V} \right) \tfrac{\partial \Phi}{\partial \overline{U}} - i \omega \overline{V} \tfrac{\partial \Phi}{\partial \overline{V}} \\ 
    & \hspace{2em} + \left( i U \left( \zeta_0 + \tfrac{\partial \Phi}{\partial U} \right) - i \overline{U} \left( \overline{\zeta_0} + \tfrac{\partial \Phi}{\partial \overline{U}} \right) \right) P + (i V P + U Q) \left( \zeta_1 + \tfrac{\partial \Phi}{\partial V} \right) \\
    & \hspace{4em} + \left( -i \overline{V} P + \overline{U} Q \right) \left( \overline{\zeta_1}+\tfrac{\partial \Phi}{\partial \overline{V}} \right) + \cdots
\end{split}
\end{equation}
where $R$ denotes the nonlinear terms of \eqref{eq:HI_TaylorSeries}. Substituting the Taylor series expansion
\[ \Phi = \displaystyle \sum\limits_{2\leq r+s+q+\ell+m \leq p} U^r V^s \overline{U}^q \overline{V}^\ell \mu^m \Phi_{rsq\ell m} \]
into the invariance equation \eqref{eq:app_invariance} and equating coefficients of like powers of $U, V, \overline{U}, \overline{V}, \mu$, and their products, we obtain a set of linear equations for the normal form coefficients. We refer to App. D.2 of \cite{HI2011} for full details. 

Finally, we use the orthogonality of $\zeta_1^*$ to simplify the calculation of the normal form coefficients, where we again refer the reader to App. D.2 of \cite{HI2011} for explicit details. We find that the normal form coefficients are given by 
\begin{equation} \label{eq:app_normalformcoefficients}
\begin{split}
    a &= -\tfrac{1}{4} \\ 
    b &= -\tfrac{1}{36\omega^2} \left( 8-30\omega^2-43 \omega^4+3\omega^6+8\omega^8 \right) \\
    c &= \tfrac{1}{216\omega^3} \left( 104-282 \omega^2+233 \omega^4+201 \omega^6+104 \omega^8 \right) \\
    \alpha &= \tfrac{1}{8\omega} \\
    \beta &= -\tfrac{1}{432\omega^3} \left( 8+270\omega^2+101\omega^4-27\omega^6+8\omega^8 \right) \\
    \gamma &= -\tfrac{1}{432 \omega^4} \left( 192+2\omega^2+153\omega^4+79\omega^6+192\omega^8 \right),
\end{split}
\end{equation}
where we recall that $\omega = \sqrt{\eps A}>0$. From these normal form coefficients, we can conclude that 
\[ a<0, \quad c>0, \quad \alpha >0, \quad \beta<0, \quad \text{ and } \quad \gamma<0. \]
The remaining coefficient $b$ encodes the criticality of the reversible Hopf bifurcation, and it can be expressed in terms of the original Brusselator parameters as
\[ b = -\frac{2}{9A}\eps^3 \left( A+\frac{2}{\eps} \right)\left( A+\frac{1}{\eps} \right) \left( A - \frac{21+\sqrt{313}}{16\eps} \right) \left( A-\frac{21-\sqrt{313}}{16\eps} \right). \]
Hence, we conclude that 
\[ b \,\,
\begin{cases} 
\,\, <0 \quad {\rm (subcritical)} & \text{ for } A \in \left( 0,\frac{21-\sqrt{313}}{16\eps} \right) \cup \left( \frac{21+\sqrt{313}}{16\eps}, \infty \right) \\ 
\,\, >0 \quad {\rm (supercritical)} & \text{ for } A \in \left( \frac{21-\sqrt{313}}{16\eps},\frac{21+\sqrt{313}}{16\eps} \right).
\end{cases} \]
This establishes the results of the weakly nonlinear analysis in Sec.~\ref{ss:weaklynonlinear}.

\section{Numerical Computation of \texorpdfstring{$W^u\!\left( S_s^{\eps} \right) \cap W^s\!\left( S_s^{\eps} \right)$}{Lg}} \label{app:homman}

In this appendix, we describe the numerical continuation method that was used to generate the homoclinic orbits in the transverse intersection $W^u(S^\eps_s) \bigcap W^s(S^\eps_s)$, shown in Fig.~\ref{fig:homoclinicmanifold}.
We consider the spatial ODEs \eqref{spatialODE-S0s-rectified} in which the saddle sheet, $S_s^0$, of the critical manifold has been rectified to the $(v,q)$-plane:
\begin{equation}  \label{eq:homoclinicmanifold}
  \begin{split}
      \hat{u}_x &= p - \eps q \tfrac{u_s^2(v)}{v \left( u_c(v) - u_s(v) \right)} \\ 
      p_x &= v \hat{u} \left( u_c(v) - u_s(v) - \hat{u} \right) \\ 
      v_x &= \eps \, q \\ 
      q_x &= \eps \left( -v \hat{u} \left( u_c(v) - u_s(v) - \hat{u} \right) + \hat{u} +u_s(v)-v\, u_s(v) u_c(v) \right),
  \end{split}
\end{equation}
where $u_c(v) = \tfrac{1+B+\sqrt{(1+B)^2-4A v}}{2v}$ and $u_s(v) = \tfrac{1+B-\sqrt{(1+B)^2-4A v}}{2v}$. For this system, the saddle slow manifold has the following asymptotic expansion
\[ S_s^{\eps} = \left\{ \begin{bmatrix} \hat{u}_{\eps}(v,q) \\ p_{\eps}(v,q) \end{bmatrix} = \begin{bmatrix} 0 \\ 0 \end{bmatrix} + \eps \begin{bmatrix} 0 \\ \tfrac{q\,u_s^2(v)}{v \left( u_c(v) - u_s(v) \right)} \end{bmatrix} + \mathcal{O}(\eps^2) \right\}. \]
The stable spectrum and stable subspace of the saddle slow manifold perturb to
\[ \sigma_s^{\eps} = \left\{ \lambda_s(v) + \mathcal{O}(\eps^2) \right\}, \qquad \mathbb{E}^s\left( S_s^{\eps} \right) = \operatorname{span} \left\{ \begin{bmatrix} 1 \\ \lambda_s(v) \\ 0 \\ 0 \end{bmatrix} + \eps \begin{bmatrix}  1 \\ \lambda_s(v) \\ 0 \\ \tfrac{1-\lambda_s^2(v)}{\lambda_s(v)} \end{bmatrix} + \mathcal{O}(\eps^2) \right\},  \]
where $\lambda_s(v) = - \sqrt{v(u_c(v)-u_s(v))}$ is the stable eigenvalue of the layer problem evaluated along the saddle sheet, $S_s^0$, of the critical manifold. 
Similarly, the unstable spectrum and unstable subspace of the saddle slow manifold perturb to 
\[ \sigma_u^{\eps} = \left\{ \lambda_u(v) + \mathcal{O}(\eps^2) \right\}, \qquad \mathbb{E}^u\left( S_s^{\eps} \right) = \operatorname{span} \left\{ \begin{bmatrix} 1 \\ \lambda_u(v) \\ 0 \\ 0 \end{bmatrix} + \eps \begin{bmatrix}  1 \\ \lambda_u(v) \\ 0 \\ \tfrac{1-\lambda_u^2(v)}{\lambda_u(v)} \end{bmatrix} + \mathcal{O}(\eps^2) \right\},  \]
where $\lambda_u(v) = \sqrt{v(u_c(v)-u_s(v))}$ is the unstable eigenvalue of the layer problem evaluated along the saddle sheet, $S_s^0$, of the critical manifold. 

To compute the manifold $W^u\!\left( S_s^{\eps} \right) \cap W^s\!\left( S_s^{\eps} \right)$, we seek solutions of the two-point boundary value problem defined by 
\[ \boldsymbol{u}_x = T \boldsymbol{F}(\boldsymbol{u}), \]
where $\boldsymbol{u} = (\hat{u},p,v,q)$, $T$ is the length of the spatial interval, and $\boldsymbol{F}$ is the vector field in \eqref{eq:homoclinicmanifold}, subject to the boundary conditions 
\begin{equation}  \label{eq:homman_bc01}
    \begin{split}
        \boldsymbol{u}(0) &= \begin{bmatrix} \hat{u}_{\eps}(v_0,q_0) \\ p_{\eps}(v_0,q_0) \\ v_0 \\ q_0 \end{bmatrix} + \delta \boldsymbol{v}_u, \quad {\rm and } \quad \boldsymbol{u}(1) = \begin{bmatrix} \hat{u}_{\eps}(v_1,q_1) \\ p_{\eps}(v_1,q_1) \\ v_1 \\ q_1 \end{bmatrix} + \delta \boldsymbol{v}_s.
    \end{split}
\end{equation}
Here, $\boldsymbol{v}_u \in \mathbb{E}^u\left( S_s^{\eps} \right)$ and $\boldsymbol{v}_s \in \mathbb{E}^s\left( S_s^{\eps} \right)$ are unit vectors aligned with the unstable and stable subspaces, respectively. The parameter $\delta \ll 1$ (independent of $\eps$) measures the distance of $\boldsymbol{u}(0)$ and $\boldsymbol{u}(1)$ from $S_s^{\eps}$ along $\boldsymbol{v}_u$ and $\boldsymbol{v}_s$, respectively. 
The free parameters $v_0$ and $q_0$ specify the base point on $S_s^{\eps}$ from which the unstable manifold emerges. Similarly, the free parameters $v_1$ and $q_1$ specify the base point on $S_s^{\eps}$ at which the stable manifold terminates. To guarantee that the computed solutions are exponentially close to the saddle slow manifold, $S_s^{\eps}$, we impose the additional boundary condition
\begin{equation} \label{eq:homman_bcT}
    T = \frac{1}{\sqrt{\eps} \lambda_u(v_0)} - \frac{1}{\sqrt{\eps} \lambda_s(v_1)}.
\end{equation} 

For a fixed $\delta \ll 1$ and $v_0>0$, a starting solution of \eqref{eq:homoclinicmanifold} subject to \eqref{eq:homman_bc01} is 
\[ \boldsymbol{u} = \left( \hat{u}_{\rm HOM}(\xi;v_0), p_{\rm HOM}(\xi;0), v_0, 0 \right), \quad \text{ for } \quad \eps = 0. \]
Numerical continuation with respect to $(T,v_0,q_0,v_1,q_1,\eps)$ can be used to homotopically grow this orbit segment out to a desired $0< \hat{\eps} \ll 1$. Next, the end solution of the previous continuation run can be continued with respect to $(T,v_0,q_0,v_1,q_1,\delta)$ to adjust the integration time so that \eqref{eq:homman_bcT} is satisfied. Finally, the end solution of the previous continuation run is used as a starting solution for the full two-point boundary value problem \eqref{eq:homoclinicmanifold} subject to \eqref{eq:homman_bc01} and \eqref{eq:homman_bcT}. Continuation with respect to $(T,v_0,q_0,v_1,q_1,\delta)$ allows the manifold $W^u\!\left( S_s^{\eps} \right) \cap W^s\!\left( S_s^{\eps} \right)$ to be swept out.

\section{The singular true and faux canards of \texorpdfstring{$M_{\rm RFS}$ \eqref{M-RFS}}{Lg}}\label{sec:App-RFS}
We derive the explicit formulas for the singular true and faux canards of the reversible folded saddle point, $M_{\rm RFS}$, recall \eqref{M-RFS}.
First, in system \eqref{spatialODE}, we translate $M_{\rm RFS}$ to the origin, recalling that it exists for all $B>1$.
Let $u=\tfrac{2A}{B+1} + U$, $p=P$, $v = \tfrac{(B+1)^2}{4A} +V$, and $q=Q$.
The ODE system is
\begin{equation*}
\begin{split}
U_y &= P, \\
P_y &=\tfrac{-(B+1)^2}{4A} U^2 - \tfrac{4A^2}{(B+1)^2} V -\tfrac{4A}{B+1} UV -U^2V, \\
V_y &= \eps Q, \\
Q_y &= \eps \left(
\tfrac{A(1-B)}{B+1} + U + \tfrac{(B+1)^2}{4A} U^2 + \tfrac{4A^2}{(B+1)^2} V + \tfrac{4A}{B+1} UV + U^2 V
\right), \\
\eps_y &= 0.
\end{split}
\end{equation*}
For each value of $B > 1$, the origin in the five-dimensional $(U,P,V,Q,\eps)$ system is a degenerate equilibrium.
We apply the method of geometric desingularization 
with $U= \bar{r}^2 \bar{U}$, $P = \bar{r}^3 \bar{P}$, $V= \bar{r}^4 \bar{V}$,
$Q=\bar{r}^2 \bar{Q}$, and $\eps=\bar{r}^3 \bar{\eps}$.
Then, in the rescaling chart $K_2$, the coordinate change is 
$U= r_2^2 U_2$, $P = r_2^3 P_2$, $V= r_2^4 V_2$,
$Q=r_2^2 Q_2$, and $\eps=r_2^3$ ({\it i.e.,} $\eps_2=1$).
Also, it is useful to rescale the independent variable to $y_2 = r_2 y$.
Then, the governing equations are 
\begin{equation*}
\begin{split}
\dot{U}_2 &= P_2, \\
\dot{P}_2 &=\tfrac{-(B+1)^2}{4A} U_2^2 - \tfrac{4A^2}{(B+1)^2} V_2 -\tfrac{4A}{B+1} r_2^2 U_2V_2 -  r_2^4U_2^2V_2, \\
\dot{V}_2 &= Q_2, \\
\dot{Q}_2 &= 
\tfrac{A(1-B)}{B+1} + r_2^2 U_2 
+ r_2^4 \left( \tfrac{(B+1)^2}{4A} U_2^2 + \tfrac{4A^2}{(B+1)^2} V_2\right) + \tfrac{4A}{B+1} r_2^6 U_2V_2 + r_2^8 U_2^2 V_2, \\
\dot{r}_2 &= 0,
\end{split}
\end{equation*}
where the overdot denotes the derivative with respect to $y_2 \in \mathbb{R}$.

On the invariant set $\{ r_2 = 0 \}$, the system has two algebraic solutions,
\begin{equation} \label{truefauxcanards-RFS}
\Gamma_0 (y_2)_\pm = 
\left(
\pm \tfrac{2\sqrt{2} A^2 \sqrt{B-1}}{(B+1)^{5/2}} y_2, \, \,
\pm \tfrac{2\sqrt{2} A^2 \sqrt{B-1}}{(B+1)^{5/2}}, \, \, 
\tfrac{-A(B-1)}{2(B+1)} y_2^2, \, \, 
\tfrac{-A(B-1)}{B+1} y_2
\right).
\end{equation}
These represent the singular true and faux canards of $M_{\rm RFS}$.
They persist for $0<\eps \ll 1$, as may be shown using analysis similar to that presented in Subsec.~\ref{subsec:persistence}.
Persistent true and faux canards of $M_{\rm RFS}$ are shown in Fig.~\ref{fig:persistence}.
For $B > B_T$, with $|B-B_T| = \mathcal{O}(1)$, the spatially-periodic patterns have slow segments along these persistent true and faux canards, {\it i.e.,} $M_{\rm RFS}$ and its canards are the mechanism that create the spatially-periodic canards in this large regime.
Moreover, by deriving the perturbation expansion for these canards to a sufficiently high order in $0<r_2 \ll 1$, one can show that the true and faux canards of $M_{\rm RFS}$ limit on the true and faux canards of $M_{\rm RFSN-II}$, recall Lemma~\ref{lem-gamma0}.

\printbibliography
\end{document}